\documentclass[11pt]{article}

\usepackage{float}
\usepackage[utf8]{inputenc}
\usepackage[T1]{fontenc}
\usepackage{lmodern}
\usepackage[margin=1.1in]{geometry}
\usepackage{amsmath,amssymb,amsthm,mathtools}
\usepackage{enumitem}
\usepackage{graphicx}
\usepackage{booktabs}
\usepackage{caption}
\usepackage[round,authoryear]{natbib}
\usepackage[hidelinks]{hyperref}
\hypersetup{pdftitle = {Confidence regions for a mean from sign-flipped residuals}, pdfauthor = {Roman Zinovev}}

\numberwithin{equation}{section}

\newtheoremstyle{ims}{8pt}{8pt}{\itshape}{\parindent}{\scshape}{.}{0.5em}{\thmname{#1}\thmnumber{ #2}\thmnote{ (\upshape #3)}}
\newtheoremstyle{imsup}{8pt}{8pt}{\normalfont}{\parindent}{\scshape}{.}{0.5em}{\thmname{#1}\thmnumber{ #2}\thmnote{ (#3)}}

\theoremstyle{ims}
\newtheorem{theorem}{Theorem}[section]
\newtheorem{proposition}[theorem]{Proposition}
\newtheorem{lemma}[theorem]{Lemma}

\theoremstyle{imsup}

\newtheorem{assumption}{Assumption}

\makeatletter
\renewenvironment{proof}[1][\proofname]{\par\pushQED{\qed}\normalfont\topsep6\p@\@plus6\p@\relax\trivlist\if\relax\detokenize{#1}\relax\item[]\else\item[\hskip\labelsep\scshape#1.]\mbox{}\par\nobreak\noindent\fi\ignorespaces}{\popQED\endtrivlist\@endpefalse}
\renewcommand\section{\@startsection{section}{1}{0pt}{3.5ex plus 1ex}{2ex}{\centering\normalfont\scshape}}
\renewcommand\subsection{\@startsection{subsection}{2}{0pt}{3ex plus 1ex}{1.5ex}{\normalfont\itshape}}
\makeatother

\newcommand{\E}{\mathbb{E}}
\newcommand{\Prob}{\mathbb{P}}
\newcommand{\R}{\mathbb{R}}
\newcommand{\eqd}{\stackrel{d}{=}}
\newcommand{\bareps}{\bar\varepsilon}

\begin{document}

\begin{center}
{\large\bfseries CONFIDENCE REGIONS FOR A MEAN FROM\\[2pt] SIGN-FLIPPED RESIDUALS}\\[12pt]
{\scshape By Roman Zinovev}\\[4pt]
{\itshape Nanyang Technological University, 50 Nanyang Avenue, Singapore 639798}\\[2pt]
{\small\texttt{roman003@e.ntu.edu.sg}}
\end{center}

\begin{abstract}

Consider a problem in which confidence regions are constructed for the mean of $n$ random vectors without any assumptions about the dependence among the coordinates. Arlot, Blanchard, and Roquain construct these regions via sign reversal and prove that, when a remainder term is added, the resampled quantile keeps the confidence level, while leaving open the question of whether this remainder can be eliminated. This paper proves that, in the general case, this is not possible, and instead we propose a threshold defined by a fixed-point equation. This threshold is valid in finite samples for any seminorm, provided only that the data are sign-symmetric, and under moment conditions it is asymptotically equivalent to the ideal threshold, including when the dimension grows with $n$.

\end{abstract}



\section{Introduction}\label{sec:intro}

Let $Y_1,\dots,Y_n$ be random vectors in $\R^K$ and let $\mu$ denote the common mean of the $Y$. The main goal is to construct a nonasymptotic confidence region in the form of

\[ \mathcal C(t)=\{x\in\R^K:\phi(\bar Y-x)\le t\} \]

with $\Prob(\mu\in\mathcal C(t))\ge1-\alpha$, where $\phi$ is a seminorm, which is fixed in advance, and $t$ is the threshold computed from the data. For example, if $\phi=\|\cdot\|_\infty$, then this region yields simultaneous confidence intervals for all coordinates. Just as in \citet{ABR2010a}, our main motivation is related to neuroimaging, where, in typical cases, $n\le100\ll K$, and the structure of the dependencies can be completely arbitrary.

Since even for $K=1$, an interval that covers $\mu$ with probability at least $1-\alpha$ for every law with a finite mean has an expected length of infinity \citep[Cor.~2]{BahadurSavage1956}, \citet{ABR2010a} had to assume that the observed random variables are i.i.d.\ and symmetric about the mean $\mu$. Under this assumption, the sign-flipped means $\bar y(\varepsilon)=n^{-1}\sum_t\varepsilon_ty_t$ of the centered observations $y_t=Y_t-\mu$ provide an exact test \citep[Lemma~3.1]{ABR2010a}, but they depend on the unknown $\mu$, and therefore \citet{ABR2010a} instead reverse the signs of the residuals $e_t = Y_t - \bar{Y}$ and prove that

\begin{equation}\label{eq:ABRintro}
\Prob\big(\phi(\bar Y-\mu)>q_{\alpha_0(1-\delta)}\{\phi(\bar e(\varepsilon))\}+\gamma_1f\big)\le\alpha_0+\Prob\big(\phi(\bar Y-\mu)>f\big),
\end{equation}

where $\bar e(\varepsilon)=n^{-1}\sum_t\varepsilon_te_t$, $\gamma_1\le\sqrt{2\log(2/(\alpha_0\delta))/n}$, and $f\ge0$ is any function \citep[Thm.~3.2]{ABR2010a}. Consequently, the right-hand side does not exceed $\alpha$ if $\Prob(\phi(\bar Y-\mu)>f)\le\alpha-\alpha_0$, and we denote the threshold in formula \eqref{eq:ABRintro} by $t_{\rm ABR}$.

The smallest valid threshold that does not depend on the data is the upper $\alpha$-quantile $Q^*_\alpha$ of $\phi(\bar Y-\mu)$, which \citet{ABR2010a} referred to as the ideal threshold, and which cannot be calculated directly because it depends on the unknown structure of the coordinate dependencies. However, in their simulations \citep[\S4.3]{ABR2010a}, the raw quantile $t_{\rm raw}=q_\alpha\{\phi(\bar e(\varepsilon))\}$, which was computed without $\gamma_1f$, was close to $Q^*_\alpha$, and at the end of their Section 3.1, the following questions were left open: whether it is possible to remove this term “possibly with some additional assumption on the distribution of $Y$” and whether it is possible to obtain “a self-contained result based on the symmetry assumption (SA) alone.”

To show that, in general, the remainder $\gamma_1f$ cannot be eliminated, we consider as an example symmetric, independent data with $K=1$ that are not identically distributed. In this case, $t_{\rm raw}$ is asymptotically a mid-$p$ test (Section~\ref{sec:midp}), and $m=\lceil\log_2(1/\alpha)\rceil$ nonzero observations are sufficient to raise the rejection probability to $2^{1-m}\in(\alpha,2\alpha]$ in the limit (Section~\ref{sec:attain}). For independent and identically distributed symmetric observations with finite fourth moments, $t_{\rm raw}$ is asymptotically valid (Section~\ref{sec:effK}), but it is still not valid for moderate values of $n$, even for normal samples (Section~\ref{sec:gauss}). To answer the second question, we replace $t_{\rm raw}$ with the threshold

\begin{equation}\label{eq:cstarintro}
c^*=\sup\big\{w\ge0:\ w\le q_\alpha\{\phi(\bar e(\varepsilon))+|\bareps|\,w\}\big\},\qquad \bareps=\frac1n\sum_{t=1}^n\varepsilon_t.
\end{equation}

In a nutshell, we replace the unknown term $\bareps\,\bar y$ in the exact sign-flip test for a candidate value $x$ with its maximum possible value $|\bareps|\,w$, and $c^*$ is the largest radius $w$ that the test does not reject after this substitution. When $K=1$, the region of $c^*$ is the outer approximation that \citet[\S VI]{CsajiCampiWeyer2015} compute for their region of sign-perturbed sums using semidefinite programming, and our construction extends it to vector observations and arbitrary seminorms. This threshold is valid in finite samples for any seminorm, assuming only sign symmetry, and requires only a single sorting of the resampled values (similar to $t_{\rm raw}$). The threshold $c^*$ can exceed $t_{\rm ABR}$ only when it exceeds $f$ (Section~\ref{sec:abr}), and in this sense the remainder is indeed “(almost) unnecessary,” although the counterexamples given above show that it cannot be completely dispensed with.

\citet[\S4.4]{ABR2010a} suggest proving that their thresholds are asymptotically optimal, and we prove this for $c^*$ in Section~\ref{sec:efficiency}: the ratio $c^*/Q^*_\alpha$ tends to one for i.i.d.\ rows with finite fourth moments, and even when $K$ grows with $n$ if the coordinates have sub-exponential tails, and if $K$ is fixed and the distribution has an absolutely continuous component, then $c^*/t_{\rm raw}=1+\sqrt{2/(\pi n)}+o_P(n^{-1/2})$ for any seminorm. For symmetric observations with sub-exponential coordinates, the rejection probability of $c^*$ converges to $\alpha$ if $\log K=o(n^{1/3})$ (just as \citet{CCKK2022} obtain this result for $t_{\rm raw}$ in this range).

The studentization of the data “renders void the assumption of independent data points,” on which the proofs in \citet[\S4.1]{ABR2010a} are based, whereas the studentized threshold $c_{\rm stud}$ of Section~\ref{sec:variants} scales the $j$-th coordinate by a factor of $\|y_{\cdot j}\|_2/\sqrt n$, which remains unchanged when the signs are flipped, so that the exact test remains exact. Both $c^*$ and $c_{\rm stud}$ can be used in step-down procedures of \citet{RomanoWolf2005}, and the procedure based on $c^*$ rejects every hypothesis rejected by the Bonferroni correction for one-dimensional thresholds, while asymptotically it rejects nearly the same hypotheses as the oracle procedure (Section~\ref{sec:stepdown}). For the uncentered approach from \citet{ABR2010b}, we show that the one-sided version also controls the family-wise error rate, using an argument that, to the best of our knowledge, has not been encountered in the literature for the maximum over coordinates, although the case $K=1$ appears to be known.

Finally, $c^*$ is infinite when $\alpha<2^{1-n}$, and in Section~\ref{sec:below} we show that, for $K=1$, there exists no interval that is valid at such levels for all symmetric laws, equivariant under translation and has finite endpoints, which proves the conjecture from \citet[Remark~1]{PaulKuchibhotla2026}. The proofs of all results are collected in the Supplementary Material.

\section{Notation}\label{sec:notation}

\subsection{Data and the symmetry assumption}

As in the previous section, we let $Y_1,\dots,Y_n$ be a sample of $n\ge2$ random vectors $Y_t=(Y_{t1},\dots,Y_{tK})$ in $\R^K$, and denote by $\mu=(\mu_1,\dots,\mu_K)$ their common mean. Unlike \citet{ABR2010a}, where a $K\times n$ data matrix was used, we have observations organized as rows of an $n\times K$ matrix. For centered observations, we write $y_t = Y_t - \mu$, and

\[ \bar Y=n^{-1}\sum_tY_t,\qquad \bar y=n^{-1}\sum_ty_t, \]

\[ e_t=Y_t-\bar Y=y_t-\bar y,\qquad t=1,\dots,n,\qquad \sum_te_t=0. \]

\begin{assumption}\label{ass:SA}
The law of $y=(y_1,\dots,y_n)$ is invariant under sign changes of the observations: $\varepsilon\cdot y\eqd y$ for every $\varepsilon\in\{\pm1\}^n$, where $\varepsilon\cdot y=(\varepsilon_1y_1,\dots,\varepsilon_ny_n)$.
\end{assumption}

This is a weaker version of the symmetry assumption (SA) from \citet{ABR2010a}, where it is further assumed that the observations are i.i.d. This allows, for example, dependent observations $y_t=\lambda_t\xi_t$, where $\xi_t$ are independent symmetric variables, and the random scales $\lambda_t$ are independent of $\xi$. In the one-dimensional case, we are dealing with orthant symmetry as described in \citet[p.~1281]{Efron1969}, in which, as he notes, the observations need not have an identical distribution. No conditions regarding moments are required, except for the results on efficiency given in Section~\ref{sec:efficiency}. Also, $\phi:\R^K\to[0,\infty)$ is a seminorm, that is, a subadditive and positively homogeneous function satisfying the condition $\phi(-v)=\phi(v)$.

\subsection{Sign patterns and resampled quantiles}

The resampled quantiles are taken over the set $\mathcal E$ of vectors $\varepsilon\in\{\pm1\}^n$, which we call sign patterns. This is either the full set $\{\pm1\}^n$ with $N=2^n$ elements, or, when $2^n$ is too large, the identity pattern $\mathbf 1$ together with $B$ patterns drawn independently and uniformly from $\{\pm1\}^n$. In the second case, $N=B+1$, and patterns are counted with regard to multiplicity, and the addition of the identity pattern ensures that the exact test from Lemma~\ref{lem:exact} below is valid at level $\alpha$ for any $B$ \citep[Thm.~2]{HemerikGoeman2018}. For $\gamma\in(0,1)$ and $v:\mathcal E\to\R$, we define the resampled quantile $q_\gamma\{v\}$ as the $(k_\gamma+1)$th largest of the $N$ values of $v(\varepsilon)$, where $k_\gamma=\lfloor\gamma N\rfloor$, that is

\begin{equation}\label{eq:q}
q_\gamma\{v\}=q_\gamma\{v(\varepsilon)\}=\inf\{x\in\R:\#\{\varepsilon\in\mathcal E:v(\varepsilon)>x\}\le k_\gamma\}.
\end{equation}

Every pattern $\varepsilon$ determines the values

\[ \bareps=\frac1n\sum_{t=1}^n\varepsilon_t,\qquad \bar y(\varepsilon)=\frac1n\sum_{t=1}^n\varepsilon_ty_t,\qquad \bar e(\varepsilon)=\frac1n\sum_{t=1}^n\varepsilon_te_t,\qquad E(\varepsilon)=\phi(\bar e(\varepsilon)). \]

\begin{lemma}[Exact sign-flip test]\label{lem:exact}
Under Assumption~\ref{ass:SA}, for every $\alpha\in(0,1)$ and every measurable function $W:(\R^K)^n\to\R$,
\[ \Prob\big(W(y)>q_\alpha\{W(\varepsilon\cdot y)\}\big)\le k_\alpha/N\le\alpha. \]

\end{lemma}

This is Lemma 3.1 from \citet{ABR2010a}, which is formulated for an i.i.d.\ sample and $\mathcal E=\{\pm1\}^n$. Their proof relies solely on the invariance $(\varepsilon_ty_t)_t\eqd(y_t)_t$, rather than on independence or the form of $W$, and is therefore applicable to dependent observations. In principle, the duality of tests and confidence regions transforms the lemma into an exact region $C_{\rm ex}$, defined below. However, as noted in \citet[\S3.1]{ABR2010a}, such a region is difficult to compute and does not take the form $\mathcal C(t)$. The unknown values $\bar y(\varepsilon)$ and the computable values $\bar e(\varepsilon)$, which \citet{ABR2010a} use in the proof of their Theorem 3.2, are related by the identity

\begin{equation}\label{eq:id}
\bar y(\varepsilon)=\bar e(\varepsilon)+\bareps\,\bar y,\qquad \varepsilon\in\{\pm1\}^n.
\end{equation}

\subsection{Confidence regions}

A threshold is a measurable function $t=t(Y,\mathcal E)\in[0,\infty]$, and it defines the region

\begin{equation}\label{eq:region}
\mathcal C(t)=\{x\in\R^K:\phi(\bar Y-x)\le t\}.
\end{equation}

The region contains $\mu$ if and only if $\phi(\bar y) \le t$, since $\bar Y - \mu = \bar y$. We call $\Prob(\phi(\bar y)>t)$ the rejection probability of $t$. Its supremum over a class of distributions is the size of $t$ over that class, and $t$ is valid at level $\alpha$ for this class if its size does not exceed $\alpha$. \citet[\S2.2]{ABR2010b} note that empirical centering of the data makes the threshold invariant under translation (their property (TI)). Since the thresholds in this paper are computed using residuals, adding $c\in\R^K$ to each $Y_t$ shifts $\mathcal C(t)$ by $c$, meaning that the regions are translation-equivariant.

Under the hypothesis $H_x:\mu=x$, the observations $Y_t-x$ satisfy Assumption~\ref{ass:SA}, and the exact test from Lemma~\ref{lem:exact} rejects the hypothesis $H_x$ with a probability of at most $k_\alpha/N$. The set $C_{\rm ex}$ of points $x$ that are not rejected, which we call the exact region, satisfies the condition $\Prob(\mu\notin C_{\rm ex})\le k_\alpha/N$. The sign-flipped means of $Y_t - x = e_t + (\bar{Y} - x)$ are equal to $\bar{e}(\varepsilon) + \bareps(\bar{Y} - x)$, and

\begin{equation}\label{eq:Cex}
C_{\rm ex}=\big\{x\in\R^K:\ \phi(\bar Y-x)\le q_\alpha\{\phi(\bar e(\varepsilon)+\bareps(\bar Y-x))\}\big\}.
\end{equation}

In order to obtain some insights, we compare the threshold $c^*$ of Section~\ref{sec:cstar} with the raw centered quantile and the threshold of \citet[Thm.~3.2]{ABR2010a},

\begin{equation}\label{eq:thresholds}
t_{\rm raw}=q_\alpha\{E(\varepsilon)\},\qquad t_{\rm ABR}=q_{\alpha_0(1-\delta)}\{E(\varepsilon)\}+\gamma_1(\alpha_0\delta)f,
\end{equation}

and with the ideal threshold $Q^*_\alpha$ of \citet{ABR2010a}. Here $\alpha_0\in(0,\alpha]$ and $\delta\in(0,1)$ are parameters, $f\ge0$ is any bound with $\Prob(\phi(\bar y)>f)\le\alpha_f=\alpha-\alpha_0$, and $\gamma_1(\eta)=(2B(n,\eta/2)-n)/n$, where $B(n,\eta)$ is the largest $k$ with $\Prob(\mathrm{Bin}(n,1/2)\ge k)\ge\eta$. Among all data-independent thresholds, $Q^*_\alpha$ is the smallest of those that are valid at the $\alpha$ level, and we evaluate the efficiency of a threshold based on its ratio to $Q^*_\alpha$.

\section{The resampled quantile without remainder term}\label{sec:size}

In order to show that the remainder in \eqref{eq:ABRintro} cannot be eliminated under Assumption~\ref{ass:SA}, it suffices to find a single law from this class for which $t_{\rm raw}$ is not valid. We restrict attention to $K=1$, $\phi=|\cdot|$, and $\mathcal E=\{\pm1\}^n$, since a counterexample for $K=1$ is also a counterexample for any $K$ and the max-norm, provided the remaining $K-1$ coordinates are set to zero. The distributions considered in Sections~\ref{sec:midp} and~\ref{sec:attain} are not i.i.d., and since we also want a counterexample in the class considered in \citet{ABR2010a}, in Section~\ref{sec:gauss} we switch to i.i.d.\ normal samples.

\subsection{Limiting rejection probability}\label{sec:midp}

Assumption~\ref{ass:SA} for $K=1$ implies that, given the values $|y_1|,\dots,|y_n|$, any combination of signs is equally likely (Lemma~\ref{lem:orbit}). Thus, when these values are fixed, the data follow

\begin{equation}\label{eq:orbitlaw}
y_t=\sigma_ta_t,\qquad t=1,\dots,n,
\end{equation}

with fixed $a\in[0,\infty)^n$ and $\sigma$ uniformly distributed on $\{\pm1\}^n$, and the size of any threshold over the class of Assumption~\ref{ass:SA} is equal to the supremum of its rejection probability over distributions of this form. We take

\[ a=(b,0,\dots,0), \ b\in[0,\infty)^m, \]

so that all observations except the first $m$ are zero, and let $n\to\infty$ for fixed $b$. Then $n\bar y=S$, where $S=\sum_{t=1}^m\sigma_tb_t$, and $\bar y(\varepsilon)$ depends on $\varepsilon$ only through $\varepsilon_1,\dots,\varepsilon_m$.

The exact test that uses the true value $\mu$ rejects the hypothesis only if the $p$-value $\Prob(|S'|\ge|S|\mid S)$ does not exceed $\alpha$, where $S'$ is an independent copy of $S$. Cases where $|S'|=|S|$ enter this $p$-value with a weight of one, whereas in the raw quantile they are accounted for with a weight close to $1/2$. Indeed, for a pattern with $|\bar y(\varepsilon)|=|\bar y|$, we have $\bar y(\varepsilon)=\pm\bar y$, and by \eqref{eq:id}, $|\bar e(\varepsilon)|=|1\mp\bareps|\,|\bar y|$. Depending on the sign of $\bareps$, this value is either greater than or equal to $|\bar y|$ or strictly less than $|\bar y|$, and for large values of $n$, the sign of $\bareps$ is determined by the remaining $n-m$ signs, and both signs of $\bareps$ occur with a probability close to $1/2$. Thus, the limit is described by the mid-$p$-value of $|S|$ \citep[eq.~(3)]{Wang2024},

\[ F_b(x)=\tfrac12\,\Prob(|S'|\ge x)+\tfrac12\,\Prob(|S'|>x),\qquad x\ge0. \]

Moreover, the inequality $\Prob(F_b(|S|)\le\alpha)\le2\alpha$ holds for mid-$p$-values at each level, as follows from \citet[\S3 and Prop.~5.1]{Wang2024}.

\begin{theorem}\label{thm:midp}
Let $K=1$, $\alpha\in(0,1/2)$, and let $y$ follow \eqref{eq:orbitlaw} with $a=(b,0,\dots,0)$ and $b\in[0,\infty)^m$ fixed. Then

\[ \Prob\big(F_b(|S|)<\alpha\big)\ \le\ \liminf_{n\to\infty}\Prob(|\bar y|>t_{\rm raw})\ \le\ \limsup_{n\to\infty}\Prob(|\bar y|>t_{\rm raw})\ \le\ \Prob\big(F_b(|S|)\le\alpha\big)\ \le\ 2\alpha. \]

In particular,
\[ \Prob\big(F_b(|S|)=\alpha\big)=0\quad\Longrightarrow\quad\lim_{n\to\infty}\Prob(|\bar y|>t_{\rm raw})=\Prob\big(F_b(|S|)<\alpha\big). \]
\end{theorem}

\subsection{A law that attains the bound}\label{sec:attain}

The bound $2\alpha$ in Theorem~\ref{thm:midp} is attained for a dyadic $\alpha$ with the simplest choice of $b$, while for other values of $\alpha$, the same choice yields a limit in the interval $(\alpha,2\alpha)$, for example, $1/16=1.25\alpha$ when $\alpha=0.05$. If we take $m=\lceil\log_2(1/\alpha)\rceil$ and $b=\mathbf 1_m$, then $|S|=m$ with probability $2^{1-m}$ and $F_b(m)=2^{-m}\le\alpha$, whereas for every smaller value of $|S|$, we have $F_b\ge2^{1-m} >\alpha$. If $\alpha>2^{-m}$, then by Theorem~\ref{thm:midp}, $\mathcal C(t_{\rm raw})$ does not contain $\mu$ in the limit exactly when all $m$ signs coincide, which is an event with probability $2^{1-m}\in(\alpha,2\alpha)$. When $\alpha=2^{-m}$ is a dyadic number, $F_b(m)$ equals $\alpha$, and the theorem specifies only the bounds $0$ and $2\alpha$ for the limit points. In the proof, we instead directly count the number of rejecting patterns, and this count also shows that the rejection probability equals its limit for any large finite number $n$.

\begin{theorem}\label{thm:attain}
Let $K=1$, $\alpha\in(0,1/2)$, $m=\lceil\log_2(1/\alpha)\rceil$, and let $y$ follow \eqref{eq:orbitlaw} with $a=(\mathbf 1_m,0,\dots,0)$. There is $n_0(\alpha)$ such that

\[ \Prob(|\bar y|>t_{\rm raw})=2^{1-m}\in(\alpha,2\alpha]\quad\text{for every }n\ge n_0(\alpha),\qquad 2^{1-m}=2\alpha\iff\alpha=2^{-m}. \]
\end{theorem}

Under this distribution, $\bar e(\varepsilon)$ is a function of the first $m$ signs and the sum of the remaining $n-m$ signs, and therefore the rejection probability reduces to binomial sums and can be calculated exactly for all $2^n$ sign patterns. For $\alpha\in\{1/4, 1/8, 1/10, 1/20, 1/64\}$, it is equal to $2^{1-m}$ for $n=10, 20, 50, 100$, and $200$, as well as for $n=m$, except in the case where $\alpha=1/64$, where it is equal to $7/32$, so in these examples the conclusion of the theorem holds for smaller values of $n$ than those required in its proof.

In his paper on posterior predictive $p$-values, \citet[after Eq.~(5.6)]{Meng1994} constructs a distribution function that attains the bound $2\alpha$, but argues that “in most practical situations the factor $2$ is simply too large.” In fact, Theorem~\ref{thm:attain} yields the same factor when using the resampling threshold itself under an explicitly specified sign-symmetric distribution with $m$ nonzero observations.

\subsection{Normal samples}\label{sec:gauss}

The distribution underlying Theorem~\ref{thm:attain} is not i.i.d., because all observations except $m$ are zero. However, for i.i.d.\ data with a normal distribution and a small $n$, the threshold $t_{\rm raw}$ is not valid either. The reason is simple: it never exceeds the largest value of $|\bar e(\varepsilon)|$, which is attained when each $\varepsilon_t$ has the same sign as $e_t$, and according to the Cauchy--Bunyakovsky--Schwarz inequality

\[ t_{\rm raw}\le\max_{\varepsilon}|\bar e(\varepsilon)|=\frac1n\sum_{t=1}^n|e_t|\le\frac{\|e\|_2}{\sqrt n}. \]

Thus, regardless of the level, the value of $t_{\rm raw}$ is exceeded as soon as $|\bar y| > \|e\|_2/\sqrt{n}$, and in terms of the standard $t$-statistic

\[ T=\frac{\sqrt n\,\bar y}{s},\qquad s^2=\frac{\|e\|_2^2}{n-1}, \]

this event is $\{|T|>\sqrt{n-1}\}$, and under normality $T$ has the $t$-distribution with $n-1$ degrees of freedom.

\begin{theorem}\label{thm:gauss}
Let $K=1$ and let $Y_1,\dots,Y_n$ be i.i.d.\ normal with mean $\mu$ and positive variance. Then for every $\alpha\in(0,1)$

\[ \Prob(|\bar y|>t_{\rm raw})\ge\Prob\big(|T|>\sqrt{n-1}\big)\quad(n\ge2),\qquad \Prob\big(|T|>\sqrt{n-1}\big)>2^{1-n}\quad(n\ge3). \]
\end{theorem}

The lower bound does not depend on $\alpha$, and the inequality from \citet{Wendel1948} for gamma function ratios yields

\[ \Prob\big(|T|>\sqrt{n-1}\big)\ge\frac{2^{-(n-1)/2}}{\sqrt{\pi n/2}}, \]

so it decreases like $2^{-n/2}$, that is, much more slowly than $2^{1-n}$. The comparison with $2^{1-n}$ is relevant because, at lower levels, there are no intervals that are translation-equivariant and have finite endpoints (Theorem~\ref{thm:imposs}). At every level $\alpha$ with $2^{1-n}\le\alpha<\Prob(|T|>\sqrt{n-1})$, the threshold $c^*$ from Section~\ref{sec:cstar} is valid for all symmetric distributions (Theorem~\ref{thm:cstar}), whereas $t_{\rm raw}$ is not valid even for normal data.

\begin{figure}[H]
\centering
\includegraphics[width=0.65\textwidth]{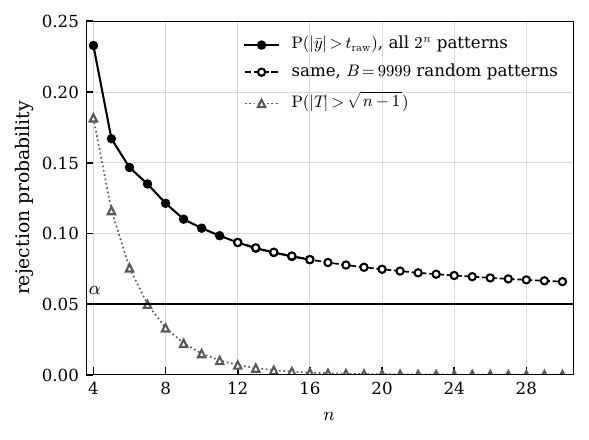}
\caption{Rejection probability of $t_{\rm raw}$ for i.i.d.\ normal data at the $\alpha=0.05$ level, with all $2^n$ sign patterns for $n\le16$ and with $9999$ random patterns for $n\ge12$.}\label{fig:gauss}
\end{figure}

At the $\alpha=0.05$ level, the bound from Theorem \ref{thm:gauss} exceeds $0.05$ only for $n\le6$, since it replaces $t_{\rm raw}$ by the largest of the values of $|\bar e(\varepsilon)|$, whereas for larger values of $n$ we compute the rejection probability itself. Given the direction of $e$, this is the probability in the tail of the $t_{n-1}$ distribution, which can be computed exactly for all $2^n$ sign patterns ($n=4,\dots,16$), as well as for the identity pattern and $B=9999$ random patterns ($n=12,\dots,30$), since the random patterns are selected independently of the data. We average this value over $2000$ simulated directions, yielding standard errors of less than $0.001$, and for $n=12,\dots,16$, these two sets of patterns yield values that differ by no more than $0.0003$. For $n=10$, the probability of a rejection is $0.104$, which is approximately equal to $2\alpha$, for $n=16$ it is $0.081$, and for $n=30$ it is still $0.066$ (Figure~\ref{fig:gauss}), so $t_{\rm raw}$ is not valid at the $0.05$ level for normally distributed samples of moderate size.

\section{A threshold valid under sign symmetry}\label{sec:cstar}

If a candidate value $x$ is given, it is easy to check whether it belongs to the exact region $C_{\rm ex}$ from \eqref{eq:Cex}: for a given $x$, we compute the quantile from \eqref{eq:Cex} and compare it with $\phi(\bar Y-x)$. Thus, a grid search for points $x$ would yield $C_{\rm ex}$, but the number of grid points grows exponentially with $K$. Following the idea in \citet[\S VI]{CsajiCampiWeyer2015}, we instead seek an outer approximation of $C_{\rm ex}$ of the form $\mathcal C(t)$, with a threshold computed from residuals at the same computational cost as for $t_{\rm raw}$, and any such region is valid provided only that Assumption~\ref{ass:SA} holds.

\subsection{The fixed-point threshold}\label{sec:fixedpoint}

Suppose that $x\in C_{\rm ex}$, and let $w=\phi(\bar Y-x)$, so that definition \eqref{eq:Cex} and the triangle inequality $\phi(\bar e(\varepsilon)+\bareps(\bar Y-x))\le E(\varepsilon)+|\bareps|\,w$, applied to each pattern, yield

\[ w\le q_\alpha\{\phi(\bar e(\varepsilon)+\bareps(\bar Y-x))\}\le q_\alpha\{E(\varepsilon)+|\bareps|\,w\}, \]

since the quantile preserves the pointwise order. The right-hand side no longer depends on the direction of $\bar Y-x$, but depends only on $w$, and we define $c^*$ as the largest value of $w$ that satisfies this inequality,

\begin{equation}\label{eq:cstar}
c^*=\sup\big\{w\ge0:\ w\le q_\alpha\{E(\varepsilon)+|\bareps|\,w\}\big\}.
\end{equation}

Then, for every $x\in C_{\rm ex}$, $\phi(\bar Y-x)\le c^*$ holds, that is, $C_{\rm ex}\subseteq\mathcal C(c^*)$.

When $|\bareps|<1$, the quantity $E(\varepsilon)+|\bareps|\,w$ increases as $w$ increases with a slope of $|\bareps|<1$, and therefore it is not less than $w$ if and only if

\[ w\le r(\varepsilon):=\frac{E(\varepsilon)}{1-|\bareps|}. \]

The constant patterns $\pm\mathbf 1$ have $E(\varepsilon)=0$ and $|\bareps|=1$, and their value is equal to $w$ for any $w$. From the definition of a quantile, it follows that $q_\alpha\{E(\varepsilon)+|\bareps|\,w\}\ge w$ if and only if $k_\alpha+1$ or more patterns have a value of $w$ or greater, and if we consider the constant patterns separately, we obtain the following result.

\begin{lemma}\label{lem:closed}
Let $m_1$ be the number of patterns $\varepsilon\in\mathcal E$ with $|\bareps|=1$, counted with multiplicity. If $k_\alpha\ge m_1$, then $c^*$ is the $(k_\alpha+1-m_1)$-th largest of the values $r(\varepsilon)$ over the patterns with $|\bareps|<1$. If $k_\alpha<m_1$, then $c^*=+\infty$.
\end{lemma}

In the complete set $\{\pm1\}^n$, the only constant patterns are $\pm\mathbf 1$, so $m_1=2$, and $c^*$ is finite when $k_\alpha\ge2$. The Monte Carlo set contains the identity and every random pattern equal to $\pm\mathbf 1$, so $m_1-1$ has a binomial distribution with mean $B\,2^{1-n}$. For values of $B=1999$ and $\alpha=0.05$, we obtain $k_\alpha=100$, whereas for $n=10$, this mean is only about $4$, and therefore the condition $k_\alpha\ge m_1$ holds with probability close to one as long as $n$ is not small. In both cases, calculating $c^*$ requires one sort of $N$ values of $r(\varepsilon)$, just as in the case of $t_{\rm raw}$. In passing, note that $c^*\ge t_{\rm raw}$, since $r(\varepsilon)\ge E(\varepsilon)$ for $|\bareps|<1$, and $E(\varepsilon)=0$ for constant patterns (see also Lemma~\ref{lem:sandwich}). In the location model with $K=1$, the region $\mathcal C(c^*)$ coincides with the outer approximation from \citet[\S VI]{CsajiCampiWeyer2015}, where each of their semi-definite programs returns the value $r(\varepsilon)^2$ of a single pattern.

\begin{theorem}\label{thm:cstar}
Under Assumption~\ref{ass:SA}, $C_{\rm ex}\subseteq\mathcal C(c^*)$ and

\[ \Prob\big(\phi(\bar y)>c^*\big)\le k_\alpha/N\le\alpha. \]

For $\mathcal E=\{\pm1\}^n$ the threshold $c^*$ is finite if and only if $\alpha\ge2^{1-n}$.
\end{theorem}

\subsection{Studentized threshold}\label{sec:variants}

Section 4.1 of \citet{ABR2010a} discusses the studentization method, in which each coordinate is divided by its own scale. The argument presented in Section~\ref{sec:fixedpoint} also applies to this case if the triangle inequality is replaced by a bound for the sign-flipped statistic.

To do this, we write $y_{\cdot j}=(y_{1j},\dots,y_{nj})$ and $e_{\cdot j}=(e_{1j},\dots,e_{nj})$ for the $j$-th column, and divide $\bar y_j$ by $\|y_{\cdot j}\|_2/\sqrt n$, a scaling factor that remains unchanged when the signs are flipped, so Lemma~\ref{lem:exact} can still be applied. As $\|y_{\cdot j}\|_2^2=\|e_{\cdot j}\|_2^2+n\bar y_j^2\ge\|e_{\cdot j}\|_2^2$, from equality \eqref{eq:id} we see that

\[ \frac{\sqrt n\,|\bar y_j(\varepsilon)|}{\|y_{\cdot j}\|_2}\le\frac{\sqrt n\,|\bar e_j(\varepsilon)|}{\|e_{\cdot j}\|_2}+|\bareps|\,\frac{\sqrt n\,|\bar y_j|}{\|y_{\cdot j}\|_2}, \]

and this inequality replaces the triangle inequality. We define

\[ c_{\rm stud}=\sup\Big\{w\ge0:\ w\le q_\alpha\Big\{\max_j\frac{\sqrt n\,|\bar e_j(\varepsilon)|}{\|e_{\cdot j}\|_2}+|\bareps|\,w\Big\}\Big\}, \]

with the convention $0/0=0$ for columns with $\|e_{\cdot j}\|_2=0$, where $\bar e_j(\varepsilon)=0$ for every $\varepsilon$.

As $\|y_{\cdot j}\|_2$ depends on $\mu$, a bound on the studentized statistic is not an interval for $\mu_j$. In the notation of \citet[p.~1279]{Efron1969}, $\sqrt n\,|\bar y_j|/\|y_{\cdot j}\|_2$ is equal to $|S_n|/\sqrt n$, and Student's statistic is a monotonic function $T_n=S_n[(n-1)/(n-S_n^2)]^{1/2}$ of $S_n$. Inverting this bound, we obtain the following interval: if $\sqrt n\,|\bar y_j|/\|y_{\cdot j}\|_2\le c<1$, then $n\bar y_j^2\le c^2\big(\|e_{\cdot j}\|_2^2 + n\bar y_j^2\big)$, that is, $|\bar Y_j - \mu_j| \le c \,\|e_{\cdot j}\|_2 / \sqrt{n(1 - c^2)}$.

\begin{theorem}\label{thm:variants}
Under Assumption~\ref{ass:SA}, $\Prob\big(\max_j\sqrt n\,|\bar y_j|/\|y_{\cdot j}\|_2>c_{\rm stud}\big)\le\alpha$ and

\[ \Prob\Big(c_{\rm stud}<1\ \text{and}\ |\bar Y_j-\mu_j|>\frac{c_{\rm stud}}{\sqrt{1-c_{\rm stud}^2}}\,\frac{\|e_{\cdot j}\|_2}{\sqrt n}\ \text{for some }j\Big)\le\alpha. \]

\end{theorem}

The condition $c_{\rm stud} < 1$ is natural and unavoidable if one wants an interval for $\mu_j$, because $\sqrt{n}\,|\bar{y}_j| \le \|y_{\cdot j}\|_2$ always holds, so a threshold equal to one or greater does not exclude any value of $\mu_j$.

\section{Nonexistence below the level \texorpdfstring{$2^{1-n}$}{2\^{}(1-n)}}\label{sec:below}

From Theorem~\ref{thm:cstar}, we see that the threshold $c^*$ is infinite if and only if $\alpha<2^{1-n}$, and in this section our goal is to show that at these levels there is no interval that is translation-equivariant and has finite endpoints, regardless of the threshold used.

Under Assumption~\ref{ass:SA} with nonzero observations, all $n$ signs of $y_1, \dots, y_n$ coincide with probability $2^{1-n}$. If this event occurs, the residuals can be identical for any value of $|\bar Y-\mu|$, as shown in the proof of Theorem~\ref{thm:imposs}(i) below. Therefore, there is no threshold computed from them that can bound $|\bar Y-\mu|$.

Theorem 2 of \citet{Lanke1974} states that when $K=1$, no locally bounded interval is valid at levels below $2^{1-n}$ for i.i.d.\ random variables drawn from a continuous distribution that is symmetric about $\mu$. However, \citet[Remark~1]{PaulKuchibhotla2026} point out that this statement is false, although the proof in \citet {Lanke1974} shows that it is impossible to construct such an interval at levels below $2^{-n}$, using a distribution in which, with probability $2^{-n}$, all observations lie in a single cluster at a large distance from $\mu$. In the same remark, the authors conjecture that for $\alpha<2^{1-n}$, there exists no locally bounded interval for the median that is translation-equivariant, a property they call shift invariance. In part (ii) of the following theorem, we prove this conjecture for every interval with finite endpoints, which is a weaker condition than local boundedness. Therefore, the level $2^{1-n}$ is optimal.

\begin{theorem}\label{thm:imposs}
Let $K=1$ and $\alpha<2^{1-n}$.
\begin{enumerate}[label=(\roman*)]
\item Let $h$ be a measurable function of the residuals $e=(e_1,\dots,e_n)$ with values in $[0,\infty]$ and $h(-e)=h(e)$. If $\Prob(|\bar y|>h(e))\le\alpha$ for every distribution that satisfies Assumption~\ref{ass:SA}, then $h(e)=+\infty$ for every $e\in\mathbf 1^\perp=\{v\in\R^n:\sum_tv_t=0\}$, that is, for every possible value of the residual vector.
\item There is no interval $[L(Y),U(Y)]$ with measurable endpoints taking finite values that is equivariant with respect to translations \citep[\S6.11]{LehmannRomano2005}, meaning that $L(Y+c\mathbf 1)=L(Y)+c$ and $U(Y+c\mathbf 1)=U(Y)+c$ for all $c\in\R$, and that satisfies $\Prob(\mu\in[L(Y),U(Y)])\ge1-\alpha$ for every i.i.d.\ distribution with a continuous, everywhere positive density symmetric about $\mu$.
\end{enumerate}
\end{theorem}

Part (i) explains what happens to both thresholds below the level $2^{1-n}$: $t_{\rm raw}$ with $\mathcal E=\{\pm1\}^n$ is always finite and, therefore, invalid there, while $c^*$ is infinite precisely at these levels. This argument also holds for a random set of patterns chosen independently of the data. Both parts rely on the event that all observations lie on one side of $\mu$, and for intervals, the two-sided version of the Lanke construction with clusters at the points $c$ and $-c$ is applied, which doubles the level $2^{-n}$ reached in Lanke’s proof to $2^{1-n}$. As these distributions have a median of $\mu$, the theorem also excludes translation-equivariant intervals with finite endpoints for the median of a continuous distribution, at levels below $2^{1-n}$.

\section{Comparison with the threshold of Arlot, Blanchard and Roquain}\label{sec:abr}

Both $t_{\rm ABR}$ and $c^*$ are derived from identity \eqref{eq:id} and differ in how they constrain the term $\bareps\,\bar y$. \citet{ABR2010a} divide the level into $\alpha=\alpha_0(1-\delta)+\alpha_0\delta+\alpha_f$, and these three parts correspond to the resampled quantile of $E(\varepsilon)$, patterns with $|\bareps|>\gamma_1(\alpha_0\delta)$, and the external bound $f$ for $\phi(\bar y)$, respectively. The main idea of the fixed point \eqref{eq:cstar} is to replace $f$ with $c^*$ itself, and the same decomposition of the level bounds $c^*$ from above.

\begin{lemma}\label{lem:sandwich}
Let $\beta\in(0,\alpha)$ be such that $k_\beta=\lfloor\beta N\rfloor$ is not less than the number $m_1$ of patterns in $\mathcal E$ for which $|\bareps|=1$. Then

\[ q_\alpha\{E(\varepsilon)\}\le c^*\le\frac{q_{\alpha-\beta}\{E(\varepsilon)\}}{1-q_\beta\{|\bareps|\}}. \]

Moreover, for $\mathcal E=\{\pm1\}^n$, $q_\beta\{|\bareps|\}\le\sqrt{2\log(2/\beta)/n}$.
\end{lemma}

In our opinion, the most interesting case is when we take $\mathcal E=\{\pm1\}^n$ and $\beta=\alpha-\alpha_0(1-\delta)=\alpha_0\delta+\alpha_f$, and assume that $k_\beta\ge2$. The proof of Theorem~\ref{thm:abr}, given below, shows that at most $k_{\alpha_0\delta}$ patterns have $|\bareps|>\gamma_1=\gamma_1(\alpha_0\delta)$, so that $q_\beta\{|\bareps|\}\le\gamma_1$, and if $\gamma_1<1$, then by Lemma~\ref{lem:sandwich} it follows that

\[ c^*\le q_{\alpha_0(1-\delta)}\{E(\varepsilon)\}+\gamma_1\,c^*, \]

which is $t_{\rm ABR}$ in \eqref{eq:thresholds} with $f$ replaced by $c^*$.

In fact, \citet[Cor.~3.3]{ABR2010a} already reduce the role of $f$ by iterating their Theorem 3.2, so that $f$ appears only with the coefficient $\gamma_J$, which is accompanied by a loss of level, since the levels of $J$ steps are summed. The fixed point \eqref{eq:cstar} can be regarded as the endpoint of this iteration, at which the threshold itself takes the place of $f$. Since the inequality $E(\varepsilon)+|\bareps|\,w\ge w$ is imposed within a single quantile at level $\alpha$, $c^*$ requires neither a bound $f$ nor a level split, and such a split appears only in Lemma~\ref{lem:sandwich} as a tool for comparing $c^*$ with $t_{\rm ABR}$.

\begin{theorem}\label{thm:abr}
Let $\mathcal E=\{\pm1\}^n$ and let $c^*$ be computed at the level $\alpha=\alpha_0+\alpha_f$ of $t_{\rm ABR}$. Then

\[ c^*\le t_{\rm ABR}+\gamma_1\,(c^*-f)_+, \]

so $c^*\le t_{\rm ABR}$ whenever $c^*\le f$. If $\gamma_1<1$, then also $c^*\le t_{\rm ABR}/(1-\gamma_1)$.
\end{theorem}

Thus, $c^*$ can exceed $t_{\rm ABR}$ only when $c^*>f$, and in our simulation results presented in Section~\ref{sec:numerics}, the median of $c^*$ lies below the median of $t_{\rm ABR}$ for all considered values of $n$ and $K$ (Figure~\ref{fig:width}).

\section{Asymptotic efficiency}\label{sec:efficiency}

It can legitimately be asked how much the regions expand when validity is required only under Assumption~\ref{ass:SA}, since $c^*$ is never less than $t_{\rm raw}$, and the guarantee implied by Theorem~\ref{thm:cstar} can lead to the appearance of regions whose width is greater than necessary. In this section, we show that the additional width asymptotically vanishes for independent rows with finite moments: $c^*/Q^* _\alpha\to1$ both for a fixed $K$ (Section~\ref{sec:effK}), where we also find the exact size of the extra width, and for $K$ growing with $n$ (Section~\ref{sec:effHD}), and for a fixed $K$, we show that this is also true when the sign patterns are chosen randomly (Section~\ref{sec:mc}).

\subsection{Fixed dimension}\label{sec:effK}

The key point in the proof of the following theorem is that, by the bound of \citet{Bentkus2005} for convex sets, the distributions of $\sqrt n\,\phi(\bar y)$ and, under sign flipping, $\sqrt n\,E(\varepsilon)$ are comparable to Gaussian distributions. Lemma~\ref{lem:sandwich} bounds $c^*$ by the two quantiles of $E(\varepsilon)$, and the Gaussian isoperimetric inequality \citep{Borell1975,SudakovTsirelson1978} translates these comparisons into relative errors of the quantiles.

\begin{theorem}\label{thm:effK}
Let the rows $y_1,\dots,y_n$ be i.i.d.\ and symmetric about $0$, with $\E\|y_t\|_2^4<\infty$ and a positive definite covariance matrix $\Sigma$. Let $K$ be fixed, $\phi$ a seminorm that is not identically zero, $\alpha\in(0,1/2)$ and $\mathcal E=\{\pm1\}^n$. Then

\[ \frac{c^*}{Q^*_\alpha}=1+O_P\Big(\sqrt{\frac{\log n}{n}}\Big)\qquad\text{and}\qquad\Prob\big(\phi(\bar y)>c^*\big)\to\alpha. \]

Moreover, $\Prob(\phi(\bar y)>(1+\kappa_n)c^*)\to\alpha$ for every deterministic sequence $\kappa_n\to0$.
\end{theorem}

The lower bound for $c^*$ is a bound for $t_{\rm raw}=q_\alpha\{E(\varepsilon)\}$, so using the same argument, it can be shown that the rejection probability of $t_{\rm raw}$ converges to $\alpha$ for i.i.d.\ rows. Thus, the failures described in Section~\ref{sec:size} occur either for small values of $n$, as in the case of normal samples from Section~\ref{sec:gauss}, or for distributions that are not i.i.d., as in Theorem~\ref{thm:midp}. 

In addition, \citet[Thm.~3]{WeyerCampiCsaji2017} show that for $K=1$, the region of sign-perturbed sums, which, for the location model, is a variant of $C_{\rm ex}$, asymptotically lies within a slightly enlarged normal confidence ellipsoid without specifying the rate. However, the main advantage of Theorem~\ref{thm:effK}, in combination with Proposition~\ref{prop:mc} for random patterns, is the rate $\sqrt{\log n/n}$ for vectors and any seminorm. The logarithm in this rate arises from bounding $|\bareps|$ by its quantile in Lemma~\ref{lem:sandwich}, whereas the averaging of $|\bareps|$ over patterns, on the contrary, eliminates it and determines the next term in the expression $c^*/t_{\rm raw}$.

\begin{theorem}\label{thm:sharp}
Let the assumptions of Theorem~\ref{thm:effK} hold, and let the law of $y_1$ have a nonzero absolutely continuous component. Then
\[ \sqrt n\,\Big(\frac{c^*}{t_{\rm raw}}-1\Big)\to\sqrt{\frac2\pi}\quad\text{in probability}\qquad\text{and}\qquad\frac{c^*}{Q^*_\alpha}=1+O_P(n^{-1/2}). \]
\end{theorem}

Note that the ratio $t_{\rm raw}/Q^*_\alpha$ is also equal to $1+O_P(n^{-1/2})$. Therefore, for distributions with an absolutely continuous component (which excludes lattice laws, such as the one in Figure~\ref{fig:level}(a)), the guarantee for a finite sample extends the region by approximately $0.80/\sqrt n$ in relative terms. The constant $\sqrt{2/\pi}$ is the limit of the expression $\sqrt n\,\E|\bareps|$ and reflects the fact that in \eqref{eq:cstar} each pattern is considered with its own $|\bareps|$, whereas $t_{\rm ABR} $ uses the quantile $\gamma_1(\alpha_0\delta)$ of $|\bareps|$, which is larger by a factor of the order of $\sqrt{\log(1/(\alpha_0\delta))}$.

\subsection{Growing dimension}\label{sec:effHD}

In applications such as neuroimaging, $K$ is usually much larger than $n$, so we will now let $K$ grow with $n$ and take the maximum norm, which will yield simultaneous intervals for the coordinates of $\mu$. Recall that the upper bound from Lemma~\ref{lem:sandwich} divides the quantile at level $\alpha-\beta$ by $1-q_\beta\{|\bareps|\}$, and for $c^*$ we must additionally show that this factor changes $\Prob(\|\bar y\|_\infty>t)$ by only $o(1)$. We do this using the anti-concentration inequality from \citet {Nazarov2003} for the Gaussian approximation of $\sqrt n\,\bar y$ \citep[Lemma~A.1]{CCK2017} or, without any Gaussian approximation, using the anti-concentration inequality for Rademacher sums in Lemma~J.4 of the supplement of \citet{CCKK2022}.

\begin{theorem}\label{thm:effHD}
Let $\phi=\|\cdot\|_\infty$, $\alpha\in(0,1/2)$ and $\mathcal E=\{\pm1\}^n$. Let the rows $y_t$ be independent and symmetric about $0$, with $n^{-1}\sum_t\E y_{tj}^2\in[\underline\sigma^2,\bar\sigma^2]$ for all $j$ and constants $0<\underline\sigma\le\bar\sigma$, and with $\E\exp(|y_{tj}|/b)\le2$ for a constant $b$ and all $t$ and $j$.
\begin{enumerate}[label=(\alph*)]
\item If $\log K=o(n^{1/3})$, then $\Prob(\|\bar y\|_\infty>c^*)\to\alpha$.
\item If $\log K=o(n^{1/5})$, then moreover $c^*/Q^*_\alpha\to1$ in probability.
\end{enumerate}
Furthermore, in (b) $\Prob(\|\bar y\|_\infty>(1+\kappa_n)c^*)\to\alpha$ for every deterministic sequence $\kappa_n$ with $|\kappa_n|\le1/2$ and $\kappa_n\log(3K)\to0$.
\end{theorem}

We note that the convergence of the rejection probability holds in the same range as the result from \citet{CCKK2022} for $t_{\rm raw}$, and instead of a Gaussian approximation, its proof uses the anti-concentration inequality for Rademacher sums, applied conditionally on the rows up to their signs. The ratio $c^*/Q^*_\alpha$ involves quantiles of the sign-flip distribution at two levels, and comparing them requires a lower bound on its density. In order to obtain it, we use the Gaussian approximations from \citet[Lemmas~4.3 and~4.6]{CCKK2022}, and the range of $K$ is the one in which these results hold.

\subsection{Efficiency of the studentized threshold}\label{sec:effS}

For coordinates with different scales, we will now consider the studentized threshold $c_{\rm stud}$ from Section~\ref{sec:variants} and recall that it compares each coordinate with its own scale. Our task will be to compare it with the upper $\alpha$-quantile $Q^{\rm stud}_\alpha$ of the expression $\max_j\sqrt n\, |\bar y_j|/\|y_{\cdot j}\|_2$. The rate $n^{-1/3}$ in the following theorem follows from the bound on the errors of the estimated scales $\|y_{\cdot j}\|_2/\sqrt n$ using Chebyshev’s inequality, and we do not claim that it is optimal.

\begin{theorem}\label{thm:effS}
Let $\mathcal E=\{\pm1\}^n$, $\alpha\in(0,1/2)$, and let the assumptions of Theorem~\ref{thm:effK} hold. Then

\[ \frac{c_{\rm stud}}{Q^{\rm stud}_\alpha}=1+O_P(n^{-1/3})\qquad\text{and}\qquad\Prob\Big(\max_j\frac{\sqrt n\,|\bar y_j|}{\|y_{\cdot j}\|_2}>c_{\rm stud}\Big)\to\alpha. \]

\end{theorem}

\subsection{Random sign patterns}\label{sec:mc}

The theorems above use all $2^n$ sign patterns, while in practice, as is often done for resampled quantities, $c^*$ is computed by a Monte Carlo approximation from the identity and $B$ random patterns (Section~\ref{sec:notation}). We write $\mathcal E_B$ for this set and $c^*(\gamma;\mathcal E)$ for the threshold \eqref{eq:cstar} at the level $\gamma$ computed from a set $\mathcal E$.

\begin{proposition}\label{prop:mc}
Let $u\ge2/\sqrt B$ and $2u/\sqrt B<\alpha<1-2u/\sqrt B$. Conditionally on the data, with probability at least $1-2e^{-2u^2}$ over the random patterns,
\[ c^*\big(\alpha+2u/\sqrt B;\{\pm1\}^n\big)\le c^*(\alpha;\mathcal E_B)\le c^*\big(\alpha-2u/\sqrt B;\{\pm1\}^n\big). \]
If $B=B_n\to\infty$, $\mathcal E=\mathcal E_{B_n}$ and the assumptions of Theorem~\ref{thm:effK} hold, then $c^*/Q^*_\alpha=1+O_P(\sqrt{\log n/n}+B_n^{-1/2})$, and the statements on the rejection probability hold as there.
\end{proposition}

The guarantee for a finite sample, as stated in Theorem~\ref{thm:cstar}, holds for any $B$, since the identity pattern is included. Therefore, the proposition applies only to the width of the region. In particular, the rate specified in Theorem~\ref{thm:effK} holds provided that $B_n\ge n/\log n$.

\section{Step-down procedures}\label{sec:stepdown}

In the companion paper \citet{ABR2010b}, thresholds obtained by the resampling method are used for multiple testing, and each threshold from Section~\ref{sec:cstar} defines a step-down procedure. \citet{ABR2010b} tested the hypotheses $H_j:\mu_j=0$, $j=1,\dots,K$, and controlled the probability of rejecting at least one true hypothesis, that is, the family-wise error rate (FWER). Since a fixed value of each $\mu_j$ is used in the hypotheses, the signs can be changed on the observations themselves rather than on the residuals. \citet{ABR2010b} call this the \emph{uncentered approach}, and we compare it with the centered thresholds from Section~\ref{sec:cstar}.

Roughly, the idea of step-down procedures can be explained as follows. Let $T_j$ be a statistic for hypothesis $H_j$, and for each set of indices $C$, let $t(C)$ be a threshold computed from the data, such that $t(C) \le t(C')$ whenever $C \subseteq C'$. The procedure begins with all indices, $C_0 = \{1, \dots, K\}$, and repeats

\begin{equation}\label{eq:sdrec}
C_{i+1}=\{j\in C_i:\ T_j\le t(C_i)\},\qquad i\ge0,
\end{equation}

until the set ceases to change, after which hypotheses not included in the final set are rejected. To fix ideas, suppose that $K=2$, that $\mu_2$ is large, and that $\mu_1>0$ is moderate. In the first stage, the threshold $t(\{1,2\})$ is computed from both coordinates, so hypothesis $H_2$ is rejected, while hypothesis $H_1$ may remain. In the second stage, the threshold is $t(\{1\})\le t(\{1,2\})$, and then $H_1$ can also be rejected, which could not have been achieved using a single-step procedure with $t(\{1,2\})$. A true hypothesis can be rejected only if $\max_{j\in I_0}T_j>t(I_0)$, where $I_0$ is the set of true hypotheses. This is the argument underlying Theorem 3 from \citet{RomanoWolf2005}, which \citet[Thm.~3.2]{ABR2010b} reformulate with a brief summary of the proof in this setting, and which reduces control over the FWER to the validity of a single threshold $t(I_0)$.

\begin{theorem}\label{thm:fwer}
Let Assumption~\ref{ass:SA} hold and let all thresholds be computed from the same set $\mathcal E$ of sign patterns. Each of the following step-down procedures controls the FWER at level $\alpha$.
\begin{enumerate}[label=(\alph*)]
\item $H_j:\mu_j=0$, $T_j=|\bar Y_j|$ and $t(C)=c^*$ computed for $\phi(v)=\max_{j\in C}|v_j|$.
\item $H_j:\mu_j=0$, $T_j=\sqrt n\,|\bar Y_j|/\|Y_{\cdot j}\|_2$ and $t(C)=c_{\rm stud}$ computed for the coordinates in $C$.
\item $H_j:\mu_j=0$, $T_j=|\bar Y_j|$ and $t(C)=q_\alpha\{\max_{j\in C}|n^{-1}\sum_t\varepsilon_tY_{tj}|\}$, or $H_j:\mu_j\le0$, $T_j=\bar Y_j$ and $t(C)=q_\alpha\{\max_{j\in C}n^{-1}\sum_t\varepsilon_tY_{tj}\}$.
\end{enumerate}
\end{theorem}

The procedures involving $c^*$ and $c_{\rm stud}$ follow from Theorems~\ref{thm:cstar} and~\ref{thm:variants}, applied to the coordinates in $I_0$, and the two-sided uncentered test corresponds to the test described in \citet{ABR2010b}, where the authors note that, to the best of their knowledge, the validity of the uncentered approach has been proven only in the two-sided case \citep[\S1.3]{ABR2010b}. When $\mu_j<0$, the observations are not symmetric, and Lemma~\ref{lem:exact} cannot be applied to them directly. However, from $\bareps\le1$ and $\mu_j\le0$ it follows that $n^{-1}\sum_t\varepsilon_tY_{tj}=\bar y_j(\varepsilon)+\bareps\,\mu_j\ge\bar y_j(\varepsilon) +\mu_j$ for $j\in I_0$, and this inequality reduces the one-sided case to Lemma~\ref{lem:exact} (Appendix~\ref{app:mt}).

As long as the FWER is kept at $\alpha$, a procedure with lower thresholds rejects more false hypotheses and, consequently, has more power. Our task is to compare procedure (a) from Theorem~\ref{thm:fwer}, whose set of rejected hypotheses we denote by $R_{c^*}$, with two benchmark methods that use the same statistic $|\bar Y_j|$.

The first benchmark is the oracle procedure, which uses the true law: let $Q^*(C)$ denote the upper $\alpha$-quantile of $\max_{j\in C} |\bar y_j|$ and denote by $R_{\rm or}(\lambda)$ the set of hypotheses rejected by the step-down procedure with thresholds $\lambda\,Q^*(C)$. When $\lambda=1$, these are the ideal thresholds from Section~\ref{sec:notation} for the seminorms $\max_{j\in C}|v_j|$.

\begin{theorem}\label{thm:poweror}
Under the assumptions of Theorem~\ref{thm:effK}, for every $u\in(0,1)$,

\[ \Prob\big(R_{\rm or}(1+u)\subseteq R_{c^*}\subseteq R_{\rm or}(1-u)\big)\to1. \]

\end{theorem}

The second benchmark uses the threshold $c^*_j(\gamma)$, that is, $c^*$ at level $\gamma$, computed from the $j$th column alone using the same set $\mathcal E$ of sign patterns, and we denote by $R_B$ the set of hypotheses rejected by the step-down procedure with thresholds

\[ t_B(C)=\max_{j\in C}c^*_j(\alpha/|C|). \]

This is the Bonferroni correction to the one-dimensional thresholds from Section~\ref{sec:cstar}, in which the level $\alpha$ is split equally among the coordinates in $C$.

\begin{theorem}\label{thm:bonf}
For every data set and every nonempty set $C$, the threshold $c^*(C)$ of Theorem~\ref{thm:fwer}(a) satisfies $c^*(C)\le t_B(C)$, and therefore $R_B\subseteq R_{c^*}$.
\end{theorem}

Theorem~\ref{thm:bonf} essentially states that the joint threshold never rejects fewer hypotheses than the Bonferroni correction, and this holds for any data set, not just with high probability.

\section{Simulation study}\label{sec:numerics}

\begin{figure}[H]
\centering
\includegraphics[width=\textwidth]{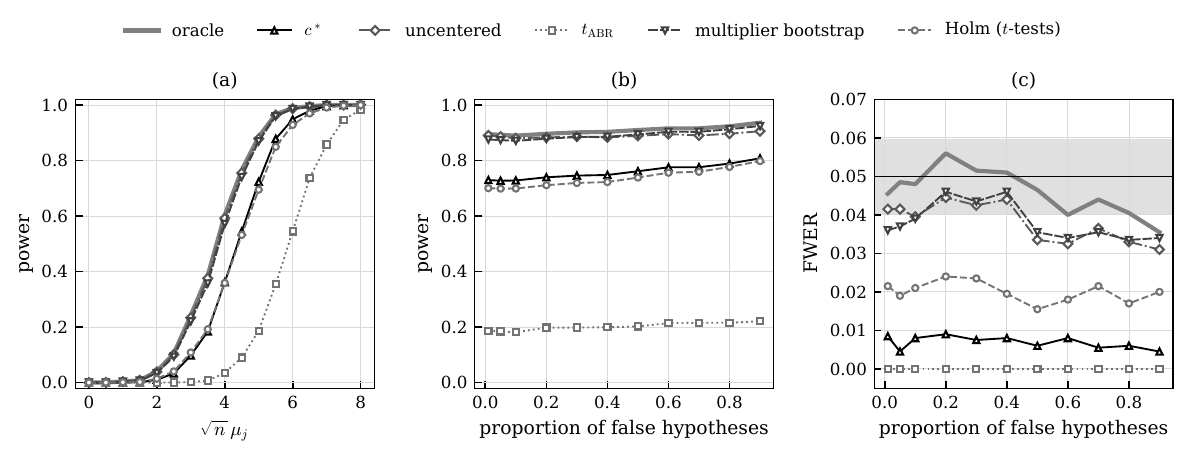}
\caption{Power and FWER of step-down procedures at $\alpha=0.05$ for $n=50$ and $K=1000$ for different effect sizes $\sqrt n\,\mu_j$ with $100$ false hypotheses (a) and for different proportions of false hypotheses at $\sqrt n\,\mu_j=5$ (b and c). The gray band in (c) marks $\alpha\pm2$ Monte Carlo standard errors.}\label{fig:mt}
\end{figure}

This section presents a small simulation study in which, for comparison, we also include two thresholds without a finite-sample guarantee under Assumption~\ref{ass:SA}, Bonferroni $t$-intervals, and the Gaussian multiplier bootstrap from \citet{CCK2013}. In all simulations, we use the maximum norm and $\alpha=0.05$, with correlation $0.5$ between the coordinates unless otherwise specified, and each threshold is computed from the identity pattern and $1999$ random sign patterns. Rejection probabilities and FWER are estimated from $2000$ samples ($500$ in Figure~\ref{fig:mt}(a)). In $t_{\rm ABR}$, we set $\delta=1/n$ and $\alpha_0=\alpha(1-1/n)$, as proposed in \citet[\S3.1]{ABR2010a}, with the Bonferroni bound $f$ computed from the true standard deviations. The (Python) code can be found at \url{https://github.com/0xKARTOD/signflip-confidence-regions} (see also \url{https://doi.org/10.5281/zenodo.24549573}).

\begin{figure}[H]
\centering
\includegraphics[width=\textwidth]{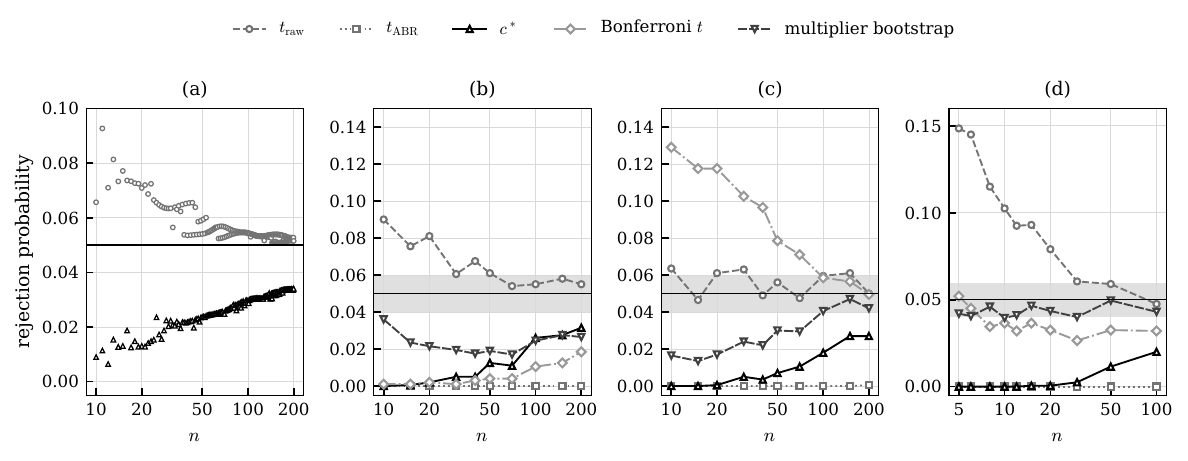}
\caption{Rejection probabilities at $\alpha=0.05$ (horizontal line), with gray bands indicating $\alpha\pm2$ Monte Carlo standard errors. (a): exact values for i.i.d.\ observations with values in $\{-1,0,1\}$. (b)--(d): simulated values with coordinates of correlation $0.5$, for Gaussian rows with $K=10$ (five of which have a standard deviation ten times larger), for centered exponential rows with $K=100$ and for i.i.d.\ Gaussian rows with $K=100$.}\label{fig:level}
\end{figure}

A general pattern can be observed across all four panels of Figure~\ref{fig:level}, where the rejection probability of $t_{\rm raw}$ exceeds $\alpha$ (at least for small values of $n$), whereas that of $c^*$ remains below $\alpha$, even in panel (c), where the coordinates are distributed asymmetrically and Assumption~\ref{ass:SA} does not hold. In the case of the lattice distribution of panel (a), the rejection probability of $t_{\rm raw}$ exceeds $\alpha$ even at $n=200$ (due to the occurrence of tied values at the quantile), and in panel (d), which most closely resembles the conditions examined by \citet{ABR2010a}, this effect persists up to $n=50$. Of course, this holds only for our designs, and does not contradict the findings of \citet[\S4]{ABR2010b}, whose simulations showed that, with $n=1000$, the raw quantile was nearly identical to the true quantile.

\begin{figure}[H]
\centering
\includegraphics[width=\textwidth]{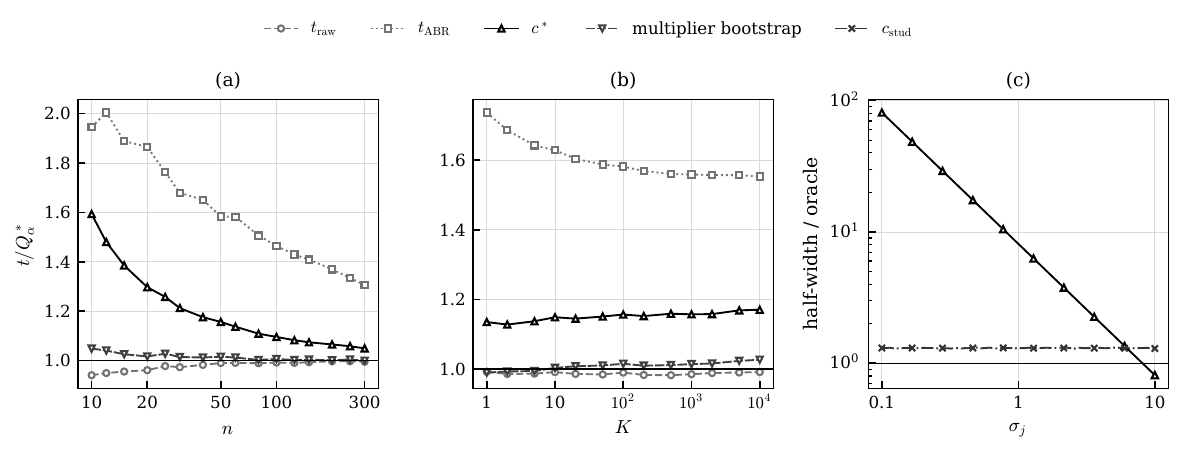}
\caption{The width of the regions at $\alpha=0.05$, see text. In (a) and (b), the median ratio of each threshold to $Q^*_\alpha$ depending on $n$ at $K=100$ and on $K$ at $n=50$, while in (c) the median half-width of the intervals from $c^*$ and $c_{\rm stud}$ relative to $\sigma_j\tilde Q_\alpha$ ($\tilde Q_\alpha$ is the upper $\alpha$-quantile of $\max_j|\bar y_j|/\sigma_j$) at $n=50$, $K=100$, and standard deviations $\sigma_j$ from $0.1$ to $10$.}\label{fig:width}
\end{figure}

The main result we can observe in Figure~\ref{fig:width} is that the regions of $c^*$ are wider than the ideal ones for small values of $n$. However, the ratio $c^*/Q^* _\alpha$ decreases as $n$ increases, as predicted by Theorems~\ref{thm:effK} and~\ref{thm:effHD}, while $t_{\rm ABR}$ turns out to be larger for any $n$ and $K$, including when $\delta=0.1$ is chosen as in \citet[\S4.1]{ABR2010b}. The second observation is that the raw quantile and the bootstrap remain close to $Q^*_\alpha$ (but without a finite-sample guarantee). In panel (c), the intervals constructed using $c_{\rm stud}$ correspond to the scale of each coordinate, whereas the intervals constructed using $c^*$ have a common half-width determined by the coordinates with the largest variances and turn out to be too wide where $\sigma_j$ is small.

For our next scenario (Figure~\ref{fig:mt}), we set $n=50$ and $K=1000$, assign the first coordinates a common mean value of $\mu_j>0$, and measure the power as the proportion of false hypotheses rejected by the procedure, averaged across samples. In addition to the step-down procedures from Section~\ref{sec:stepdown} with $c^*$, the uncentered threshold, and $t_{\rm ABR}$, we include the “oracle” step-down procedure, the procedure from \citet{RomanoWolf2005} with bootstrap thresholds, and the procedure from \citet{Holm1979} with $t$-tests.

Thus, our simulation study shows that all procedures control the FWER up to Monte Carlo error, and the uncentered approach and bootstrap are almost as powerful as the oracle. The procedure based on $c^*$ loses some power but rejects about as many false hypotheses as the procedure of \citet{Holm1979}, while $t_{\rm ABR}$ has the lowest power. We link the conservativeness of $c^*$ and $t_{\rm ABR}$ to the worst-case bound on the term $\bareps\,\bar y$ in \eqref{eq:id}, which is not required for the uncentered test.

\section{Discussion}\label{sec:discussion}

We have concluded that the remainder of \citet{ABR2010a} matters only in finite samples. For i.i.d.\ rows with a fixed $K$ and a distribution with an absolutely continuous component, $t_{\rm raw}$ is asymptotically valid and $c^*/t_{\rm raw}=1+\sqrt{2/(\pi n)}+o_P(n^{-1/2})$, that is, there is a difference between the two thresholds of the order of magnitude $n^{-1/2}$. However, for small $n$ or for observations that are not identically distributed, $t_{\rm raw}$ may yield a rejection probability well above $\alpha$ (Section~\ref{sec:size}), while $c^*$ maintains the confidence level.

Of course, $c^*$ is conservative, since it controls the level nonasymptotically, that is, for every $n$. In particular, when $K\ge2$, the rejection probability of $c^*$ turns out to be lower than $\alpha$ by more than the value determined by the discreteness of the sign-flip law, mainly because the proof uses only the fact that $\bar y$ lies in $\mathcal C(c^*)$, while this region contains points outside $C_{\rm ex}$. A region that does not contain $C_{\rm ex}$, for example, one obtained by splitting $\alpha$ between the exact regions for the max-norm and the Euclidean norm, may be less conservative, but we do not know how to compute it as quickly as $c^*$. The overall conclusion of the simulations in Section~\ref{sec:numerics} is that $c^*$ is suitable for coordinates with comparable scales, while $c_{\rm stud}$ is suitable for coordinates with different scales. Furthermore, when testing the hypothesis $H_j:\mu_j=0$, the uncentered step-down test is almost as powerful as the oracle procedure. Finally, Theorem~\ref{thm:imposs} showed that at a level below $2^{1-n}$, that is, when $n\le5$ and $\alpha=0.05$, and under Assumption~\ref{ass:SA} alone, there does not exist a translation-equivariant region of finite size.

Further study is needed because some questions remain open. Although the rejection probability of $c^*$ tends to $\alpha$ when $\log K=o(n^{1/3})$, as is the case with $t_{\rm raw}$ in the work of \citet{CCKK2022}, the ratio $c^*/Q^*_\alpha$ is controlled only when $\log K=o(n^{1/5})$. It would be interesting to know if the sign-flip law of $\sqrt n\|\bar y\|_\infty$ has a density bounded away from zero near its upper $\alpha$-quantile without a Gaussian approximation, since this would extend control over the ratio to the case of $\log K=o(n^{1/3})$. Likewise, Theorems~\ref{thm:sharp} and~\ref{thm:effS} concern only a fixed $K$, and the question whether the term $\sqrt{2/(\pi n)}$ remains unchanged when $K$ grows with $n$ thus remains open. Regarding the sign-perturbed sums from the work of \citet{CsajiCampiWeyer2015} in linear regression, we expect that Lemma~\ref{lem:closed} will retain its form, and the mean value will be replaced by a vector of coefficients so that their outer approximation could be computed using a single sorting operation rather than semidefinite programming.

\section*{Statement on AI use}

In preparing this article, we used the AI assistants Claude Opus~5.5 and Claude Sonnet~5.5 (Anthropic). All results and proofs were obtained by us, but the assistants were used to check our arguments for errors and to search the literature for results that could be cited rather than re-proven. They also helped to correct grammatical and stylistic errors in the text, to verify that the code for the numerical experiments calculates the threshold values specified in the article, and to format the figures.

\bibliographystyle{plainnat}
\bibliography{refs}

\clearpage
\appendix
\begin{center}
{\large\bfseries SUPPLEMENT TO ``CONFIDENCE REGIONS FOR A MEAN FROM\\[2pt] SIGN-FLIPPED RESIDUALS''}\\[12pt]
{\scshape By Roman Zinovev}\\[4pt]
{\itshape Nanyang Technological University, 50 Nanyang Avenue, Singapore 639798}
\end{center}

\section{Proofs for Section~\ref{sec:size}}\label{app:size}

In the proofs for Section~\ref{sec:size}, we will first need to reduce the sign-symmetric data to random signs multiplied by fixed values, and then count the combinations of signs for which $|\bar e(\varepsilon)|$ is not less than $|\bar y|$. Throughout this section, we assume $K=1$, $\phi=|\cdot|$, and $\mathcal E=\{\pm1\}^n$, so that $N=2^n$ and $k_\alpha=\lfloor\alpha2^n\rfloor$, and we write $\sigma\cdot a=(\sigma_1a_1,\dots,\sigma_na_n)$ for $\sigma\in\{\pm1\}^n$ and $a\in\R^n$. Definition \eqref{eq:q} makes $t_{\rm raw}$ the $(k_\alpha+1)$-th largest value among the $2^n$ values of $|\bar e(\varepsilon)|$, so the event $|\bar y|>t_{\rm raw}$ occurs if and only if no more than $k_\alpha$ of these values are not less than $|\bar y|$:

\begin{equation}\label{eq:rawrej}
|\bar y|>t_{\rm raw}\iff\#\big\{\varepsilon\in\{\pm1\}^n:\ |\bar e(\varepsilon)|\ge|\bar y|\big\}\le k_\alpha .
\end{equation}

The proofs presented below rely on two estimates for the sum $D=\sum_{i=1}^{n'}\varepsilon_i$ of $n'\ge1$ independent Rademacher variables, i.e., $\Prob(\varepsilon_i=1)=\Prob(\varepsilon_i=-1)=1/2$. The first is Hoeffding’s inequality \citep[Theorem 2]{Hoeffding1963} for variables with values in $[-1,1]$, which gives $\Prob(D\ge u)\le\exp(-u^2/(2n'))$ for $u>0$. Applied to $-D$, which has the same distribution as $D$, the inequality yields

\begin{equation}\label{eq:hoeff}
\Prob(|D|\ge u)\le2\exp\Big(-\frac{u^2}{2n'}\Big),\qquad u>0.
\end{equation}

The second bound concerns the point probabilities of $D$, which we express in terms of $c_k=\binom{2k}k4^{-k}$, $k\ge0$. We use Stirling’s formula with error estimates from the work of \citet[equations (1)--(2)]{Robbins1955} in the form

\[ \sqrt{2\pi}\,k^{k+\frac12}e^{-k}e^{1/(12k+1)}<k!<\sqrt{2\pi}\,k^{k+\frac12}e^{-k}e^{1/(12k)},\qquad k\ge1, \]

which, applied to $(2k)!$ and $k!$ in $c_k=(2k)!/(4^k(k!)^2)$, gives, for $k\ge1$,

\[ c_k<\frac1{\sqrt{\pi k}}\exp\Big(\frac1{24k}-\frac2{12k+1}\Big)<\frac1{\sqrt{\pi k}},\qquad c_k>\frac1{\sqrt{\pi k}}\exp\Big(\frac1{24k+1}-\frac1{6k}\Big)>\frac{e^{-1/6}}{\sqrt{\pi k}}. \]

For the point probabilities of $D$ this gives

\begin{align*}
n'&=2k: && \max_j\Prob(D=j)=\Prob(D=0)=c_k<(\pi n'/2)^{-1/2},\\
n'&=2k+1,\ k\ge1: && \max_j\Prob(D=j)=\Prob(D=1)=c_k\tfrac{2k+1}{2k+2}<(\pi k)^{-1/2}\le n'^{-1/2},\\
& && \Prob(|D|\le1)=2\Prob(D=1)=c_k\tfrac{2k+1}{k+1}\ge c_k,
\end{align*}

where $(\pi k)^{-1/2}\le n'^{-1/2}$ because $\pi k\ge2k+1$. Combining $k\le n'/2$, the lower bound for $c_k$, and the inequality $e^{-1/6}(2/\pi)^{1/2}>1/2$, as well as the case $n'=1$, where $\Prob(D=\pm1)=1/2$, we have

\begin{equation}\label{eq:binom}
\max_j\Prob(D=j)\le n'^{-1/2},\qquad \Prob(|D|\le1)\ge\tfrac12\,n'^{-1/2}.
\end{equation}

\subsection*{Technical lemma}

Lemma~\ref{lem:orbit} below is a version of the lemma of \citet[p.~1283]{Efron1969}, which shows that under orthant symmetry, the distribution of the Student's $t$-statistic is a mixture of generalized binomial distributions. Initially, Efron assumes that $y=0$ has probability zero \citep[p. 1278, footnote 2]{Efron1969}, and we adapt his argument by assuming that the uniform sign vector $\sigma$ is chosen independently of $y$, which also covers the zero observations that we consider in Section~\ref{sec:midp}.

\begin{lemma}\label{lem:orbit}
Let $K=1$, let $y$ satisfy Assumption~\ref{ass:SA}, and let $\sigma$ be uniform on $\{\pm1\}^n$ and independent of $y$. Then $y\eqd\sigma\cdot|y|$, and for every measurable $A\subseteq\R^n$

\[ \sup\big\{\Prob(y\in A):\ y\text{ satisfies Assumption~\ref{ass:SA}}\big\}=\sup_{a\in[0,\infty)^n}\Prob(\sigma\cdot a\in A). \]

\end{lemma}

\begin{proof}

First, we fix $a\in[0,\infty)^n$ and let $\varepsilon\in\{\pm1\}^n$, then $\varepsilon\cdot(\sigma\cdot a)=(\varepsilon\cdot\sigma)\cdot a$, and $\varepsilon\cdot\sigma\eqd\sigma$, since multiplication by $\varepsilon$ is a bijection of the set $\{\pm1\}^n$. It follows from this that $\varepsilon\cdot(\sigma\cdot a)\eqd\sigma\cdot a$, that is, the law of $\sigma\cdot a$ satisfies Assumption~\ref{ass:SA} and is one of the laws in the left-hand side. Therefore

\[ \sup\big\{\Prob(y\in A):\ y\text{ satisfies Assumption~\ref{ass:SA}}\big\}\ge\sup_{a\in[0,\infty)^n}\Prob(\sigma\cdot a\in A). \]

Conversely, consider $y$ satisfying Assumption~\ref{ass:SA}, and write $y=s\cdot|y|$, where $s_t=\mathrm{sign}(y_t)$, and $\mathrm{sign}(0)=1$. Taking into account $\sigma$, which is independent of $y$ and takes each value with probability $2^{-n}$, and applying Assumption~\ref{ass:SA} to each term, we obtain that for any measurable $B\subseteq\R^n$

\[ \Prob(\sigma\cdot y\in B)=\sum_{\varepsilon\in\{\pm1\}^n}\Prob(\sigma=\varepsilon)\,\Prob(\varepsilon\cdot y\in B)=2^{-n}\sum_{\varepsilon\in\{\pm1\}^n}\Prob(y\in B)=\Prob(y\in B), \]

which shows that $y\eqd\sigma\cdot y=(\sigma\cdot s)\cdot|y|$. Furthermore, given that $s_t^2=1$, the event $\{\sigma\cdot s=v\}$ is equal to the event $\{\sigma=v\cdot s\}$ for any $v\in\{\pm1\}^n$, and provided that $y$ is known, the vector $v\cdot s$ is fixed, while $\sigma$ remains uniformly distributed. Hence, for a given $|y|$, each combination of signs $v$ has a conditional probability of $2^{-n}$, i.e., for any measurable set $C\subseteq[0,\infty)^n$,

\[ \Prob\big(\sigma\cdot s=v,\ |y|\in C\big)=\E\Big[\mathbf 1\{|y|\in C\}\,\Prob\big(\sigma=v\cdot s\,\big|\,y\big)\Big]=2^{-n}\,\Prob(|y|\in C). \]

Consequently, $\sigma'=\sigma\cdot s$ is uniformly distributed on $\{\pm1\}^n$ and independent of $|y|$, and $y\eqd\sigma'\cdot|y|\eqd\sigma\cdot|y|$. Finally, $\sigma$ and $|y|$ are independent, so by Fubini’s theorem on the product of their distributions (where $\Prob_{|y|}$ is the distribution of $|y|$)

\[ \Prob(y\in A)=\Prob(\sigma\cdot|y|\in A)=\int_{[0,\infty)^n}\Prob(\sigma\cdot a\in A)\,d\Prob_{|y|}(a)\le\sup_{a\in[0,\infty)^n}\Prob(\sigma\cdot a\in A). \]

Taking the supremum over all $y$ that satisfy Assumption~\ref{ass:SA} we get the reverse inequality, which, together with the first inequality, gives

\[ \sup\big\{\Prob(y\in A):\ y\text{ satisfies Assumption~\ref{ass:SA}}\big\}=\sup_{a\in[0,\infty)^n}\Prob(\sigma\cdot a\in A). \qedhere \]

\end{proof}

\subsection*{Proof of Theorem~\ref{thm:midp}}

\begin{proof}[]

Throughout the proof, we assume that $n>m$, $y=\sigma\cdot(b,0,\dots,0)$ and $S=\sum_{t\le m}\sigma_tb_t$, so that $\bar y=S/n$. We also use the following notation, for $\varepsilon\in\{\pm1\}^n$,

\[ \tau_t=\varepsilon_t\sigma_t\ (t\le m),\qquad S(\tau)=\sum_{t\le m}\tau_tb_t,\qquad g=\sum_{t\le m}\varepsilon_t=\sum_{t\le m}\tau_t\sigma_t,\qquad D=\sum_{t>m}\varepsilon_t. \]

Then $n\bar y(\varepsilon)=\sum_t\varepsilon_ty_t=S(\tau)$ and $n\bareps=g+D$, so by \eqref{eq:id}

\begin{equation}\label{eq:nebar}
n\bar e(\varepsilon)=n\bar y(\varepsilon)-n\bareps\,\bar y=S(\tau)-r,\qquad r=\frac{(g+D)S}n.
\end{equation}

The observations with $t>m$ are zero, so the event $\{|\bar y|>t_{\rm raw}\}$ depends only on $(\sigma_1,\dots,\sigma_m)$, which we fix, and we set $x=|S|$. If $\varepsilon$ is uniformly distributed on $\{\pm1\}^n$, then the blocks $(\varepsilon_1,\dots,\varepsilon_m)$ and $(\varepsilon_{m+1},\dots,\varepsilon_n)$ are independent, and $\tau$ is a bijective function of the first block, so $\tau$ is uniformly distributed on $\{\pm1\}^m$ and is independent of $D$, while $g$ is a function of $\tau$ with $|g|\le m$. In view of \eqref{eq:nebar}, the inequality $|\bar e(\varepsilon)|\ge|\bar y|$ is equivalent to the inequality $|S(\tau)-r|\ge x$, and criterion \eqref{eq:rawrej} takes the following form

\begin{equation}\label{eq:pn}
|\bar y|>t_{\rm raw}\iff\#\big\{\varepsilon:\ |S(\tau)-r|\ge x\big\}\le k_\alpha\iff p_n\le k_\alpha2^{-n},
\end{equation}

where, with $\varepsilon$ uniform on $\{\pm1\}^n$,

\[ p_n=2^{-n}\#\big\{\varepsilon:\ |S(\tau)-r|\ge x\big\}=\Prob_\varepsilon\big(|S(\tau)-r|\ge x\big)=2^{-m}\sum_{\tau\in\{\pm1\}^m}\Prob_D\big(|S(\tau)-r|\ge x\big). \]

Here, we denote $\Prob_\varepsilon$ as the uniform distribution of $\varepsilon$ on $\{\pm1\}^n$ with fixed data, and $\Prob_D$ as the distribution law of $D=\sum_{t>m}\varepsilon_t$ under $\Prob_\varepsilon$. The independence of $\tau$ and $D$ makes $\Prob_D(\,\cdot\,)$ with a fixed $\tau$ a conditional probability given $\tau$. As in Section~\ref{sec:midp}, $S'$ denotes an independent copy of $S$.

If $x=0$, then

\[ r=0 \text{ and } p_n=1>k_\alpha2^{-n} \implies |\bar y|\le t_{\rm raw}. \]

Also

\[ F_b(0)=\Prob(|S'|>0)+\tfrac12\Prob(S'=0)\ge\tfrac12>\alpha, \]

and the event $\{F_b(|S|)\le\alpha\}$ also does not contain these values of $(\sigma_1,\dots,\sigma_m)$.

If $x>0$, we claim that each term in \eqref{eq:pn} converges to $\mathbf 1\{|S(\tau)|>x\}$ when $|S(\tau)|\ne x$, and to $1/2$ when $|S(\tau)|=x$, and we prove this separately for each case.

\emph{Case $|S(\tau)|\ne x$.} Let us denote $\Delta=\big||S(\tau)|-x\big|>0$ and note that, according to \eqref{eq:nebar}, $|r|\le(|g|+|D|)x/n\le(m+|D|)x/n$, so $|r|<\Delta$ provided that $|D|<n\Delta/x-m$. At the same time, the triangle inequality gives

\[ |S(\tau)|\ge x+\Delta\ \Rightarrow\ |S(\tau)-r|\ge|S(\tau)|-|r|>x \]

in the case $|S(\tau)|>x$, and

\[ |S(\tau)|\le x-\Delta\ \Rightarrow\ |S(\tau)-r|\le|S(\tau)|+|r|<x \]

in the case $|S(\tau)|<x$, so, on this event,

\[ \mathbf 1\{|S(\tau)-r|\ge x\}=\mathbf 1\{|S(\tau)|>x\}. \]

Using \eqref{eq:hoeff} with $n'=n-m$, for $n>mx/\Delta$,

\[ \Big|\Prob_D\big(|S(\tau)-r|\ge x\big)-\mathbf 1\{|S(\tau)|>x\}\Big|\le\Prob\Big(|D|\ge\frac{n\Delta}x-m\Big)\le2\exp\Big(-\frac{(n\Delta/x-m)^2}{2(n-m)}\Big)\longrightarrow0. \]

\emph{Case $|S(\tau)|=x$.} Then $S(\tau)=sS$ with $s=S(\tau)/S\in\{\pm1\}$, and by \eqref{eq:nebar}

\[ |S(\tau)-r|=\Big|sS-\frac{(g+D)S}n\Big|=x\,\Big|1-\frac{s(g+D)}n\Big|. \]

The number $1-s(g+D)/n$ lies in $[0,2]$ because $|g+D|\le m+(n-m)=n$, so

\[ |S(\tau)-r|\ge x\iff1-\frac{s(g+D)}n\ge1\iff sD\le-sg. \]

Here, $g$ is fixed, $|sg|\le m$, and $sD\eqd D$ because $D$ is symmetric. Symmetry guarantees that $\Prob(D\le c)=1-\Prob(D<-c)$ for every integer $c$, and hence $|2\Prob(D\le c)-1|=|\Prob(D\le c)-\Prob(D<-c)|\le\Prob(|D|\le|c|)$. Taking $c = -sg$ and using \eqref{eq:binom}, we obtain

\[ \Big|\Prob_D\big(|S(\tau)-r|\ge x\big)-\frac12\Big|\le\frac12\Prob(|D|\le m)\le\frac{2m+1}{2\sqrt{n-m}}\longrightarrow0. \]

\emph{Completion of the proof.} Summing the two cases over $\tau$ in \eqref{eq:pn} and using that $S(\tau)$ with $\tau$ uniform has the law of $S'$, we obtain

\[ p_n\longrightarrow2^{-m}\sum_{\tau\in\{\pm1\}^m}\Big(\mathbf 1\{|S(\tau)|>x\}+\tfrac12\mathbf 1\{|S(\tau)|=x\}\Big)=\Prob(|S'|>x)+\tfrac12\Prob(|S'|=x)=F_b(x). \]

From $\lfloor\alpha2^n\rfloor>\alpha2^n-1$ we get $k_\alpha2^{-n}\in(\alpha-2^{-n},\alpha]$ and $k_\alpha2^{-n}\to\alpha$.

If $F_b(x)<\alpha$, then 

\[ p_n-k_\alpha2^{-n}\to F_b(x)-\alpha<0\]

and $|\bar y|>t_{\rm raw}$ for all large $n$. Similarly, if $F_b(x)>\alpha$, then $|\bar y|\le t_{\rm raw}$ for all large $n$. There are finitely many values $(\sigma_1, \dots, \sigma_m)$, so there exists $n_1$ such that for $n \ge n_1$

\[ \mathbf 1\{F_b(|S|)<\alpha\}\le\mathbf 1\{|\bar y|>t_{\rm raw}\}\le\mathbf 1\{F_b(|S|)\le\alpha\} \]

for every value of $(\sigma_1,\dots,\sigma_m)$. Taking expectations over the uniform $(\sigma_1,\dots,\sigma_m)$, for which $S\eqd S'$, we get

\[ \Prob\big(F_b(|S|)<\alpha\big)\le\Prob(|\bar y|>t_{\rm raw})\le\Prob\big(F_b(|S|)\le\alpha\big),\qquad n\ge n_1, \]

which implies the two middle inequalities of the theorem. If $\Prob(F_b(|S|)=\alpha)=0$, then the two outer probabilities coincide, and the rejection probability is equal to $\Prob(F_b(|S|)<\alpha)$ for $n\ge n_1$.

It remains to prove the final inequality of the theorem: $\Prob(F_b(|S|)\le\alpha)\le2\alpha$. To do so, we consider the function $p(x)=\Prob(|S'|\ge x)$, which is nonincreasing and left-continuous, with $p(x)=1$ when $x\le0$ and $p(x)=0$ when $x>\max|S'|$. Let $u\in(0,1)$ and denote $x_u=\inf\{x:p(x)\le u\}\in[0,\infty)$. Due to monotonicity, the set $\{x:p(x)\le u\}$ is the interval $[x_u,\infty)$ or $(x_u,\infty)$.

In the first case

\[ \Prob\big(p(|S|)\le u\big)=\Prob(|S|\ge x_u)=p(x_u)\le u, \]

and in the second

\[ \Prob\big(p(|S|)\le u\big)=\Prob(|S|>x_u)=\lim_{x\downarrow x_u}\Prob(|S|\ge x)=\lim_{x\downarrow x_u}p(x)\le u. \]

Since $F_b(x)=\tfrac12p(x)+\tfrac12\Prob(|S'|>x)\ge\tfrac12p(x)$ and $2\alpha<1$,

\[ \Prob\big(F_b(|S|)\le\alpha\big)\le\Prob\big(p(|S|)\le2\alpha\big)\le2\alpha. \qedhere \]

\end{proof}

\subsection*{Proof of Theorem~\ref{thm:attain}}

\begin{proof}[]

In this proof we use the same notation as in the proof of Theorem~\ref{thm:midp}, and in addition

\[ b=\mathbf 1_m,\qquad S=\sum_{t\le m}\sigma_t,\qquad S(\tau)=\sum_{t\le m}\tau_t,\qquad x=|S|\in\{m,m-2,\dots\}. \]

As $\alpha<1/2$, we have $m\ge2$, and the definition of $m$ gives $2^{-m}\le\alpha<2^{1-m}$. Consider, for $n>m^2/2$, the event $E_n$ and the number $\pi_n$ given by

\[ E_n=\Big\{|D|<\frac{2n}m-m\Big\},\qquad \pi_n=2\exp\Big(-\frac{(2n/m-m)^2}{2(n-m)}\Big), \]

so that $\Prob(E_n^c)\le\pi_n$ by \eqref{eq:hoeff}. As in the previous proof, $|\bar y|\le t_{\rm raw}$ when $x=0$, so it remains to consider the case $x\ne0$.

\emph{Case $x=m$.} Here, $\sigma_1=\dots=\sigma_m$, and $|S(\tau)|=m$ holds only when $\tau=\pm\mathbf 1_m$. For these two patterns, $S(\tau)=sS$ and $g=\sum_t\tau_t\sigma_t=sm$, therefore, according to the previous proof for the case $|S(\tau)|=x$, the inequality $|S(\tau)-r|\ge x$ holds if and only if $sD\le-m$, that is, with probability

\[ \Prob(D\ge m)=\tfrac12\big(1-\Prob(|D|\le m-1)\big)\le\tfrac12\big(1-\Prob(|D|\le1)\big)\stackrel{\eqref{eq:binom}}{\le}\tfrac12\big(1-\tfrac12(n-m)^{-1/2}\big). \]

For every other $\tau$, it holds that $|S(\tau)|\le m-2$, and since $|r|=|g+D|\,m/n$ and $|g|\le m$,

\[ |S(\tau)-r|\ge m\ \Rightarrow\ |r|\ge m-|S(\tau)|\ge2\ \Rightarrow\ |g+D|\ge\frac{2n}m\ \Rightarrow\ |D|\ge\frac{2n}m-m. \]

Thus, the inequality $|S(\tau)-r|\ge x$ holds for these $\tau$ only on $E_n^c$, and by combining these two estimates above with \eqref{eq:pn}, we observe

\[ p_n\le2\cdot2^{-m}\Prob(D\ge m)+\Prob(E_n^c)\le2^{-m}-2^{-m-1}(n-m)^{-1/2}+\pi_n. \]

Due to the fact that $k_\alpha2^{-n}>\alpha-2^{-n}\ge2^{-m}-2^{-n}$, the event $|\bar y|>t_{\rm raw}$ occurs whenever

\begin{equation}\label{eq:top}
\pi_n+2^{-n}\le2^{-m-1}(n-m)^{-1/2}.
\end{equation}

\emph{Case $0<x\le m-2$.} The patterns $\tau=\pm\mathbf 1_m$ have $|S(\tau)|=m\ge x+2$, and on $E_n$

\[ |r|\le\frac{(|g|+|D|)x}n<\frac{(m+2n/m-m)x}n=\frac{2x}m<2, \]

so both patterns satisfy the condition $|S(\tau)-r|\ge m-|r|>m-2\ge x$ on $E_n$, which in \eqref{eq:pn} gives $p_n\ge2\cdot2^{-m}\Prob(E_n)\ge2^{1-m}(1-\pi_n)$. This together with $k_\alpha2^{-n}\le\alpha$ ensures that $|\bar y|\le t_{\rm raw}$ when

\begin{equation}\label{eq:low}
2^{1-m}(1-\pi_n)>\alpha.
\end{equation}

As $n\to\infty$, 

\[ (\pi_n+2^{-n})(n-m)^{1/2}\to0, \]

\[ 2^{1-m}(1-\pi_n)\to2^{1-m}>\alpha, \]

so, there exists $n_0(\alpha)>m^2/2$ such that the inequalities \eqref{eq:top} and \eqref{eq:low} hold for any $n\ge n_0(\alpha)$. Then the inequality $|\bar y|>t_{\rm raw}$ holds for these $n$ if and only if $x=m$, and

\[ \Prob(|\bar y|>t_{\rm raw})=\Prob(\sigma_1=\dots=\sigma_m)=2\cdot2^{-m}=2^{1-m}, \]

which proves the theorem, since $2^{1-m}=2\alpha$ exactly when $\alpha=2^{-m}$.

\end{proof}

\subsection*{Proof of Theorem~\ref{thm:gauss}}

If the observations are i.i.d.\ $N(\mu,\varsigma^2)$, then 

\[ Z=\sqrt n\,\bar y/\varsigma\sim N(0,1),\qquad R=\|e\|_2^2/\varsigma^2\sim\chi^2_{n-1},\qquad T=\sqrt n\,\bar y/s\sim t_{n-1} \]

as in \citep[Example~5.1.1 and eq.~(5.16)]{LehmannRomano2005}, with independent $Z$ and $R$. Thus $Z^2$ and $R$ are independent random variables following chi-squared distributions, and the quantity $Z^2/(Z^2+R)$ follows a beta distribution \citep[Problem~5.18]{LehmannRomano2005}:

\[ T^2/(T^2+n-1)=Z^2/(Z^2+R)\sim{\rm Beta}(\tfrac12,\tfrac{n-1}2). \]

\begin{proof}[]

The threshold $t_{\rm raw}$ is one of the values of $|\bar e(\varepsilon)|$, and the largest of them is achieved when $\varepsilon_t=\mathrm{sign}(e_t)$. According to the Cauchy--Bunyakovsky--Schwarz inequality and the definition $\|e\|_2^2=(n-1)s^2$,

\[ t_{\rm raw}\le\max_\varepsilon|\bar e(\varepsilon)|=\frac1n\sum_t|e_t|\le\frac{\|e\|_2}{\sqrt n}=\sqrt{\frac{n-1}n}\;s. \]

Hence

\[ \{|T|>\sqrt{n-1}\}=\Big\{|\bar y|>\sqrt{\tfrac{n-1}n}\;s\Big\}\subseteq\{|\bar y|>t_{\rm raw}\}\ \Rightarrow\ \Prob(|\bar y|>t_{\rm raw})\ge\Prob\big(|T|>\sqrt{n-1}\big). \]

With $a=(n-1)/2$, the variable $X=T^2/(T^2+n-1)\sim{\rm Beta}(\tfrac12,a)$ exceeds $1/2$ exactly when $|T|>\sqrt{n-1}$, so that

\[ \Prob\big(|T|>\sqrt{n-1}\big)=\frac1{B(\frac12,a)}\int_{1/2}^1x^{-1/2}(1-x)^{a-1}\,dx\ge\frac1{B(\frac12,a)}\int_{1/2}^1(1-x)^{a-1}\,dx=\frac{2^{-a}}{a\,B(\frac12,a)}. \]

Inequality (7) of \citet{Wendel1948} states $(x/(x+\lambda))^{1-\lambda}\le\Gamma(x+\lambda)/(x^\lambda\Gamma(x))\le1$ for $x>0$ and $0<\lambda<1$. Taking here $x=a$ and $\lambda=1/2$, we get $\Gamma(a+\frac12)/\Gamma(a)\ge a/\sqrt{a+\frac12}$, so

\[ B\big(\tfrac12,a\big)=\frac{\Gamma(\frac12)\Gamma(a)}{\Gamma(a+\frac12)}\le\frac{\sqrt{\pi(a+\frac12)}}a,\qquad \Prob\big(|T|>\sqrt{n-1}\big)\ge\frac{2^{-a}}{\sqrt{\pi(a+\frac12)}}=\frac{2^{-(n-1)/2}}{\sqrt{\pi n/2}}. \]

The right-hand side exceeds $2^{1-n}=4^{-a}$ if and only if $h(a)=4^a-\pi(a+\frac12)>0$, and it is easy to check that

\[ h\big(\tfrac32\big)=8-2\pi>0\qquad\text{and}\qquad h'(a)=4^a\log4-\pi>0\quad\text{for }a\ge\tfrac32. \]

Thus $h$ is positive on $[\tfrac32,\infty)$, and

\[ \Prob\big(|T|>\sqrt{n-1}\big)>2^{1-n}\qquad\text{for }n\ge4. \]

The case $n=3$ follows from $\Prob(|t_2|>u)=1-u/\sqrt{2+u^2}$, which gives

\[ \Prob\big(|t_2|>\sqrt2\big)=1-\frac1{\sqrt2}>\frac14=2^{1-3}=2^{-2}. \]

\end{proof}

\section{Proofs for Sections~\ref{sec:cstar} and~\ref{sec:below}}\label{app:cstar}

In the proofs in this section, we will first need to establish two facts about \eqref{eq:q}, because most of the proofs given below reduce the statement about the resampled quantile to the counting of sign patterns. If $v,v':\mathcal E\to\R$ satisfy the condition $v\le v'$ pointwise, then $\#\{\varepsilon:v(\varepsilon)>x\}\le\#\{\varepsilon:v'(\varepsilon)>x\}$ for any $x$, that is, the infimum in \eqref{eq:q} is taken over a larger set for $v$ than for $v'$, and

\begin{equation}\label{eq:qmono}
v\le v'\ \Rightarrow\ q_\gamma\{v\}\le q_\gamma\{v'\}.
\end{equation}

Also, by definition, $q_\gamma\{v\}$ is the $(k_\gamma+1)$-th largest value among the $N$ values of $v(\varepsilon)$, so for any $w\in\R$

\begin{equation}\label{eq:qchar}
q_\gamma\{v\}\ge w\iff\#\{\varepsilon\in\mathcal E:\ v(\varepsilon)\ge w\}\ge k_\gamma+1.
\end{equation}

These facts will be useful in the next proofs.

\subsection*{Proof of Lemma~\ref{lem:closed}}

\begin{proof}[]

Using $ \mathcal W =\{w\ge0:\ w\le q_\alpha\{E(\varepsilon)+|\bareps|\,w \}\}$ we have $c^*=\sup\mathcal W$ according to \eqref{eq:cstar}, and \eqref{eq:qchar} shows that

\[ w\in\mathcal W\iff G(w)\ge k_\alpha+1,\qquad G(w)=\#\big\{\varepsilon\in\mathcal E:\ E(\varepsilon)+|\bareps|\,w\ge w\big\},\qquad w\ge0. \]

Next, note that a pattern with $|\bareps|=1$ is $\varepsilon=\pm\mathbf 1$, so $\bar e(\varepsilon)=\pm n^{-1}\sum_te_t=0$ and $E(\varepsilon)=\phi(0)=0$. Its value is $E(\varepsilon)+|\bareps|\,w=w$, and it is counted in $G(w)$ for each $w$. A pattern with $|\bareps|<1$ satisfies

\[ E(\varepsilon)+|\bareps|\,w\ge w\iff w\le\frac{E(\varepsilon)}{1-|\bareps|}=r(\varepsilon). \]

Hence, for $\mathcal E_0=\{\varepsilon\in\mathcal E:\ |\bareps|<1\}$ (which contains $N-m_1$ elements, counted with multiplicity)

\begin{equation}\label{eq:Gw}
G(w)=m_1+\#\{\varepsilon\in\mathcal E_0:\ r(\varepsilon)\ge w\},\qquad w\ge0.
\end{equation}

It remains to find the supremum of $\{w\ge0:\ G(w)\ge k_\alpha+1\}$, and we separately consider two cases.

\emph{Case $k_\alpha<m_1$.} Here, \eqref{eq:Gw} gives $G(w)\ge m_1\ge k_\alpha+1$ for any $w\ge0$, so $\mathcal W=[0,\infty)$ and $c^*=+\infty$.

\emph{Case $k_\alpha\ge m_1$.} Here, the index $j=k_\alpha+1-m_1$ satisfies the condition $1\le j\le N-m_1$, since $k_\alpha\le\alpha N<N$, so there exists the $j$-th largest value $r_{(j)}$ of the function $r$ on $\mathcal E_0$. Moreover, it is nonnegative because $E\ge0$, and \eqref{eq:Gw} gives

\[ w\in\mathcal W\iff\#\{\varepsilon\in\mathcal E_0:\ r(\varepsilon)\ge w\}\ge j\iff w\le r_{(j)}, \]

so $\mathcal W=[0,r_{(j)}]$ and $c^*=r_{(j)}$.

\end{proof}

\subsection*{Proof of Theorem~\ref{thm:cstar}}

\begin{proof}[]

To show that $C_{\rm ex}\subseteq\mathcal C(c^*)$, it suffices to show that $\phi(\bar Y-x)\le c^*$ for every $x\in C_{\rm ex}$. To do so, fix any $x \in C_{\rm ex}$ and set $w = \phi(\bar Y - x)$. From the triangle inequality and the homogeneity of $\phi$, it follows that for any pattern

\[ \phi\big(\bar e(\varepsilon)+\bareps(\bar Y-x)\big)\le\phi(\bar e(\varepsilon))+|\bareps|\,\phi(\bar Y-x)=E(\varepsilon)+|\bareps|\,w, \]

and therefore

\[ w\stackrel{\eqref{eq:Cex}}{\le}q_\alpha\big\{\phi\big(\bar e(\varepsilon)+\bareps(\bar Y-x)\big)\big\}\stackrel{\eqref{eq:qmono}}{\le}q_\alpha\{E(\varepsilon)+|\bareps|\,w\}, \]

so $w$ belongs to the set from \eqref{eq:cstar} and $w\le c^*$, which, for an arbitrary $x\in C_{\rm ex}$, means $C_{\rm ex}\subseteq\mathcal C(c^*)$.

Next, we bound the rejection probability $\Prob(\phi(\bar y)>c^*)$. Recalling that $\bar Y-\mu=\bar y$, we take $x=\mu$ and, according to \eqref{eq:id}, obtain $\bar e(\varepsilon)+\bareps\,\bar y=\bar y(\varepsilon)$, that is

\[ \mu\notin C_{\rm ex}\iff\phi(\bar y)>q_\alpha\{\phi(\bar y(\varepsilon))\}\iff W(y)>q_\alpha\{W(\varepsilon\cdot y)\},\qquad W(y)=\phi(\bar y), \]

where we used the fact that $W(\varepsilon\cdot y)=\phi(n^{-1}\sum_t\varepsilon_ty_t)=\phi(\bar y(\varepsilon))$. Lemma~\ref{lem:exact} and the inclusion $C_{\rm ex}\subseteq\mathcal C(c^*)$ yield

\[ \Prob\big(\phi(\bar y)>c^*\big)=\Prob\big(\mu\notin\mathcal C(c^*)\big)\le\Prob(\mu\notin C_{\rm ex})\le k_\alpha/N\le\alpha. \]

Finally, in the full set $\{\pm1\}^n$, the patterns with $|\bareps|=1$ are $\pm\mathbf 1$, so $m_1=2$, and Lemma~\ref{lem:closed} shows that $c^*<\infty$ if and only if $k_\alpha\ge2$, where

\[ k_\alpha=\lfloor\alpha2^n\rfloor\ge2\iff\alpha2^n\ge2\iff\alpha\ge2^{1-n}. \qedhere \]

\end{proof}

\subsection*{Proof of Theorem~\ref{thm:variants}}

\begin{proof}[]

Below, we use the convention $0/0=0$, which applies only to columns with $\|y_{\cdot j}\|_2=0$ or $\|e_{\cdot j}\|_2=0$, where the numerator is also equal to zero. Let

\[ W(y)=\max_j\frac{\sqrt n\,|\bar y_j|}{\|y_{\cdot j}\|_2},\qquad V(\varepsilon)=\max_j\frac{\sqrt n\,|\bar e_j(\varepsilon)|}{\|e_{\cdot j}\|_2}. \]

Note that changing the sign does not affect $\|y_{\cdot j}\|_2$, so $W(\varepsilon\cdot y)=\max_j\sqrt n\,|\bar y_j(\varepsilon)|/\|y_{\cdot j}\|_2$, and Lemma~\ref{lem:exact} gives $\Prob(\Omega^c)\le\alpha$ for $\Omega=\{W(y)\le q_\alpha\{W(\varepsilon\cdot y)\}\}$. Then the identity $\sum_te_{tj}=0$ gives

\[ \|y_{\cdot j}\|_2^2=\sum_t(e_{tj}+\bar y_j)^2=\|e_{\cdot j}\|_2^2+n\bar y_j^2\ge\|e_{\cdot j}\|_2^2, \]

and by \eqref{eq:id}, for every $j$,

\[ \frac{\sqrt n\,|\bar y_j(\varepsilon)|}{\|y_{\cdot j}\|_2}\le\frac{\sqrt n\,|\bar e_j(\varepsilon)|}{\|y_{\cdot j}\|_2}+|\bareps|\,\frac{\sqrt n\,|\bar y_j|}{\|y_{\cdot j}\|_2}\le V(\varepsilon)+|\bareps|\,W(y). \]

This gives $W(\varepsilon\cdot y)\le V(\varepsilon)+|\bareps|\, W(y)$, and using \eqref{eq:qmono} (on $\Omega$), $w=W(y)$ satisfies

\[ w\le q_\alpha\{W(\varepsilon\cdot y)\}\le q_\alpha\{V(\varepsilon)+|\bareps|\,w\}, \]

so that $w\le c_{\rm stud}$, which proves the first claim.

For the second part, note that on $\Omega$, for $c=c_{\rm stud}<1$, the inequality $\sqrt n\,|\bar y_j|\le c\,\|y_{\cdot j}\|_2$ and the identity for $\|y_{\cdot j}\|_2^2$ for each $j$ give

\[ n\bar y_j^2\le c^2\big(\|e_{\cdot j}\|_2^2+n\bar y_j^2\big)\ \Rightarrow\ n\bar y_j^2(1-c^2)\le c^2\|e_{\cdot j}\|_2^2\ \Rightarrow\ |\bar Y_j-\mu_j|=|\bar y_j|\le\frac c{\sqrt{1-c^2}}\,\frac{\|e_{\cdot j}\|_2}{\sqrt n}. \]

Therefore, the event in which $c_{\rm stud}<1$ and this bound fails for some $j$ is contained in $\Omega^c$ and has a probability of no more than $\alpha$.

\end{proof}

\subsection*{Proof of Theorem~\ref{thm:imposs}}

\begin{proof}[]

\emph{(i).} Let us fix $e_0$ in $\mathbf 1^\perp=\{v\in\R^n:\sum_tv_t=0\}$, the space in which the residual vector always lies, and suppose, to the contrary, that $h(e_0)<\infty$. We will set $c=h(e_0)+1$ and consider

\[ y_t=\sigma_t(c+e_{0t}),\qquad t=1,\dots,n, \]

with $\sigma$ uniformly distributed on $\{\pm1\}^n$, which satisfies Assumption~\ref{ass:SA}, because $\varepsilon\cdot y=(\varepsilon\cdot\sigma)\cdot(c\mathbf 1+e_0)$ and $\varepsilon\cdot\sigma\eqd\sigma$. From $\sum_te_{0t}=0$ it follows that

\[ \sigma=\mathbf 1\ \Rightarrow\ \bar y=c,\ \ e=e_0,\qquad \sigma=-\mathbf 1\ \Rightarrow\ \bar y=-c,\ \ e=-e_0, \]

and in both cases, $h(e)=h(e_0)<c=|\bar y|$ (let us recall that $h(-e_0)=h(e_0)$). Therefore

\[ \Prob\big(|\bar y|>h(e)\big)\ge\Prob(\sigma=\pm\mathbf 1)=2^{1-n}>\alpha, \]

which contradicts the validity of $h$ and shows that $h(e_0)=+\infty$ for any $e_0\in\mathbf 1^\perp$.

\emph{(ii).} As in the previous part, suppose, contrary to (ii), that $[L,U]$ possesses the specified properties. Translation equivariance with $c=\bar Y$, applied to the residual vector $e=Y-\bar Y\mathbf 1$, yields

\[ L(Y)=\bar Y+L(e),\qquad U(Y)=\bar Y+U(e), \]

and with $\bar y=\bar Y-\mu$ and $\Lambda(e)=\max\{|L(e)|,|U(e)|\}<\infty$,

\[ \mu\notin[L(Y),U(Y)]\iff\bar y\notin[-U(e),-L(e)]\ \Leftarrow\ |\bar y|>\Lambda(e). \]

Let $g$ denote the standard normal density, and for $c>0$, let $Y_1,\dots,Y_n$ be i.i.d.\ with density $\tfrac12g(\cdot-c)+\tfrac12g(\cdot+c)$, which is continuous, positive everywhere, and symmetric with respect to $\mu=0$. We can write $Y_t=c\sigma_t+Z_t$ with $\sigma$ uniform on $\{\pm1\} ^ n$ and $Z_1,\dots,Z_n$ i.i.d.\ $N(0,1)$, and they are all independent. Then on $\{\sigma=\pm\mathbf 1\}$,

\[ \bar y=\pm c+\bar Z,\qquad e=e_Z=Z-\bar Z\mathbf 1, \]

and the distribution of $(\bar Z,e_Z)$ does not depend on $c$. Hence

\[ \Prob\big(\mu\notin[L(Y),U(Y)]\big)\ge2^{-n}\Big\{\Prob\big(|c+\bar Z|>\Lambda(e_Z)\big)+\Prob\big(|-c+\bar Z|>\Lambda(e_Z)\big)\Big\}. \]

Now, if $c\to\infty$, then both indicators $\mathbf 1\{|\pm c+\bar Z|>\Lambda(e_Z)\}$ pointwise tend to $1$, since $\Lambda(e_Z)<\infty$, so by the dominated convergence theorem, the right-hand side tends to $2^{1-n}$. The limit $2^{1-n}$ is greater than $\alpha$, so for sufficiently large $c$, $\Prob(\mu\notin[L(Y),U(Y)])$ exceeds $\alpha$, which contradicts the condition.

\end{proof}

\section{Proofs for Section~\ref{sec:abr}}\label{app:abr}

Throughout the proofs in this appendix, $\mathcal W=\{w\ge0:\ w\le q_\alpha\{E(\varepsilon)+|\bareps|\,w\}\}$ is the set from \eqref{eq:cstar}, so that $c^*=\sup\mathcal W$. For any $w\in\mathcal W$, the inequality $E(\varepsilon)\ge(1-|\bareps|)w$ holds for at least $k_\alpha+1$ patterns. To prove Lemma~\ref{lem:sandwich} and Theorem~\ref{thm:abr}, we discard patterns with a large $|\bareps|$ and compare the remaining ones with a quantile of $E(\varepsilon)$.

\begin{lemma}\label{lem:split}
Let $\beta,\beta'>0$ with $\beta+\beta'\le\alpha$, and let $\lambda\in[0,1)$ satisfy $\#\{\varepsilon\in\mathcal E:\ |\bareps|>\lambda\}\le k_\beta$. Then every $w\in\mathcal W$ satisfies

\[ (1-\lambda)\,w\le q_{\beta'}\{E(\varepsilon)\}. \]

\end{lemma}

\begin{proof}

Using the inequality $\lfloor a+b\rfloor\ge\lfloor a\rfloor+\lfloor b\rfloor$ for real $a,b$ and the fact that $\lfloor\cdot\rfloor$ is nondecreasing, we have

\[ k_\alpha=\lfloor\alpha N\rfloor\ge\lfloor\beta N+\beta'N\rfloor\ge\lfloor\beta N\rfloor+\lfloor\beta'N\rfloor=k_\beta+k_{\beta'}. \]

If $w\in\mathcal W$, then it follows from \eqref{eq:qchar} that at least $k_\alpha+1$ patterns satisfy the condition $E(\varepsilon)+|\bareps|\,w\ge w$. No more than $k_\beta$ of them have $|\bareps|>\lambda$, so at least $k_\alpha+1-k_\beta\ge k_{\beta'}+1$ of them have $|\bareps|\le\lambda$, and for them

\[ E(\varepsilon)\ge(1-|\bareps|)\,w\ge(1-\lambda)\,w. \]

Thus, at least $k_{\beta'}+1$ patterns satisfy $E(\varepsilon)\ge(1-\lambda)w$, which, together with \eqref{eq:qchar}, means that $q_{\beta'}\{E(\varepsilon)\}\ge(1-\lambda)w$.

\end{proof}

\subsection*{Proof of Lemma~\ref{lem:sandwich}}

\begin{proof}[]

\emph{Lower bound.} In the first step, we note that the value $w_0=q_\alpha\{E(\varepsilon)\}\ge0$ satisfies the condition $E(\varepsilon)+|\bareps|\,w_0\ge E(\varepsilon)$ for every pattern, and \eqref{eq:qmono} gives

\[ q_\alpha\{E(\varepsilon)+|\bareps|\,w_0\}\ge q_\alpha\{E(\varepsilon)\}=w_0, \]

so $w_0\in\mathcal W$, which gives 

\[ c^*\ge w_0=q_\alpha\{E(\varepsilon)\}. \]

\emph{Upper bound.} In this step, when $v=|\bareps|$, the set of $x$ in \eqref{eq:q} is closed (because $x\mapsto\#\{\varepsilon:\ |\bareps|>x\}$ is nonincreasing and right-continuous), so the infimum is attained, and $\lambda=q_\beta\{|\bareps|\}$ satisfies $\#\{\varepsilon:\ |\bareps|>\lambda\}\le k_\beta$. Only the $m_1 \le k_\beta$ patterns with $|\bareps|=1$ have the maximum possible value of $1$, so the $(k_\beta+1)$-th largest value $\lambda$ is $<1$. This justifies the application of Lemma~\ref{lem:split} with $\beta'=\alpha-\beta$, and since the resulting bound $(1-\lambda)w\le q_{\alpha-\beta}\{E(\varepsilon)\}$ holds for every $w \in \mathcal W$, it also holds for $c^*=\sup\mathcal W$, that is,

\[ c^*\le\frac{q_{\alpha-\beta}\{E(\varepsilon)\}}{1-q_\beta\{|\bareps|\}}. \]

\emph{The bound on $q_\beta\{|\bareps|\}$.} Here, for $\mathcal E=\{\pm1\}^n$ and $u_\beta=\sqrt{2\log(2/\beta)/n}$, the sum $n\bareps$ over the uniform pattern is the sum of $n$ independent Rademacher variables. Applying \eqref{eq:hoeff} with $n'=n$ and $u=nu_\beta$, we have

\[ 2^{-n}\#\{\varepsilon:\ |\bareps|>u_\beta\}\le\Prob\big(|n\bareps|\ge nu_\beta\big)\le2\exp\Big(-\frac{n^2u_\beta^2}{2n}\Big)=2\exp\big(-\log(2/\beta)\big)=\beta. \]

The count is an integer not exceeding $\beta2^n$, and therefore not exceeding $k_\beta$, and based on \eqref{eq:q}, this means that $q_\beta\{|\bareps|\} \le u_\beta$.

\end{proof}

\subsection*{Proof of Theorem~\ref{thm:abr}}

\begin{proof}[]

Let $\eta=\alpha_0\delta$, $\gamma_1=\gamma_1(\eta)$ and $h=B(n,\eta/2)$, and let $J$ be the number of signs $+1$ in a uniform pattern, so that $J\sim{\rm Bin}(n,1/2)$ and $n\bareps=2J-n$. From $\Prob(J\ge\lceil n/2\rceil)\ge1/2\ge\eta/2$ it follows that $h\ge n/2$ and $\gamma_1 =(2h-n)/n\ge0$, and by the definition of $B(n,\eta/2)$, $\Prob(J\ge h+1)<\eta/2$. Since $|2J-n|>2h-n$ means that $J\ge h+1$ or $J\le n-h-1$,

\[ 2^{-n}\#\{\varepsilon:\ |\bareps|>\gamma_1\}=\Prob\big(|2J-n|>2h-n\big)=\Prob(J\ge h+1)+\Prob(J\le n-h-1), \]

and due to the symmetry of ${\rm Bin} (n,1/2)$ the two probabilities on the right are equal, so

\[ 2^{-n}\#\{\varepsilon:\ |\bareps|>\gamma_1\}=2\,\Prob(J\ge h+1)<\eta. \]

Hence $\#\{\varepsilon:\ |\bareps|>\gamma_1\}<\eta2^n$, and as the count is an integer, it does not exceed $k_\eta$.

\emph{Case $\gamma_1<1$.} Consider this case and note that $\eta+\alpha_0(1-\delta)=\alpha_0\le\alpha$. Then the conditions of Lemma~\ref{lem:split} are fulfilled for $\beta=\eta$, $\beta'=\alpha_0(1-\delta)$ and $\lambda=\gamma_1$, so any $w\in\mathcal W$ satisfies

\begin{equation}\label{eq:abrkey}
w\le q_{\alpha_0(1-\delta)}\{E(\varepsilon)\}+\gamma_1w.
\end{equation}

Therefore, $\mathcal W$ is bounded, and by definition $t_{\rm ABR}=q_ {\alpha_0(1-\delta)}\{E(\varepsilon)\}+\gamma_1f$, which is at least $q_ {\alpha_0 (1-\delta)}\{E(\varepsilon)\}$,

\[ c^*=\sup\mathcal W\le\frac{q_{\alpha_0(1-\delta)}\{E(\varepsilon)\}}{1-\gamma_1}\le\frac{t_{\rm ABR}}{1-\gamma_1}. \]

It is interesting to note that inequality \eqref{eq:abrkey} also holds for $w=c^*$, since $c^*$ is finite according to the last expression, and for a finite $c^*$, the proof of Lemma~\ref{lem:closed} shows that $\mathcal W=[0,c^*]$. Bounding the multiplier $c^*$ in the last term by $f+(c^*-f)_+$, we get

\[ c^*\le q_{\alpha_0(1-\delta)}\{E(\varepsilon)\}+\gamma_1f+\gamma_1(c^*-f)_+=t_{\rm ABR}+\gamma_1(c^*-f)_+. \]

\emph{Case $\gamma_1\ge1$.} Here we consider two subcases. If $c^*=+\infty$, then the right-hand side $t_{\rm ABR}+\gamma_1(c^*-f)_+$ is also equal to $+\infty$. If $c^*<\infty$, then $c^*\le\gamma_1c^*\le\gamma_1f+\gamma_1(c^*-f)_+\le t_{\rm ABR}+\gamma_1(c^*-f)_+$.

To summarize, in both cases $c^*\le t_{\rm ABR}+\gamma_1(c^*-f)_+$, and this value is $t_{\rm ABR}$ when $c^*\le f$.

\end{proof}

\section{Proofs for Section~\ref{sec:efficiency}}\label{app:eff}

The proofs of Theorems~\ref{thm:effK}, \ref{thm:effHD}(b), and~\ref{thm:effS} are based on a single comparison: the distribution law of $\sqrt n\,\phi(\bar y)$ and the sign-flip law of $\sqrt n\,E(\varepsilon)$ are compared with the same Gaussian distribution, and therefore we need Lemma~\ref{lem:gq} below, which converts these distances into relative errors of quantiles. Additionally, the proof of Theorem~\ref{thm:effHD}(a) does not use the Gaussian law and works by conditioning on the rows without their signs.

The upper $\gamma$-quantile of the random variable $X$ is given by 

\[ Q_\gamma(X)=\inf\{x\in\R:\ \Prob(X>x)\le\gamma\},\qquad \gamma\in(0,1), \]

so that $Q^*_\alpha=Q_\alpha(\phi(\bar y))$. The function $\gamma\mapsto Q_\gamma(X)$ is nonincreasing, and $\Prob(X>Q_\gamma(X))\le\gamma$ due to the continuity of the distribution function from the right. Let $\Phi$ denote the standard normal distribution function, $\Phi'$ its density, $z_\gamma=\Phi^{-1}(1-\gamma)$, $\|A\|_{\rm op}=\sup_{\|u\|_2\le1}\|Au\|_2$ the operator norm of a matrix, and $\lambda_{\min}(B)$ the smallest eigenvalue of a symmetric matrix $B$. For a function $v:\mathcal E\to\R$, we define

\[ \Prob_\varepsilon\big(v(\varepsilon)\le x\big)=\frac1N\,\#\{\varepsilon\in\mathcal E:\ v(\varepsilon)\le x\},\qquad x\in\R, \]

which for $\mathcal E=\{\pm1\}^n$ is the probability for a uniformly distributed pattern with fixed data, as in Section~\ref{sec:effK}.

\begin{lemma}\label{lem:kolm}
Let $X,X'$ be random variables with $\sup_x|\Prob(X\le x)-\Prob(X'\le x)|\le\Delta$. Then $Q_{\gamma+\Delta}(X')\le Q_\gamma(X)\le Q_{\gamma-\Delta}(X')$ for $\Delta<\gamma<1-\Delta$.
\end{lemma}

\begin{proof}

For this short lemma, we set $x_0 = Q_\gamma(X)$. Then $\Prob(X'>x_0) \le \Prob(X>x_0) + \Delta \le \gamma + \Delta$, which implies that $Q_{\gamma+\Delta}(X') \le x_0$. For $x<x_0$, we have $\Prob(X>x)>\gamma$, which implies $\Prob(X'>x)>\gamma-\Delta$ and $x<Q_{\gamma-\Delta}(X')$. And as $x$ approaches $x_0$, we obtain $x_0\le Q_{\gamma-\Delta}(X')$.

\end{proof}

\begin{lemma}\label{lem:kolmsign}
Let $v:\mathcal E\to\R$ and let $X'$ be a random variable with $\sup_x|\Prob_\varepsilon(v(\varepsilon)\le x)-\Prob(X'\le x)|\le\Delta$. Then for $\Delta+1/N<\gamma<1-\Delta$

\begin{equation}\label{eq:kolmq}
Q_{\gamma+\Delta}(X')\le q_\gamma\{v\}\le Q_{\gamma-1/N-\Delta}(X').
\end{equation}

\end{lemma}

\begin{proof}

If $X_v$ has the distribution of $v(\varepsilon)$ under $\Prob_\varepsilon$, then $\Prob(X_v>x)=N^{-1}\#\{\varepsilon:\ v(\varepsilon)>x\}$, therefore it follows from \eqref{eq:q} that $q_\gamma\{v\}=Q_{k_\gamma/N}(X_v)$, and since $\gamma-1/N<k_\gamma/N\le\gamma$, then $Q_\gamma(X_v)\le q_\gamma\{v\}\le Q_{\gamma-1/N}(X_v)$. It remains to apply Lemma~\ref{lem:kolm} to $X_v$ and $X'$ at levels $\gamma$ and $\gamma-1/N$.

\end{proof}

The next Lemma~\ref{lem:gq} follows from the Gaussian isoperimetric inequality \citep[Theorem 3.1]{Borell1975}, and we present a proof, because its constants $L$ and $m$ are included in all the following estimates.

\begin{lemma}\label{lem:gq}
Let $G\sim N(0,\Gamma)$ in $\R^K$ with $\Gamma$ positive semidefinite, let $\phi$ be a seminorm with $L=\sup_{\|u\|_2\le1}\phi(\Gamma^{1/2}u)>0$, and let $g_\gamma=Q_\gamma(\phi(G))$ and $m=g_{1/2}$. Then
\begin{enumerate}[label=(\roman*)]
\item the distribution function $F$ of $\phi(G)$ is continuous, $F(0)=0$, $F$ is strictly increasing on $[0,\infty)$, and $\Prob(\phi(G)>g_\gamma)=\gamma$ for $\gamma\in(0,1)$,
\item $|g_\gamma-g_{\gamma'}|\le L\,|z_\gamma-z_{\gamma'}|$ for $\gamma,\gamma'\in(0,1)$,
\item $0<m\le g_\gamma\le m+L\,z_\gamma$ for $\gamma\in(0,1/2]$.
\end{enumerate}
In particular, for $\alpha\in(0,1/2)$ and $\gamma,\gamma'\in[\alpha/2,1/2]$,

\begin{equation}\label{eq:gqlip}
|g_\gamma-g_{\gamma'}|\le\frac{L}{\Phi'(z_{\alpha/2})}\,|\gamma-\gamma'|.
\end{equation}

\end{lemma}

\begin{proof}

Throughout the proof, we write $Z\sim N(0,I_K)$ and $f(u)=\phi(\Gamma^{1/2}u)$, so that $\phi(G)\eqd f(Z)$, and by the definition of $L$

\begin{equation}\label{eq:gqlipf}
|f(u)-f(u')|\le f(u-u')\le L\,\|u-u'\|_2,\qquad u,u'\in\R^K.
\end{equation}

(i) Here, the zero set of the function $f$ is a linear subspace different from $\R^K$, since $L>0$, so $F(0)=0$. For $x>0$, the level set $\{f=x\}$ is the boundary of the convex set $A_x=\{f\le x\}$, which contains the ball $\{\|u\|_2\le x/L\}$, so $\{f=x\}$ has Lebesgue measure zero, and the function $F$ is continuous. According to Theorem 3.1 from the work of \citet{Borell1975}, for the standard Gaussian measure $\mathcal N_K$, its inner measure $\mathcal N_{K*}$, and the Euclidean unit ball $O$, the condition $\mathcal N_K(A)=\Phi(a)$ implies that $\mathcal N_{K*}(A+u\,O)\ge\Phi(a+u)$ for $u\ge0$. Further, according to \eqref{eq:gqlipf}, $A_x+u\,O\subseteq A_{x+Lu}$, therefore

\begin{equation}\label{eq:gqiso}
F(x+Lu)\ge\Phi\big(\Phi^{-1}(F(x))+u\big),\qquad x>0,\ u\ge0.
\end{equation}

Since $0<F(x)<1$ for $x>0$ (the set $\{f>x\}$ is open and nonempty), it follows from \eqref{eq:gqiso} that $F$ is strictly increasing on $[0,\infty)$, and $g_\gamma$ is the only value of $x$ for which $F(x)=1-\gamma$, which proves (i).

(ii) In this part, let $\gamma'<\gamma$, so that $g_{\gamma'}\ge g_\gamma>0$ and $F(g_\gamma)=1-\gamma$. Substituting $x=g_\gamma$ and $u=z_{\gamma'}-z_\gamma$ into \eqref{eq:gqiso}, we obtain $F(g_\gamma+Lu)\ge\Phi(z_{\gamma'})=1-\gamma'$, that is, $g_{\gamma'}\le g_\gamma+L(z_{\gamma'}-z_\gamma)$.

(iii) Here, $m>0$ follows from part (i), and from $\gamma\le1/2$ we get $g_\gamma\ge g_{1/2}=m$, while from (ii) with $\gamma'=1/2$ it follows that $g_\gamma\le m+Lz_\gamma$.

Finally, $\Phi'(z_\gamma)\ge\Phi'(z_{\alpha/2})$ for $\gamma\in[\alpha/2,1/2]$, so the mean value theorem gives $|z_\gamma-z_{\gamma'}|\le|\gamma-\gamma'|/\Phi'(z_{\alpha/2})$, and \eqref{eq:gqlip} follows from (ii).

\end{proof}

\subsection*{Proof of Theorem~\ref{thm:effK}}

\begin{proof}[]

Throughout the proof, $Z\sim N(0,I_K)$, $G=\Sigma^{1/2}Z$, $g_\gamma=Q_\gamma(\phi(G))$ and $L=\sup_{\|u\|_2\le1}\phi(\Sigma^{1/2}u)$, which is positive since $\Sigma^{1/2}$ is invertible and $\phi$ is not identically equal to zero. Given $L_0=L/\Phi'(z_{\alpha/2})$, Lemma~\ref{lem:gq} for $\Gamma=\Sigma$ yields

\begin{equation}\label{eq:Kgauss}
|g_\gamma-g_{\gamma'}|\le L_0|\gamma-\gamma'|\quad\text{for }\gamma,\gamma'\in[\alpha/2,1/2],\qquad g_\gamma\ge m>0\quad\text{for }\gamma\le1/2.
\end{equation}

Theorem~1.1 from \citet{Bentkus2005} deals with independent random vectors $X_1,\dots,X_n$ in $\R^d$ with zero mean, whose sum $S$ has an invertible covariance operator $C^2$, and a centered Gaussian vector $Z'$ with the same covariance. According to the theorem,

\[ \sup_{A}\big|\Prob(S\in A)-\Prob(Z'\in A)\big|\le c\,d^{1/4}\sum_{k=1}^n\E\|C^{-1}X_k\|_2^3, \]

where the supremum is taken over all convex sets, and $c$ is an absolute constant. For $X_t=y_t/\sqrt n$, for which $S=\sqrt n\,\bar y$ and $C=\Sigma^{1/2}$, and for the convex sets $\{v:\ \phi(v)\le x\}$, we have

\[ \rho_n:=\sup_x\big|\Prob(\sqrt n\,\phi(\bar y)\le x)-\Prob(\phi(G)\le x)\big|\le c\,K^{1/4}\,\frac{\E\|\Sigma^{-1/2}y_1\|_2^3}{\sqrt n}<\infty, \]

and then Lemma~\ref{lem:kolm} gives

\begin{equation}\label{eq:Koracle}
g_{\gamma+\rho_n}\le\sqrt n\,Q_\gamma(\phi(\bar y))\le g_{\gamma-\rho_n},\qquad\rho_n<\gamma<1-\rho_n.
\end{equation}

Under $\Prob_\varepsilon$, the vector $\sqrt n\,\bar e(\varepsilon)=\sum_t\varepsilon_te_t/\sqrt n$ is the sum of independent vectors with zero mean and covariance $\hat\Sigma=n^{-1}\sum_te_te_t^\top$. Let $Z$ be independent of the data, and $\Prob$ refer to $Z$ given fixed data. If $\hat\Sigma$ is invertible, then the same theorem and the inequality $\|\hat\Sigma^{-1/2}e_t\|_2\le\|\hat\Sigma^{-1/2}\|_{\rm op}\|e_t\|_2$ together give

\[ \Delta_n:=\sup_x\big|\Prob_\varepsilon(\sqrt n\,E(\varepsilon)\le x)-\Prob(\phi(\hat\Sigma^{1/2}Z)\le x)\big|\le\frac{c\,K^{1/4}\|\hat\Sigma^{-1/2}\|_{\rm op}^3}{\sqrt n}\cdot\frac1n\sum_{t=1}^n\|e_t\|_2^3, \]

otherwise we set $\Delta_n=1$. According to the strong law of large numbers (as well as the inequality $(a+b)^3\le4(a^3+b^3)$)

\[ \frac1n\sum_{t=1}^n\|e_t\|_2^3\le\frac4n\sum_{t=1}^n\|y_t\|_2^3+4\|\bar y\|_2^3\to4\,\E\|y_1\|_2^3,\qquad\hat\Sigma=\frac1n\sum_{t=1}^ny_ty_t^\top-\bar y\bar y^\top\to\Sigma, \]

almost surely, so that

\begin{equation}\label{eq:KDelta}
\Delta_n=O_P(n^{-1/2}).
\end{equation}

Now all that remains is to replace $\hat\Sigma$ with $\Sigma$ in the Gaussian distribution. The difference $X=\hat\Sigma^{1/2}-\Sigma^{1/2}$ solves the Sylvester equation $X\Sigma^{1/2}+\hat\Sigma^{1/2}X=\hat\Sigma-\Sigma$, and as $\Sigma^{1/2}\ge\lambda_{\min}(\Sigma)^{1/2}I_K$ and $\hat\Sigma^{1/2}\ge0$, Lemma~2.1 from \citet{Schmitt1992} gives

\[ \|X\|_{\rm op}\le\|\hat\Sigma-\Sigma\|_{\rm op}/\lambda_{\min}(\Sigma)^{1/2}. \]

As each entry of $n^{-1}\sum_t(y_ty_t^\top-\Sigma)$ has a mean of zero and a variance of no more than $\E\|y_1\|_2^4/n$ and $\E\|\bar y\|_2^2=\mathrm{tr}\,\Sigma/n$, we have

\[ \hat\kappa_n:=\big\|\Sigma^{-1/2}(\hat\Sigma^{1/2}-\Sigma^{1/2})\big\|_{\rm op}\le\frac{\|\hat\Sigma-\Sigma\|_{\rm op}}{\lambda_{\min}(\Sigma)}=O_P(n^{-1/2}), \]

and every $u\in\R^K$ satisfies $|\phi(\hat\Sigma^{1/2}u)-\phi(\Sigma^{1/2}u)|\le L\,\hat\kappa_n\|u\|_2$.

By Markov’s inequality, 

\[ \Prob(\|Z\|_2^2\ge2K)\le1/2, \]

so the median of $\|Z\|_2$ does not exceed $\sqrt{2K}$, and $1-\Phi(u)\le e^{-u^2/2}$ gives $z_{1/n}\le\sqrt{2\log n}$. Lemma~\ref{lem:gq}(iii), taking $\Gamma=I_K$ and $\phi=\|\cdot\|_2$ (i.e., $L=1$) then gives

\[ \Prob(\|Z\|_2>s_n)\le1/n \quad \text{for} \quad s_n=\sqrt{2K}+\sqrt{2\log n}. \]

Consequently, for $\hat g_\gamma=Q_\gamma(\phi(\hat\Sigma^{1/2}Z))$

\[ \Prob\big(\phi(\hat\Sigma^{1/2}Z)>g_{\gamma-1/n}+L\hat\kappa_ns_n\big)\le\Prob\big(\phi(G)>g_{\gamma-1/n}\big)+\Prob(\|Z\|_2>s_n)\le\gamma, \]

and for $x<g_{\gamma+1/n}-L\hat\kappa_ns_n$

\[ \Prob\big(\phi(\hat\Sigma^{1/2}Z)>x\big)\ge\Prob\big(\phi(G)>x+L\hat\kappa_ns_n\big)-\Prob(\|Z\|_2>s_n)>\gamma+\frac1n-\frac1n=\gamma, \]

and therefore

\begin{equation}\label{eq:Kcov}
g_{\gamma+1/n}-L\hat\kappa_ns_n\le\hat g_\gamma\le g_{\gamma-1/n}+L\hat\kappa_ns_n,\qquad\frac1n<\gamma<1-\frac1n.
\end{equation}

Lemma~\ref{lem:kolmsign} with $v=\sqrt n\,E$, $X'=\phi(\hat\Sigma^{1/2}Z)$ and $N=2^n$, combined with \eqref{eq:Kcov}, gives

\begin{equation}\label{eq:Kresamp}
g_{\gamma+\Delta_n+1/n}-L\hat\kappa_ns_n\le\sqrt n\,q_\gamma\{E(\varepsilon)\}\le g_{\gamma-2^{-n}-\Delta_n-1/n}+L\hat\kappa_ns_n
\end{equation}

for $\Delta_n+2^{-n}+1/n<\gamma<1-\Delta_n-1/n$.

Choosing $\beta_n=n^{-1/2}$, we get $k_{\beta_n}=\lfloor2^n/\sqrt n\rfloor\ge2=m_1$ for $n\ge3$, and $\beta_n<\alpha$ for large $n$, so using Lemma~\ref{lem:sandwich} with $\beta=\beta_n$, we have

\[ q_\alpha\{E(\varepsilon)\}\le c^*\le\frac{q_{\alpha-\beta_n}\{E(\varepsilon)\}}{1-\eta_n},\qquad\eta_n=\sqrt{\frac{2\log(2/\beta_n)}{n}}=\sqrt{\frac{\log(4n)}{n}}<1\quad(n\ge3). \]

Substituting \eqref{eq:Kresamp} at the levels $\alpha$ and $\alpha-\beta_n$, we get

\begin{equation}\label{eq:Kcstar}
g_{\alpha+\Delta_n+1/n}-L\hat\kappa_ns_n\le\sqrt n\,c^*\le\frac{g_{\alpha-\theta_n}+L\hat\kappa_ns_n}{1-\eta_n},\qquad\theta_n=\beta_n+2^{-n}+\Delta_n+\frac1n.
\end{equation}

All levels in \eqref{eq:Koracle} and \eqref{eq:Kcstar} on the event $\Omega_n=\{\theta_n+\rho_n\le\min(\alpha/2,1/2-\alpha)\}$ (the probability of which tends to one) lie in the range $[\alpha/2,1/2]$, and \eqref{eq:Kgauss} gives

\[ g_\alpha-L_0\theta_n-L\hat\kappa_ns_n\le\sqrt n\,c^*\le\frac{g_\alpha+L_0\theta_n+L\hat\kappa_ns_n}{1-\eta_n},\qquad g_\alpha-L_0\rho_n\le\sqrt n\,Q^*_\alpha\le g_\alpha+L_0\rho_n. \]

By dividing and using the inequalities $(g_\alpha-x)/(g_\alpha+y)\ge1-(x+y)/m$ and $(g_\alpha+x)/(g_\alpha-y)\le(1+x/m)/(1-y/m)$ for $x,y\ge0$ and $y<m\le g_\alpha$, on the event $\Omega_n\cap\{\delta_n<1\}$ we have

\[ 1-\delta_n\le\frac{c^*}{Q^*_\alpha}\le\frac{1+\delta_n}{(1-\eta_n)(1-\delta_n)},\qquad\delta_n=\frac{L_0(\theta_n+\rho_n)+L\hat\kappa_ns_n}{m}. \]

Combining \eqref{eq:KDelta} with the bounds for $\rho_n$ and $\hat\kappa_n$, we get $\theta_n+\rho_n=O_P(n^{-1/2})$ and $\hat\kappa_ns_n=O_P(\sqrt{\log n/n})$, so that $\delta_n=O_P(\sqrt{\log n/n})$. As $\eta_n=O(\sqrt{\log n/n})$ and

\[ \frac{1+\delta_n}{(1-\eta_n)(1-\delta_n)}-1=\frac{2\delta_n+\eta_n-\eta_n\delta_n}{(1-\eta_n)(1-\delta_n)}=O_P(\delta_n+\eta_n), \]

we finally get the ratio to the ideal quantile

\[ \frac{c^*}{Q^*_\alpha}=1+O_P(\delta_n+\eta_n)=1+O_P\Big(\sqrt{\frac{\log n}{n}}\Big). \]

Here, we prove that $\Prob(\phi(\bar y)>(1+\kappa_n)c^*)\to\alpha$ for every deterministic sequence $\kappa_n\to0$. Let $a_n=n^{-1/4}$ and $\Omega'_n=\Omega_n\cap\{\Delta_n\le a_n,\ L\hat\kappa_ns_n\le a_n\}$, then $\Prob(\Omega'_n)\to1$. As the quantile is nonincreasing in the level, \eqref{eq:Kcstar} provides deterministic bounds on $\Omega'_n$

\[ \ell_n:=g_{\alpha+a_n+1/n}-a_n\le\sqrt n\,c^*\le\frac{g_{\alpha-\beta_n-2^{-n}-a_n-1/n}+a_n}{1-\eta_n}=:u_n, \]

and $\ell_n,u_n\to g_\alpha$ by \eqref{eq:Kgauss}. Using the definition of $\rho_n$ for a deterministic sequence $\kappa_n\to0$, we have

\[ \Prob\big(\phi(G)>(1+\kappa_n)u_n\big)-\rho_n-\Prob(\Omega_n'^c)\le\Prob\big(\phi(\bar y)>(1+\kappa_n)c^*\big)\le\Prob\big(\phi(G)>(1+\kappa_n)\ell_n\big)+\rho_n+\Prob(\Omega_n'^c). \]

Both $(1+\kappa_n)\ell_n$ and $(1+\kappa_n)u_n$ tend to $g_\alpha$, and it follows from Lemma~\ref{lem:gq}(i) that the distribution function of $\phi(G)$ is continuous and $\Prob(\phi(G)>g_\alpha)=\alpha$, so both bounds tend to $\alpha$, and

\begin{equation}\label{eq:Klevel}
\Prob\big(\phi(\bar y)>(1+\kappa_n)c^*\big)\to\alpha\qquad\text{for every deterministic }\kappa_n\to0,
\end{equation}

which, for $\kappa_n=0$, is the second statement of the theorem.

\end{proof}

\subsection*{Proof of Theorem~\ref{thm:sharp}}

\begin{proof}[]

In this proof, we use the notation of the proof of Theorem~\ref{thm:effK} and, as usual for almost sure (a.s.) limits, regard $y_1,\dots,y_n$ as the first $n$ terms of one i.i.d.\ sequence $y_1,y_2,\dots$. We also write $\Prob$ for the probability over $Z$ given fixed data and use

\[ \hat G=\hat\Sigma^{1/2}Z,\qquad\bar F_n(s)=\Prob(\phi(\hat G)\ge s),\qquad\bar F(s)=\Prob(\phi(G)\ge s), \]

\[ x=\sqrt n\,c^*,\qquad x_0=\sqrt n\,t_{\rm raw},\qquad I=[g_\alpha/2,2g_\alpha],\qquad I'=[2g_\alpha/3,3g_\alpha/2]. \]

By \eqref{eq:Kresamp} at the level $\alpha$, \eqref{eq:Kcstar} and \eqref{eq:Kgauss}, both $x$ and $x_0$ tend to $g_\alpha$ in probability.

First, we note that the vectors $v$ for which $\phi(v)=0$ form a linear subspace different from $\R^K$. In its orthogonal complement $V$ of dimension $d\ge1$, $\phi$ is a norm with a bounded unit ball $B_V$, and $\phi(v)=\phi(P_Vv)$ for the orthogonal projection $P_V$ onto $V$, due to subadditivity. Second, if $h_\Gamma$ is the density of the distribution $P_VN(0,\Gamma)$ on $V$ for a positive definite matrix $\Gamma$, then

\[ \Prob\big(\phi(N(0,\Gamma))\le s\big)=\int_{sB_V}h_\Gamma=s^d\int_{B_V}h_\Gamma(su)\,du,\qquad s>0, \]

which is continuously differentiable for $s>0$, and by the dominated convergence theorem, its derivative $f_\Gamma(s)$ is jointly continuous in $(s,\Gamma)$. This derivative is positive because $h_\Gamma$ is bounded below on $(s+a)B_V$ and ${\rm vol}((s+a)B_V\setminus sB_V)\ge d\,s^{d-1}a\,{\rm vol}(B_V)$. With $f = f_\Sigma$ and $f_n = f_{\hat\Sigma}$, we have $\bar F_n' = -f_n$, and $\hat\Sigma \to \Sigma$ implies $f_n \to f$ uniformly on $[g_\alpha/3, 3g_\alpha]$ a.s. Consequently, there exist $C$ and, for almost every data sequence, $n_1$ such that

\begin{equation}\label{eq:Sdens}
\sup_{s\in I}\Prob\big(|\phi(\hat G)-s|\le a\big)\le C\,a,\qquad0<a\le g_\alpha/6,\ n\ge n_1.
\end{equation}

Let $J$ be the number of $+1$ signs in a uniform pattern $\varepsilon$, $\bareps_j=(2j-n)/n$, and $\mathcal J_n=\{j:\ |2j-n|\le\sqrt n\log n\}$. Here, we prove that

\begin{equation}\label{eq:Scond}
\sup_{j\in\mathcal J_n}\ \sup_{s\in I}\Big|\Prob_\varepsilon\big(\sqrt n\,E(\varepsilon)\ge s\sqrt{1-\bareps_j^2}\ \big|\ J=j\big)-\bar F_n(s)\Big|=o(n^{-1/2}),
\end{equation}

a.s., and the same applies with $>$ instead of $\ge$. When $J=j$, the set $S=\{t:\ \varepsilon_t=1\}$ is uniformly distributed among the subsets of size $j$, and $\sum_te_t=0$ and $4j(n-j)/n^2=1-\bareps_j^2$ give

\[ \sqrt n\,\bar e(\varepsilon)=\frac2{\sqrt n}\sum_{t\in S}e_t=\sqrt{1-\bareps_j^2}\;Z_j,\qquad Z_j=(jq_j)^{-1/2}\sum_{t\in S}e_t,\quad q_j=1-\frac jn. \]

The important point to note is that in the notation of \citet{BabuSingh1985}, where the size of the population $N$ (not the number of patterns) is equal to $n$, $Z_j$ is the normalized sum of a sample of $j$ elements taken without replacement from the population $\{e_1,\dots,e_n\}$ with a mean of zero and a covariance of $\hat\Sigma$ (their $\Sigma_N$). Their Theorem 1 provides an expansion of the distribution law of $Z_j$ with a remainder term of $o(j^{-1/2})$ provided that three conditions are met for some $\delta'>0$: (1) the moments of order $3+\delta'$ of the population are bounded, (2) the sampling fraction is no more than $1-\delta'$, and (3) the empirical distribution law of the population weakly converges to a strongly non-lattice law, i.e., to a law whose characteristic function has an absolute value less than one for any $t\ne0$. It remains to verify these conditions with $\delta'=1/4$ for almost every data sequence, so that the expansion applies to $Z_j$:

(1) holds because $n^{-1}\sum_t\|e_t\|_2^4\le8n^{-1}\sum_t\|y_t\|_2^4+8\|\bar y\|_2^4$ is bounded by the strong law of large numbers,

(2) holds because $j/n \leq 3/4$ for $j \in \mathcal J_n$ and large $n$, and only such $j$ matter, since patterns with $J \notin \mathcal J_n$ have a total probability of $o(n^{-1/2})$ (shown below),

(3) holds because the empirical distribution law of $e_1,\dots,e_n$ converges weakly to the distribution law $P$ of $y_1$ as $\bar y\to0$, and $P$ is strongly non-lattice, as we will prove below.

Condition (3) is the reason why Theorem~\ref{thm:sharp} assumes that $P$ has an absolutely continuous component: for a lattice distribution such as the one shown in Figure~\ref{fig:level}(a), the expansion with the remainder $o(j^{-1/2})$ may fail. To verify condition (3), we write $P=p_0P_{\rm ac}+(1-p_0)P_{\rm s}$ with $p_0>0$ and $P_{\rm ac}$ absolutely continuous, and let $\hat P$, $\hat P_{\rm ac}$ be the characteristic functions. If $|\hat P_{\rm ac}(t)|=1$ for some $t\ne0$, then $P_{\rm ac}$ is concentrated on the Lebesgue-null set $\{v:\ t^\top v\in\theta+2\pi\mathbb Z\}$, which is impossible, so $|\hat P(t)|\le p_0|\hat P_{\rm ac}(t)|+1-p_0<1$ for $t\ne0$.

For a bounded measurable function $g$, the theorem compares $\E g(Z_j)$ with $\int g\,dP_j$, where $P_j$ is a signed measure with density

\[ \varphi_{\hat\Sigma}-\frac{q_j-p_j}{\sqrt{jq_j}}\,\varpi,\qquad p_j=\frac jn, \]

where $\varphi_{\hat\Sigma}$ is the density of the distribution $N(0,\hat\Sigma)$, and $\varpi$ (their $p_N$) is a linear combination of the third-order partial derivatives of $\varphi_{\hat\Sigma}$, i.e., an odd function. Their estimate states that

\[ \Big|\E g(Z_j)-\int g\,dP_j\Big|\le c\,\bar\omega(g,d_j,\Phi_{\hat\Sigma})+o(j^{-1/2}),\qquad d_j=o(j^{-1/2}), \]

where $d_j$ (their $\delta_n$), $c$, and the term $o(j^{-1/2})$ depend on $g$ only through its bound. For $g=\mathbf 1\{\phi\ge s\}$ and $g=\mathbf 1\{\phi>s\}$, the sets are symmetric, and $\varpi$ is an odd function, so $\int g\,dP_j=\bar F_n(s)$. The oscillations of the function $g$ in a ball of radius $d$ centered at point $v$ are zero unless $|\phi(v)-s|\le L_\phi d$, where $L_\phi=\sup_{\|u\|_2\le1}\phi(u)$, and from \eqref{eq:Sdens} it follows that

\[ \bar\omega(g,d_j,\Phi_{\hat\Sigma})\le CL_\phi d_j=o(n^{-1/2})\qquad\text{uniformly in }s\in I. \]

All these indicators are bounded by 1, so the theorem gives \eqref{eq:Scond} for any pair of sequences $j_n\in\mathcal J_n$ and $s_n\in I$. Uniformity in $j$ and $s$ follows from a proof by contradiction: if for some sequence of data from the set on which the probability is equal to one, the supremum in \eqref{eq:Scond} exceeded $\epsilon n^{-1/2}$ for infinitely many $n$, we could choose $j_n$ and $s_n$ for which this value would be achieved, and then the theorem would not hold for these sequences.

For $w\ge0$ let

\[ \mathcal G_n(w)=\Prob_\varepsilon\big(\sqrt n\,E(\varepsilon)\ge(1-|\bareps|)w\big),\qquad\mathcal G^0_n(w)=\Prob_\varepsilon\big(\sqrt n\,E(\varepsilon)\ge w\big), \]

and denote $\mathcal G_n(w+)$ and $\mathcal G^0_n(w+)$ as their right limits. For large $n$, the value $c^*$ is finite, and the set $\mathcal W$ from the proof of Lemma~\ref{lem:closed} is equal to $[0,c^*]$, so \eqref{eq:qchar} and the definition of $t_{\rm raw}$ give

\[ \mathcal G_n(x)\ge\frac{k_\alpha+1}{2^n}>\alpha\ge\frac{k_\alpha}{2^n}\ge\mathcal G_n(x+),\qquad\mathcal G^0_n(x_0)>\alpha\ge\mathcal G^0_n(x_0+). \]

By \eqref{eq:hoeff}, the probability of patterns with $J\notin\mathcal J_n$ does not exceed $2e^{-(\log n)^2/2}=o(n^{-1/2})$. If $|\bareps|<1$, then

\[ 1-|\bareps|=\sqrt{1-\bareps^2}\,\varrho(|\bareps|) \quad \text{with} \quad \varrho(a)=\sqrt{(1-a)/(1+a)}, \]

and for $j\in\mathcal J_n$ and $w\in I'$, the numbers $w\varrho(|\bareps_j|)$ and $w/\sqrt{1-\bareps_j^2}$ are equal to $w(1+O(\log n/\sqrt n))$ and lie in $I$ for large $n$. Then, from \eqref{eq:Scond}, it follows that, uniformly in $w\in I'$,

\begin{align}
\mathcal G_n(w)&=\E\Big[\bar F_n\big(w\varrho(|\bareps_J|)\big)\mathbf 1\{J\in\mathcal J_n\}\Big]+o(n^{-1/2}),\label{eq:Sav}\\
\mathcal G^0_n(w)&=\E\Big[\bar F_n\Big(\frac w{\sqrt{1-\bareps_J^2}}\Big)\mathbf 1\{J\in\mathcal J_n\}\Big]+o(n^{-1/2}),\nonumber
\end{align}

a.s., and the version of \eqref{eq:Scond} with $>$ and the continuity of $\bar F_n$ give the same expansions for $\mathcal G_n(w+)=\Prob_\varepsilon(\sqrt n\,E(\varepsilon)>(1-|\bareps|)w)$ and $\mathcal G^0_n(w+)=\Prob_\varepsilon(\sqrt n\,E(\varepsilon)>w)$.

On $\mathcal J_n$, the number $a=|\bareps_J|$ does not exceed $\log n/\sqrt n$, so 

\[ \varrho(a)=1-a+O(a^2),\quad (1-a^2)^{-1/2}=1+ O(a^2), \]

and the mean value theorem with the uniform convergence $f_n\to f$ and the uniform continuity of $f$ on $I$ gives

\[ \bar F_n\big(w\varrho(a)\big)=\bar F_n(w)+f(w)\,w\,a+o(a)+O(a^2),\qquad\bar F_n\Big(\frac w{\sqrt{1-a^2}}\Big)=\bar F_n(w)+O(a^2), \]

where $o(a)$ and $O(a^2)$ are uniform over $w\in I'$. Here $\E\bareps_J^2=1/n$, and

\[ \sqrt n\,\E|\bareps_J|=\E|2J-n|/\sqrt n\to\sqrt{2/\pi} \]

according to the central limit theorem, since $|2J-n|/\sqrt n$ has a second moment equal to one and is therefore uniformly integrable. Substituting into \eqref{eq:Sav}, we obtain that, uniformly with respect to $w\in I'$,

\[ \mathcal G_n(w)=\bar F_n(w)+f(w)\,w\,\sqrt{\frac2{\pi n}}+o(n^{-1/2}),\qquad\mathcal G^0_n(w)=\bar F_n(w)+o(n^{-1/2}), \]

a.s., and the same holds for the right limits.

\emph{Completion of the proof.} On the event $\{x, x_0 \in I'\}$, the probability of which tends to one,

\[ \bar F_n(x)+f(x)\,x\,\sqrt{\frac2{\pi n}}=\alpha+o(n^{-1/2}),\qquad\bar F_n(x_0)=\alpha+o(n^{-1/2}). \]

By the mean value theorem

\[ \bar F_n(x_0)-\bar F_n(x)=f_n(\xi)(x-x_0) \]

for some $\xi$ between $x$ and $x_0$, and $f_n(\xi)\to f(g_\alpha)>0$, $f(x)\,x\to f(g_\alpha)\,g_\alpha$ in probability. Hence

\[ \sqrt n\,(x-x_0)\to g_\alpha\sqrt{2/\pi} \quad \text{in probability}, \]

and dividing by $x_0\to g_\alpha$ gives $\sqrt n\,(c^*/t_{\rm raw}-1)\to\sqrt{2/\pi}$.

For the second statement, let $H=\Sigma^{-1/2}\hat\Sigma\Sigma^{-1/2}-I_K$, so that

\[ \|H\|_{\rm op}\le\|\hat\Sigma-\Sigma\|_{\rm op}/\lambda_{\min}(\Sigma)=O_P(n^{-1/2}) \]

by the proof of Theorem~\ref{thm:effK}. If $\|H\|_{\rm op} \le 1/2$, then the eigenvalues $\lambda_i$ of $H$ satisfy $\lambda_i - \log(1 + \lambda_i) \le \lambda_i^2$, so the Kullback--Leibler divergence of $N(0,\hat\Sigma)$ from $N(0,\Sigma)$ is

\[ \frac12\sum_i(\lambda_i-\log(1+\lambda_i))\le\frac12\|H\|_F^2, \]

and Pinsker's inequality gives, on $\{\|H\|_{\rm op}\le1/2\}$,

\[ \sup_s|\bar F_n(s)-\bar F(s)|\le d_{\rm TV}\big(N(0,\hat\Sigma),N(0,\Sigma)\big)\le\tfrac12\|H\|_F\le\tfrac12\sqrt K\,\|H\|_{\rm op}=O_P(n^{-1/2}). \]

Thus 

\[ \bar F(x_0)=\alpha+O_P(n^{-1/2}). \]

Since $\bar F(g_\alpha)=\alpha$ by Lemma~\ref{lem:gq}(i), it follows from the mean value theorem that $\bar F(g_\alpha)-\bar F(x_0)=f(\xi')(x_0-g_\alpha)$, where $\xi'$ lies between $x_0$ and $g_\alpha$, and $f(\xi') \to f(g_\alpha)>0$, so

\[ x_0=g_\alpha+O_P(n^{-1/2}). \]

Taking into account \eqref{eq:Koracle}, \eqref{eq:Kgauss} and using $\rho_n=O(n^{-1/2})$, $\sqrt n\,Q^*_\alpha=g_\alpha+O(n^{-1/2})$, so

\[ t_{\rm raw}/Q^*_\alpha=1+O_P(n^{-1/2}) \implies \frac{c^*}{Q^*_\alpha}=\frac{c^*}{t_{\rm raw}}\cdot\frac{t_{\rm raw}}{Q^*_\alpha}=1+O_P(n^{-1/2}). \qedhere \]

\end{proof}

\subsection*{Proof of Theorem~\ref{thm:effHD}}

\begin{proof}[]

We first prove part (b), that is, for $\log K=o(n^{1/5})$,

\[ \frac{c^*}{Q^*_\alpha}\to1\ \text{in probability},\qquad\Prob\big(\|\bar y\|_\infty>(1+\kappa_n)c^*\big)\to\alpha, \]

and then we prove part (a) for $\log K=o(n^{1/3})$,

\[ \Prob\big(\|\bar y\|_\infty>c^*\big)\to\alpha. \]

\emph{Proof of (b).} Write

\[ T=\sqrt n\,\|\bar y\|_\infty,\qquad\Sigma_n=n^{-1}\sum_t\E y_ty_t^\top,\qquad G\sim N(0,\Sigma_n),\qquad T_G=\|G\|_\infty,\qquad g_\gamma=Q_\gamma(T_G). \]

For $\phi=\|\cdot\|_\infty$ and $\Gamma=\Sigma_n$, the constant from Lemma~\ref{lem:gq} is equal to $L=\max_j(\Sigma_n)_{jj}^{1/2}\in[\underline\sigma,\bar\sigma]$. The median $m$ of $T_G$ is no less than the median $(\Sigma_n)_{11}^{1/2}z_{1/4}$ of $|G_1|$, since $T_G\ge|G_1|$. For $\lambda>0$, Jensen’s inequality and $\E e^{\pm\lambda G_j}\le e^{\lambda^2\bar\sigma^2/2}$ yield

\[ e^{\lambda\E T_G}\le\E e^{\lambda T_G}\le\sum_{j=1}^K\big(\E e^{\lambda G_j}+\E e^{-\lambda G_j}\big)\le2K\,e^{\lambda^2\bar\sigma^2/2}, \]

and $\lambda=\sqrt{2\log(2K)}/\bar\sigma$ gives the first bound in

\[ \E T_G\le\bar\sigma\sqrt{2\log(2K)},\qquad m\le2\E T_G, \]

and the second follows from Markov's inequality. Combining Lemma~\ref{lem:gq}(iii) and \eqref{eq:gqlip} with $L_0=\bar\sigma/\Phi'(z_{\alpha/2})$, we get

\begin{equation}\label{eq:HDgauss}
\underline\sigma z_{1/4}\le g_\gamma\le\bar x_K:=\bar\sigma\big(2\sqrt{2\log(2K)}+z_{\alpha/2}\big),\qquad|g_\gamma-g_{\gamma'}|\le L_0|\gamma-\gamma'|,\qquad\gamma,\gamma'\in[\alpha/2,1/2].
\end{equation}

Lemma~A.1 of \citet{CCK2017} (which the authors attribute to \citet{Nazarov2003}) states that a centered Gaussian vector $Y\in\R^p$ with $p\ge3$ \citep[\S1.1]{CCK2017} and $\E Y_j^2\ge b$ for all $j$ satisfies

\[ \Prob(Y\le y+a)-\Prob(Y\le y)\le C\,a\sqrt{\log p},\qquad y\in\R^p,\ a>0, \]

where the inequalities hold for each coordinate, and $C$ depends only on $b$. As in the proof of their Corollary 5.1, we apply it to $Y=(G_1,\dots,G_K,-G_1,\dots,-G_K,G_1)\in\R^{2K+1}$ with $b=\underline\sigma^2$, for which $\{Y\le x\mathbf 1\}=\{T_G\le x\}$, and, with $C_0$ depending only on $\underline\sigma$, we have

\begin{equation}\label{eq:HDnaz}
\Prob(x<T_G\le x+a)\le C_0\,a\sqrt{\log(2K+1)}\le C_0\,a\sqrt{\log(3K)},\qquad x\in\R,\ a>0.
\end{equation}

We compare the distribution laws of $T$ and $\sqrt n\,E(\varepsilon)$ with the distribution law of $T_G$ through

\[ \rho_n=\sup_x\big|\Prob(T\le x)-\Prob(T_G\le x)\big|,\qquad D_n=\sup_x\big|\Prob(T\le x)-\Prob_\varepsilon(\sqrt n\,E(\varepsilon)\le x)\big|. \]

Lemmas~4.3 and~4.6 of \citet{CCKK2022} concern independent centered vectors $X_t\in\R^p$, $p\ge2$, with their (1) Condition~E, where $\E\exp(|X_{tj}|/B_n)\le2$ for a constant $B_n\ge1$ and the (2) Condition~M, where $b_1^2\le n^{-1}\sum_t\E X_{tj}^2$ and $n^{-1}\sum_t\E X_{tj}^4\le B_n^2b_2^2$ for constants $0<b_1\le b_2$. For

\[ T_n=\max_jn^{-1/2}\sum_tX_{tj},\qquad T^G_n=\max_jG'_j,\qquad G'\sim N\Big(0,n^{-1}\sum_t\E X_tX_t^\top\Big), \]

and the multiplier bootstrap statistic $T^*_n=\max_jn^{-1/2}\sum_tW_t(X_{tj}-\bar X_j)$ with weights satisfying their condition (12), which is satisfied by Rademacher weights, they show that

\begin{align*}
\sup_x\big|\Prob(T_n\le x)-\Prob(T^G_n\le x)\big|&\le C\Big(\frac{B_n^2\log^5(pn)}{n}\Big)^{1/4},\\
\sup_x\big|\Prob(T_n\le x)-\Prob(T^*_n\le x\mid X_1,\dots,X_n)\big|&\le C\Big(\frac{B_n^2\log^5(pn)}{n}\Big)^{1/4},
\end{align*}

and the second with probability at least $1-2/n-3\upsilon_n$ ($\upsilon_n=(B_n^2\log^3(pn)/n)^{1/2}$, and $C$ depends only on $b_1$ and $b_2$). We take $X_t=(y_t,-y_t)\in\R^{2K}$, so that $p=2K$, $T_n=T$, $T^G_n\eqd T_G$, and $T^*_n=\sqrt n\,E(\varepsilon)$ under $\Prob_\varepsilon$. From

\[ x^k\le k!\,b^ke^{x/b}\implies\E|y_{tj}|^k\le2\,k!\,b^k,\quad\text{in particular}\quad\E y_{tj}^4\le48b^4, \]

so Conditions~E and~M are satisfied with $B_n=\max\{1,b,7b^2\}$, $b_1=\underline\sigma$ and $b_2=\max\{\underline\sigma,1\}$. With $\zeta_n=C(B_n^2\log^5(2Kn)/n)^{1/4}$ and $\Omega_n=\{D_n\le\zeta_n\}$, we obtain

\[ \rho_n\le\zeta_n\to0,\qquad\pi_n:=\Prob(\Omega_n^c)\le2/n+3\upsilon_n\to0, \]

because $\log K=o(n^{1/5})$, where $C$ is the larger of the two constants.

On the event $\Omega_n$, the distance between the distribution of $\sqrt n\,E(\varepsilon)$ under $\Prob_\varepsilon$ and the distribution of $T_G$ does not exceed $\tau_n=\zeta_n+\rho_n$, and Lemma~\ref{lem:kolmsign} gives

\[ g_{\gamma+\tau_n}\le\sqrt n\,q_\gamma\{E(\varepsilon)\}\le g_{\gamma-2^{-n}-\tau_n}. \]

By the same argument as that in the proof of Theorem~\ref{thm:effK}, with Lemma~\ref{lem:sandwich} for $\beta_n=n^{-1/2}$, $\eta_n=\sqrt{\log(4n)/n}$, and Lemma~\ref{lem:kolm} for $Q^*_\gamma$, we conclude that

\begin{equation}\label{eq:HDcstar}
g_{\alpha+\tau_n}\le\sqrt n\,c^*\le\frac{g_{\alpha-\theta_n}}{1-\eta_n}\quad\text{on }\Omega_n,\qquad g_{\gamma+\rho_n}\le\sqrt n\,Q^*_\gamma\le g_{\gamma-\rho_n},
\end{equation}

where $\theta_n=\beta_n+2^{-n}+\tau_n$, and the sequences $\tau_n$, $\theta_n$, and $\eta_n$ are deterministic and tend to zero. For large $n$, all levels in \eqref{eq:HDcstar} lie in the range $[\alpha/2,1/2]$, and \eqref{eq:HDgauss}, together with the two elementary inequalities from the proof of Theorem~\ref{thm:effK} with $m$ replaced by $\underline\sigma z_{1/4}\le g_\alpha$ and without the term $L\hat\kappa_ns_n$, gives on $\Omega_n$

\[ 1-\delta_n\le\frac{c^*}{Q^*_\alpha}\le\frac{1+\delta_n}{(1-\eta_n)(1-\delta_n)},\qquad\delta_n=\frac{L_0(\theta_n+\rho_n)}{\underline\sigma z_{1/4}}\to0. \]

Hence, for every $u>0$ and large $n$, 

\[ \{|c^*/Q^*_\alpha-1|>u\}\subseteq\Omega_n^c, \]

which has probability $\pi_n\to0$, so $c^*/Q^*_\alpha\to1$ in probability.

It remains to show that, for every deterministic sequence $\kappa_n$ with $|\kappa_n|\le1/2$ and $\kappa_n\log(3K)\to0$,

\begin{equation}\label{eq:HDlevel}
\Prob\big(\|\bar y\|_\infty>(1+\kappa_n)c^*\big)\to\alpha,
\end{equation}

which, when $\kappa_n=0$, gives $\Prob(\|\bar y\|_\infty>c^*)\to\alpha$. By \eqref{eq:HDcstar}, $\{\|\bar y\|_\infty>(1+\kappa_n)c^*\}\cap\Omega_n\subseteq\{T>(1+\kappa_n)x_n\}$ with $x_n=g_{\alpha+\tau_n}\le\bar x_K$, and Lemma~\ref{lem:gq}(i) and \eqref{eq:HDnaz} give

\[ \Prob\big(\|\bar y\|_\infty>(1+\kappa_n)c^*\big)\le\Prob\big(T_G>(1+\kappa_n)x_n\big)+\rho_n+\pi_n\le\alpha+\tau_n+C_0|\kappa_n|\,\bar x_K\sqrt{\log(3K)}+\rho_n+\pi_n. \]

From below, let $x'_n=g_{\alpha-\theta_n}$ and $1+\lambda_n=(1+|\kappa_n|)/(1-\eta_n)$, so that $0\le\lambda_n\le2(|\kappa_n|+\eta_n)$ when $\eta_n\le1/2$. Then $\{T>(1+\lambda_n)x'_n\}\cap\Omega_n\subseteq\{\|\bar y\|_\infty>(1+\kappa_n)c^*\}$ by \eqref{eq:HDcstar}, and

\[ \Prob\big(\|\bar y\|_\infty>(1+\kappa_n)c^*\big)\ge\Prob\big(T_G>(1+\lambda_n)x'_n\big)-\rho_n-\pi_n\ge\alpha-\theta_n-C_0\lambda_n\,\bar x_K\sqrt{\log(3K)}-\rho_n-\pi_n. \]

Finally, $\bar x_K\sqrt{\log(3K)}\le C\log(3K)$ and

\[ \eta_n\log(3K)\le\sqrt{\frac{\log(4n)}{n}}\,\big(\log3+o(n^{1/5})\big)\to0,\qquad\lambda_n\log(3K)\le2(|\kappa_n|+\eta_n)\log(3K)\to0, \]

so the terms with $\kappa_n$ and $\lambda_n$ in both bounds tend to zero, both bounds tend to $\alpha$, and \eqref{eq:HDlevel} follows.

\emph{Proof of (a).} In this part, we keep $T:=\sqrt n\,\|\bar y\|_\infty$, write $T(\varepsilon)=\sqrt n\,\|\bar y(\varepsilon)\|_\infty$, that is, $T=T(\mathbf 1)$, and take $\beta_n=n^{-1/2}$, $\gamma_n=\alpha-2\beta_n$, and $\eta_n=\sqrt{2\log(2/\beta_n)/n}$. From Theorem~\ref{thm:cstar}, it follows that $\Prob(\|\bar y\|_\infty>c^*)\le\alpha$, so it is sufficient to show that

\[ \liminf_n\Prob(\|\bar y\|_\infty>c^*)\ge\alpha. \]

Write $y_t=\sigma_tw_t$, where $w_t=y_t$ if the first nonzero coordinate of $y_t$ is positive, $w_t=-y_t$ if it is negative, and $w_t=0$ if $y_t=0$, and $\sigma_t\in\{\pm1\}$ is chosen uniformly and independently of everything else when $y_t=0$. The independence of the rows together with $y_t\eqd-y_t$ implies that $\sigma=(\sigma_1,\dots,\sigma_n)$ is uniformly distributed on $\{\pm1\}^n$ and independent of $w=(w_1,\dots,w_n)$. Since $\bar y(\varepsilon)=n^{-1}\sum_t\varepsilon_t\sigma_tw_t$ and $\varepsilon\mapsto(\varepsilon_t\sigma_t)_t$ is a bijection of the set $\{\pm1\}^n$, the distribution law of $T(\varepsilon)$ under $\Prob_\varepsilon$, and along with it $q^0_\gamma=q_\gamma\{T(\varepsilon)\}$, depends only on $w$, and for a given $w$, the statistic $T$ has this distribution law. Thus, for every measurable function $x(w)$

\begin{equation}\label{eq:Ccond}
\Prob\big(T>x(w)\mid w\big)=\Prob_\varepsilon\big(T(\varepsilon)>x(w)\big).
\end{equation}

Put $c_n=5b\log(Kn)$, $s_j^2=n^{-1}\sum_ty_{tj}^2=n^{-1}\sum_tw_{tj}^2$ and

\[ W_n=\Big\{\max_{t,j}|y_{tj}|\le c_n\Big\}\cap\Big\{\frac{\underline\sigma^2}2\le s_j^2\le2\bar\sigma^2\ \text{for all }j\Big\}, \]

an event that depends on $w$ only. Markov’s inequality gives $\Prob(|y_{tj}|>x)\le2e^{-x/b}$, so the first event fails with probability no more than $2nK(Kn)^{-5}$. As for the second, let $\tilde y_{tj}=y_{tj}\mathbf 1\{|y_{tj}|\le c_n\}$. By $\E y_{tj}^4\le48b^4$ and the Cauchy--Bunyakovsky--Schwarz inequality, $0\le\E y_{tj}^2-\E\tilde y_{tj}^2\le(48b^4)^{1/2}(2(Kn)^{-5})^{1/2}$, and also the variables $\tilde y_{tj}^2-\E\tilde y_{tj}^2$ are bounded by $c_n^2$ with variance no greater than $48b^4$. From Bernstein's inequality with $u=\underline\sigma^2/4$ and the union bound over $j$, we observe

\[ \Prob\Big(\max_j\Big|\frac1n\sum_t\big(\tilde y_{tj}^2-\E\tilde y_{tj}^2\big)\Big|\ge u\Big)\le2K\exp\Big(-\frac{nu^2}{96b^4+2c_n^2u/3}\Big)\to0, \]

because $\log K\cdot\log^2(Kn)=o(n)$ when $\log K=o(n^{1/3})$. Finally, $\Prob(W_n)\to1$, because on the first event $y_{tj}=\tilde y_{tj}$, and $n^{-1}\sum_t\E y_{tj}^2\in[\underline\sigma^2,\bar\sigma^2]$.

For $w \in W_n$, Hoeffding's inequality for the $2K$ Rademacher sums $\pm n^{-1/2} \sum_t \varepsilon_t w_{tj}$ with $s_j^2 \leq 2\bar\sigma^2$ and the union bound give

\begin{equation}\label{eq:Choeff}
\Prob_\varepsilon\big(T(\varepsilon)>x\big)\le2K\exp\big(-x^2/(4\bar\sigma^2)\big),\qquad x>0.
\end{equation}

Note that when $x=h_n=2\bar\sigma\sqrt{\log(4Kn)}$, the right-hand side equals $1/(2n)$, so, according to \eqref{eq:Ccond}, the event $S_n=\{T\le h_n\}$ satisfies the condition $\Prob(S_n^c\cap W_n)\le1/(2n)$. And when $x=g_n=2\bar\sigma\sqrt{\log(4K/\gamma_n)}$, it equals $\gamma_n/2$, so fewer than $k_{\gamma_n}+1$ patterns have $T(\varepsilon)>g_n$, and $q^0_{\gamma_n}\le g_n$.

On $W_n\cap S_n$, Hoeffding's inequality shows that no more than $2^n\cdot2e^{-n\eta_n^2/2}=\beta_nN$ patterns satisfy $|\bareps|>\eta_n$, and $\bar e(\varepsilon)=\bar y(\varepsilon)-\bareps\,\bar y$ gives

\[ \sqrt n\,E(\varepsilon)\le T(\varepsilon)+|\bareps|\,T\le T(\varepsilon)+|\bareps|\,h_n. \]

Thus, every pattern with $\sqrt n\,E(\varepsilon)>q^0_{\gamma_n}+\eta_nh_n$ has $T(\varepsilon)>q^0_{\gamma_n}$ or $|\bareps|>\eta_n$, and there are at most $k_{\gamma_n}+\beta_nN\le(\alpha-\beta_n)N<k_{\alpha-\beta_n}+1$, that is, $\sqrt n\,q_{\alpha-\beta_n}\{E(\varepsilon)\}\le q^0_{\gamma_n}+\eta_nh_n$. From Lemma~\ref{lem:sandwich} with $\beta=\beta_n$, it follows that for $\eta_n\le1/2$

\[ \sqrt n\,c^*\le\frac{q^0_{\gamma_n}+\eta_nh_n}{1-\eta_n}\le q^0_{\gamma_n}+s_n,\qquad s_n=2\eta_n(g_n+h_n). \]

Lemma~J.4 of the supplement of \citet{CCKK2022} concerns vectors $Z_1,\dots,Z_n\in\R^p$, $p\ge2$, with $\sum_tZ_{tj}^2\ge b'n$ for all $j$ and $\|Z_t\|_\infty\le B'$ for all $t$, as well as independent Rademacher variables $\varepsilon_1,\dots,\varepsilon_n$. It states that $Y=n^{-1/2}\sum_t\varepsilon_tZ_t$ satisfies

\[ \Prob(Y\le y+r)-\Prob(Y\le y)\le C\,(r+B'/\sqrt n)\sqrt{\log p},\qquad y\in\R^p,\ r>0, \]

with coordinatewise inequalities and $C$ depending only on $b'$. For $w\in W_n$, let $Z_t=(w_t,-w_t)\in\R^{2K}$, then $b'=\underline\sigma^2/2$, $B'=c_n$, and, by the bijection $\varepsilon\mapsto(\varepsilon_t\sigma_t)_t$ above, $\Prob(Y\le x\mathbf 1)=\Prob_\varepsilon(T(\varepsilon)\le x)$. If we take $y=(x-u)\mathbf 1$ and $r=s+u$, and let $u\downarrow0$, we get

\[ \Prob_\varepsilon\big(x\le T(\varepsilon)\le x+s\big)\le C\,(s+c_n/\sqrt n)\sqrt{\log(2K)},\qquad x\in\R,\ s>0. \]

As at least $k_{\gamma_n}+1>\gamma_nN$ patterns satisfy $T(\varepsilon)\ge q^0_{\gamma_n}$,

\[ \Prob_\varepsilon\big(T(\varepsilon)>q^0_{\gamma_n}+s_n\big)\ge\gamma_n-C\,(s_n+c_n/\sqrt n)\sqrt{\log(2K)}. \]

By the bound on $c^*$ above, $\{T>q^0_{\gamma_n}+s_n\}\cap W_n\cap S_n\subseteq\{\|\bar y\|_\infty>c^*\}$. Then, conditioning on $w$ (both $q^0_{\gamma_n}$ and $W_n$ depend on $w$ only) and using \eqref{eq:Ccond}, we obtain

\[ \Prob\big(\|\bar y\|_\infty>c^*\big)\ge\Big(\gamma_n-C\,(s_n+c_n/\sqrt n)\sqrt{\log(2K)}\Big)\Prob(W_n)-\frac1{2n}. \]

Here $\gamma_n\to\alpha$, $\Prob(W_n)\to1$, and

\[ s_n\sqrt{\log(2K)}=O\Big(\log(Kn)\sqrt{\frac{\log n}n}\Big),\qquad\frac{c_n\sqrt{\log(2K)}}{\sqrt n}=O\Big(\frac{\log^{3/2}(Kn)}{\sqrt n}\Big), \]

because $\eta_n=O(\sqrt{\log n/n})$ and $g_n+h_n=O(\sqrt{\log(Kn)})$. Both tend to zero since $\log K=o(n^{1/3})$, so 

\[ \liminf_n\Prob(\|\bar y\|_\infty>c^*)\ge\alpha. \]

\end{proof}

\subsection*{Proof of Theorem~\ref{thm:effS}}

\begin{proof}[]

We write $D={\rm diag}(\sigma_1,\dots,\sigma_K)$, $z_t=D^{-1}y_t$, $S_j=\|y_{\cdot j}\|_2/\sqrt n$, and $\hat S_j=\|e_{\cdot j}\|_2/\sqrt n$, so that $n^{-1}\sum_t\E z_{tj}^2=1$ for each $j$. The statistic and the values obtained by resampling in $c_{\rm stud}$ are

\[ M_S=\max_j\frac{|\bar y_j|}{S_j},\qquad E_S(\varepsilon)=\max_j\frac{|\bar e_j(\varepsilon)|}{\hat S_j}, \]

and if we replace $S_j$ and $\hat S_j$ with $\sigma_j$, we get $\tilde M=\|\bar z\|_\infty$ and $E_z(\varepsilon)=\|D^{-1}\bar e(\varepsilon)\|_\infty$, that is, the function $E$ of the rows $z_t$ for $\phi=\|\cdot\|_\infty$. For $V:\mathcal E\to[0,\infty)$, let

\[ c(V)=\sup\big\{w\ge0:\ w\le q_\alpha\{V(\varepsilon)+|\bareps|\,w\}\big\}, \]

so that $c_{\rm stud}=c(E_S)$, and $\tilde c=c(E_z)$ is the threshold $c^*$ for rows $z_t$. Let $\tilde Q_\gamma=Q_\gamma(\tilde M)$ and

\[ \nu_n=\max_j\max\Big\{\Big|\frac{\hat S_j}{\sigma_j}-1\Big|,\ \Big|\frac{S_j}{\sigma_j}-1\Big|\Big\}. \]

From the fact that $q_\alpha\{\lambda v\}=\lambda q_\alpha\{v\}$ for $\lambda>0$, we have $c(\lambda V)=\lambda c(V)$, and \eqref{eq:qmono}, if $V\le V'$ pointwise, implies $c(V)\le c(V')$. On the set $\{\nu_n<1\}$, all $\hat S_j$ are positive, $E_z/(1+\nu_n)\le E_S\le E_z/(1-\nu_n)$, and the same inequalities relate $M_S$ to $\tilde M$, so

\begin{equation}\label{eq:Scomp}
\frac{\tilde c}{1+\nu_n}\le c_{\rm stud}\le\frac{\tilde c}{1-\nu_n},\qquad\frac{\tilde M}{1+\nu_n}\le M_S\le\frac{\tilde M}{1-\nu_n}\qquad\text{on }\{\nu_n<1\}.
\end{equation}

Let us fix a deterministic $ \bar\nu_n \in (0,1/3]$ and set $p_n = \Prob(\nu_n > \bar\nu_n)$. For $x \in \R$, from \eqref{eq:Scomp} we observe

\[ \Prob\big(\tilde M>(1+\bar\nu_n)x\big)-p_n\le\Prob(M_S>x)\le\Prob\big(\tilde M>(1-\bar\nu_n)x\big)+p_n, \]

so $\Prob(M_S>x)\le\gamma$ for $x=\tilde Q_{\gamma-p_n}/(1-\bar\nu_n)$ and $\Prob(M_S>x)>\gamma$ for $x<\tilde Q_{\gamma+p_n}/(1+\bar\nu_n)$. As in the proof of Lemma~\ref{lem:kolm}, for $p_n<\gamma<1-p_n$,

\begin{equation}\label{eq:SQ}
\frac{\tilde Q_{\gamma+p_n}}{1+\bar\nu_n}\le Q^{\rm stud}_\gamma\le\frac{\tilde Q_{\gamma-p_n}}{1-\bar\nu_n}.
\end{equation}

Let $1+\kappa_n^\pm=(1\pm\bar\nu_n)/(1\mp\bar\nu_n)$, so that $|\kappa_n^\pm|\le3\bar\nu_n$. If $\nu_n\le\bar\nu_n$, then from \eqref{eq:Scomp} it follows that

\begin{align*}
\tilde M>(1+\kappa_n^+)\tilde c\ &\Rightarrow\ M_S\ge\frac{\tilde M}{1+\bar\nu_n}>\frac{\tilde c}{1-\bar\nu_n}\ge c_{\rm stud},\\
M_S>c_{\rm stud}\ &\Rightarrow\ \tilde M\ge(1-\bar\nu_n)M_S>\frac{1-\bar\nu_n}{1+\bar\nu_n}\,\tilde c,
\end{align*}

and hence

\begin{equation}\label{eq:Slevel}
\Prob\big(\tilde M>(1+\kappa_n^+)\tilde c\big)-p_n\le\Prob(M_S>c_{\rm stud})\le\Prob\big(\tilde M>(1+\kappa_n^-)\tilde c\big)+p_n.
\end{equation}

By $|\sqrt x-1|\le|x-1|$ for $x\ge0$, $\hat S_j^2=S_j^2-\bar y_j^2$ and $S_j^2/\sigma_j^2=n^{-1}\sum_tz_{tj}^2$,

\begin{equation}\label{eq:Snu}
\nu_n\le\max_j\Big(\Big|\frac1n\sum_{t=1}^n(z_{tj}^2-\E z_{tj}^2)\Big|+\bar z_j^2\Big).
\end{equation}

The rows $z_t$ are i.i.d.\ and symmetric (with $\E\|z_1\|_2^4<\infty$ and positive definite covariance $D^{-1}\Sigma D^{-1}$), therefore Theorem~\ref{thm:effK} applies to them with $\phi=\|\cdot\|_\infty$: $\tilde c/\tilde Q_\alpha=1+O_P(\sqrt{\log n/n})$, and \eqref{eq:Klevel} holds for the rows $z_t$. Since $ \E z_{1j} ^2=1$, \eqref{eq:Snu}, the Chebyshev and Markov inequalities give

\[ p_n\le\sum_{j=1}^K\Big(\Prob\Big(\Big|\frac1n\sum_t(z_{tj}^2-1)\Big|>\frac{\bar\nu_n}2\Big)+\Prob\Big(\bar z_j^2>\frac{\bar\nu_n}2\Big)\Big)\le\sum_{j=1}^K\Big(\frac{4\E z_{1j}^4}{n\bar\nu_n^2}+\frac2{n\bar\nu_n}\Big)\le\frac{C_z}{n\bar\nu_n^2}, \]

where $C_z=\sum_j(4\E z_{1j}^4+2)$. Note that for $\bar\nu_n=Cn^{-1/2}$ and large $C$, this means that $\nu_n=O_P(n^{-1/2})$, and for quantiles we take $\bar\nu_n=n^{-1/3}$, which balances $\bar\nu_n$ and $p_n\le C_zn^{-1/3}$.

As in the proof of Theorem~\ref{thm:effK} for the rows $z_t$, take $L_0$, let $g_\gamma$ be the upper quantiles of $\|G\|_\infty$ with $G\sim N(0,D^{-1}\Sigma D^{-1})$, and let $\rho_n=O(n^{-1/2})$. Using \eqref{eq:Koracle} and \eqref{eq:Kgauss}, $\sqrt n\,\tilde Q_\alpha$ and $\sqrt n\,\tilde Q_{\alpha\pm p_n}$ differ from $g_\alpha\ge m>0$ by at most $L_0(p_n+\rho_n)$ for large $n$, and \eqref{eq:SQ} gives

\[ \tilde Q_\alpha/Q^{\rm stud}_\alpha=1+O(\bar\nu_n+p_n+\rho_n)=1+O(n^{-1/3}),\qquad\frac{c_{\rm stud}}{Q^{\rm stud}_\alpha}=\frac{c_{\rm stud}}{\tilde c}\cdot\frac{\tilde c}{\tilde Q_\alpha}\cdot\frac{\tilde Q_\alpha}{Q^{\rm stud}_\alpha} \]

Here the first factor is $1+O_P(n^{-1/2})$ from \eqref{eq:Scomp} and $\nu_n=O_P(n^{-1/2})$, the second is $1+O_P(\sqrt{\log n/n})$, and the third is $1+O(n^{-1/3})$, so the product is $1+O_P(n^{-1/3})$. Finally, substituting \eqref{eq:Klevel} for the rows $z_t$ with $\kappa_n=\kappa_n^\pm\to0$ into \eqref{eq:Slevel}, we observe

\[ \Prob(M_S>c_{\rm stud})\to\alpha. \]

\end{proof}

\subsection*{Proof of Proposition~\ref{prop:mc}}

\begin{proof}[]

We denote by $q_\gamma$ and $q^B_\gamma$ the resampled quantiles over $\{\pm1\}^n$ and over $\mathcal E_B$, where $k^B_\gamma=\lfloor\gamma(B+1)\rfloor$, and let $\varepsilon^1,\dots,\varepsilon^B$ be the random patterns, and let $\tilde r(\varepsilon)=+\infty$ if $|\bareps|=1$, and $\tilde r=r$ in all other cases. With the value $+\infty$ allowed in \eqref{eq:q}, the $(k_\gamma+1)$-th largest value of $\tilde r$ is $+\infty$ if $k_\gamma<m_1$, and the $(k_\gamma+1-m_1)$-th largest value of $r$ among patterns with $|\bareps|<1$ in all other cases. Thus, Lemma~\ref{lem:closed} reads

\[ c^*(\gamma;\mathcal E)=q_\gamma\{\tilde r\}\qquad\text{for every set }\mathcal E, \]

and \eqref{eq:qchar} holds for $\tilde r$, also when $w=+\infty$.

With the data fixed, let

\[ F(x)=\frac1{2^n}\#\{\varepsilon\in\{\pm1\}^n:\ \tilde r(\varepsilon)\le x\},\qquad F_B(x)=\frac1B\#\{b\le B:\ \tilde r(\varepsilon^b)\le x\}, \]

$D=\sup_{x\in\R}|F_B(x)-F(x)|$ and $\varsigma_B=u/\sqrt B$, so that $B\varsigma_B\ge2$. Given the data, the values $\tilde r(\varepsilon^b)$ are i.i.d.\ with a distribution function $F$ on $(-\infty,+\infty]$, and after an increasing mapping of $(-\infty,+\infty]$ onto a bounded interval of $\R$, they can be written as $F^{-1}(U_b)$, where $U_1,\dots,U_B$ are i.i.d.\ random variables with a uniform distribution on $(0,1)$. For the empirical distribution function $G_B$ of the random variables $U_b$, by Corollary~1 of \citet{Massart1990}, we observe

\[ \Prob(D>\varsigma_B)\le\Prob\Big(\sup_{v\in(0,1)}|G_B(v)-v|>\varsigma_B\Big)\le2e^{-2B\varsigma_B^2}=2e^{-2u^2}, \]

and on $\{D\le\varsigma_B\}$ also $|F_B(x-)-F(x-)|\le\varsigma_B$ for $x\in(-\infty,+\infty]$.

\emph{Lower bound.} Let $w=q_{\alpha+2\varsigma_B}\{\tilde r\}$, and then on $\{D\le\varsigma_B\}$, using \eqref{eq:qchar},

\begin{align*}
&1-F(w-)\ge\frac{k_{\alpha+2\varsigma_B}+1}{2^n}>\alpha+2\varsigma_B\\
\implies{}&\#\{b:\ \tilde r(\varepsilon^b)\ge w\}\ge B\big(1-F(w-)-\varsigma_B\big)>B\alpha+B\varsigma_B\ge B\alpha+2>k^B_\alpha+1,
\end{align*}

so $q^B_\alpha\{\tilde r\}\ge w$ by \eqref{eq:qchar}.

\emph{Upper bound.} In this step, let $w>q_{\alpha-2\varsigma_B}\{\tilde r\}$ be a real number, and then on the set $\{D\le\varsigma_B\}$, by combining \eqref{eq:qchar} with the identity pattern,

\begin{align*}
&1-F(w-)\le\frac{k_{\alpha-2\varsigma_B}}{2^n}\le\alpha-2\varsigma_B\\
\implies{}&\#\{\varepsilon\in\mathcal E_B:\ \tilde r(\varepsilon)\ge w\}\le B(\alpha-\varsigma_B)+1\le B\alpha-1<k^B_\alpha+1,
\end{align*}

so $q^B_\alpha\{\tilde r\}<w$ by \eqref{eq:qchar}. Since $w$ can be any value, $q^B_\alpha\{\tilde r\} \le q_{\alpha-2\varsigma_B}\{\tilde r\}$, which proves the first statement.

\emph{The second statement.} To prove this part, we fix $M \ge 1$ and take $u = M$, and then with a probability of at least $1 - 2e^{-2M^2}$ for large $n$,

\[ c^*\big(\alpha+2M/\sqrt{B_n};\{\pm1\}^n\big)\le c^*(\alpha;\mathcal E_{B_n})\le c^*\big(\alpha-2M/\sqrt{B_n};\{\pm1\}^n\big), \]

and \eqref{eq:Kresamp} at these two levels gives, as in the derivation of \eqref{eq:Kcstar},

\[ g_{\alpha+\Delta_n+1/n+2M/\sqrt{B_n}}-L\hat\kappa_ns_n\le\sqrt n\,c^*(\alpha;\mathcal E_{B_n})\le\frac{g_{\alpha-\theta_n-2M/\sqrt{B_n}}+L\hat\kappa_ns_n}{1-\eta_n}. \]

This is \eqref{eq:Kcstar} with all levels shifted by $2M/\sqrt{B_n}$, so the remaining part of that proof yields $c^*/Q^*_\alpha=1+O_P(\sqrt{\log n/n}+B_n^{-1/2})$, since $M$ can be chosen so that $2e^{-2M^2}$ is arbitrarily small. For the level, let $v_n\to\infty$, where $v_n=o(\sqrt{B_n})$, and $b_n=n^{-1/4}+2v_n/\sqrt{B_n}$, and take $u=v_n$. On the intersection of $\Omega'_n$ with the event described in the first statement, the probability of which is no less than $\Prob(\Omega'_n)-2e^{-2v_n^2}\to1$, the bounds $\ell_n$ and $u_n$ become

\[ g_{\alpha+b_n+1/n}-b_n\le\sqrt n\,c^*(\alpha;\mathcal E_{B_n})\le\frac{g_{\alpha-\beta_n-2^{-n}-b_n-1/n}+b_n}{1-\eta_n}, \]

and both bounds tend to $g_\alpha$, so the rest of the argument regarding the level does not change.

\end{proof}

\section{Proofs for Section~\ref{sec:stepdown}}\label{app:mt}

Recall that the recursion \eqref{eq:sdrec}, where $C_{i+1}=\emptyset$ if $C_i=\emptyset$, is the step-down procedure described in \citet[Definition~3.1]{ABR2010b}. The sets in it decrease and become constant starting from some $i\le K$. We denote by $C_\infty$ the final set, and by $R=\{1,\dots,K\}\setminus C_\infty$ the rejection set. If $I_0\ne\emptyset$ and the function $t$ is nondecreasing, i.e., $t(C)\le t(C')$ when $C\subseteq C'$, then \citet[Thm.~3.2]{ABR2010b}, who reformulate \citet[Thm.~3]{RomanoWolf2005} for \eqref{eq:sdrec}, state that

\begin{equation}\label{eq:RW}
\Prob(R\cap I_0\ne\emptyset)\le\Prob\Big(\max_{j\in I_0}T_j>t(I_0)\Big).
\end{equation}

Their proof shows that no true hypothesis is rejected on the event $\{\max_{j\in I_0}T_j\le t(I_0)\}$ for every realization, so \eqref{eq:RW} also holds in the case where the thresholds depend on the random patterns in $\mathcal E$.

All thresholds in Section~\ref{sec:stepdown} are fixed points of the form \eqref{eq:cstar} or quantiles, and their monotonicity with respect to $C$ follows from \eqref{eq:qmono}. For a function $V:\mathcal E\to(-\infty,+\infty]$, a level $\gamma$, and a function $a:\mathcal E\to[0,\infty)$, we define, taking into account the conventions $0\cdot\infty=0$ and $\sup\emptyset=0$,

\[ \mathcal W(V)=\big\{w\ge0:\ w\le q_\gamma\{V(\varepsilon)+a(\varepsilon)\,w\}\big\}. \]

If $V\le V'$ pointwise, then $V+aw\le V'+aw$ for any $w$, and from \eqref{eq:qmono} it follows that $\mathcal W(V)\subseteq\mathcal W(V')$, and hence

\begin{equation}\label{eq:Wmono}
V\le V'\ \text{pointwise}\quad\Longrightarrow\quad\sup\mathcal W(V)\le\sup\mathcal W(V').
\end{equation}

\begin{lemma}\label{lem:dom}
Consider two step-down procedures with identical statistics $T_j$ and thresholds $t_1$ and $t_2$. Let $C^{(1)}_i$, $C^{(2)}_i$ and $R_1$, $R_2$ be their sets and rejection sets. If $t_1$ is nondecreasing, that is, $t_1(C)\le t_1(C')$ for $C\subseteq C'$, and $t_1(C)\le t_2(C)$ for any nonempty set $C$, then $R_2 \subseteq R_1$.
\end{lemma}

\begin{proof}

We show that $C^{(1)}_i\subseteq C^{(2)}_i$ by induction on $i$, where for $i=0$ both sets are equal to $\{1,\dots,K\}$. If $C^{(1)}_i\subseteq C^{(2)}_i$, then using \eqref{eq:sdrec} and the two assumptions about $t_1$, we observe

\[ j\in C^{(1)}_{i+1}\implies j\in C^{(2)}_i\ \text{and}\ T_j\le t_1(C^{(1)}_i)\le t_1(C^{(2)}_i)\le t_2(C^{(2)}_i)\implies j\in C^{(2)}_{i+1}. \]

Hence, $C^{(1)}_\infty\subseteq C^{(2)}_\infty$, and if we take the complements, then $R_2\subseteq R_1$.

\end{proof}

\subsection*{Proof of Theorem~\ref{thm:fwer}}

\begin{proof}[]

Without loss of generality, we may assume that $I_0 \neq \emptyset$, since otherwise there is nothing to prove. From \eqref{eq:RW}, it follows that it is sufficient to show that each threshold does not decrease as $C$ increases and that

\begin{equation}\label{eq:fwerI0}
\Prob\Big(\max_{j\in I_0}T_j>t(I_0)\Big)\le\alpha.
\end{equation}

The rows $(y_{tj})_{j\in I_0}$ satisfy Assumption~\ref{ass:SA}, since changing the sign in a row changes the sign of each of its coordinates, so Theorems~\ref{thm:cstar} and \ref{thm:variants}, as well as Lemma~\ref{lem:exact}, are applicable to the observations $(Y_{tj})_{j\in I_0}$ with the mean $(\mu_j)_{j\in I_0}$ and the residuals $(e_{tj})_{j\in I_0}$.

\emph{(a).} Here, $t(C)=\sup\mathcal W(V_C)$ with $\gamma=\alpha$, $a=|\bareps|$, and $V_C(\varepsilon)=\max_{j\in C}|\bar e_j(\varepsilon)|$, which does not decrease with $C$, so $t$ does not decrease by \eqref{eq:Wmono}. Under $H_j$, $j\in I_0$, we have $\bar Y_j=\bar y_j$, and Theorem~\ref{thm:cstar} with $\phi(v)=\max_{j\in I_0}|v_j|$ gives

\[ \Prob\Big(\max_{j\in I_0}|\bar Y_j|>t(I_0)\Big)=\Prob\Big(\max_{j\in I_0}|\bar y_j|>c^*(I_0)\Big)\le\alpha. \]

\emph{(b).} Now, $t(C)=\sup\mathcal W(V_C)$ with $\gamma=\alpha$, $a=|\bareps|$, and $V_C(\varepsilon)=\max_{j\in C}\sqrt n\,|\bar e_j(\varepsilon)|/\|e_{\cdot j}\|_2$, with the convention $0/0=0$ of Appendix~\ref{app:cstar}, and the function $t$ is non-decreasing (using \eqref{eq:Wmono}). Under $H_j$, $j\in I_0$, we have $Y_{\cdot j}=y_{\cdot j}$, and it follows from Theorem~\ref{thm:variants} that

\[ \Prob\Big(\max_{j\in I_0}\frac{\sqrt n\,|\bar y_j|}{\|y_{\cdot j}\|_2}>t(I_0)\Big)\le\alpha. \]

\emph{(c).} In this part, $t(C)=q_\alpha\{U_C(\varepsilon)\}$ with $U_C(\varepsilon)=\max_{j\in C}|n^{-1}\sum_t\varepsilon_tY_{tj}|$ in the two-sided case and $U_C(\varepsilon)=\max_{j\in C}n^{-1}\sum_t\varepsilon_tY_{tj}$ in the one-sided case. In both cases, $U_C \le U_{C'}$ when $C \subseteq C'$, and the function $t$ is non-decreasing according to \eqref{eq:qmono}, as in \citet[(23)--(24)]{RomanoWolf2005}. In the two-sided case, with $W(z)=\max_{j\in I_0}|\bar z_j|$

\[ \mu_j=0\ \ (j\in I_0)\implies U_{I_0}(\varepsilon)=W(\varepsilon\cdot y),\quad\max_{j\in I_0}T_j=W(y), \]

and Lemma~\ref{lem:exact} gives

\[ \Prob\Big(\max_{j\in I_0}|\bar Y_j|>t(I_0)\Big)=\Prob\big(W(y)>q_\alpha\{W(\varepsilon\cdot y)\}\big)\le\alpha. \]

In the one-sided case, with $W(z)=\max_{j\in I_0}(\bar z_j+\mu_j)$, which does not depend on the data,

\[ \mu_j\le0\ \ (j\in I_0),\ \ \bareps\le1\implies\frac1n\sum_{t=1}^n\varepsilon_tY_{tj}=\bar y_j(\varepsilon)+\bareps\,\mu_j\ge\bar y_j(\varepsilon)+\mu_j\implies U_{I_0}(\varepsilon)\ge W(\varepsilon\cdot y), \]

while $\max_{j\in I_0}\bar Y_j=W(y)$, so from \eqref{eq:qmono} and Lemma~\ref{lem:exact} we observe

\[ \Prob\Big(\max_{j\in I_0}\bar Y_j>t(I_0)\Big)\le\Prob\big(W(y)>q_\alpha\{W(\varepsilon\cdot y)\}\big)\le\alpha. \]

Thus, in each of parts (a)--(c), the threshold does not decrease as $C$ increases, and \eqref{eq:fwerI0} holds, so it follows from \eqref{eq:RW} that $\Prob(R\cap I_0\ne\emptyset)\le\alpha$.

\end{proof}

\subsection*{Proof of Theorem~\ref{thm:poweror}}

\begin{proof}[]

To prove this theorem, let $c^*(C)$ be the threshold from Theorem~\ref{thm:fwer}(a), that is, $c^*$ for $\phi_C(v)=\max_{j\in C}|v_j|$, and let $Q^*(C)$ be the upper $\alpha$-quantile of $\phi_C(\bar y)$. Observe that each $\phi_C$ is a seminorm that is not identically equal to zero, so Theorem~\ref{thm:effK} applies to each of the $2^K-1$ nonempty sets $C$. And for the event

\[ A_n=\big\{(1-u)\,Q^*(C)\le c^*(C)\le(1+u)\,Q^*(C)\ \text{for every nonempty }C\big\} \]

the union bound gives

\[ \Prob(A_n^c)\le\sum_{C\ne\emptyset}\Prob\Big(\Big|\frac{c^*(C)}{Q^*(C)}-1\Big|>u\Big)\to0. \]

The threshold $c^*(C)$ does not decrease with $C$, according to the proof of Theorem~\ref{thm:fwer}(a), and $Q^*(C)$ does not decrease, since $\phi_C(\bar y)\le\phi_{C'}(\bar y)$ for $C\subseteq C'$. On $A_n$, Lemma~\ref{lem:dom} gives

\begin{align*}
(t_1,t_2)&=\big(c^*,(1+u)Q^*\big)\implies R_{\rm or}(1+u)\subseteq R_{c^*},\\
(t_1,t_2)&=\big((1-u)Q^*,c^*\big)\implies R_{c^*}\subseteq R_{\rm or}(1-u),
\end{align*}

and therefore

\[ \Prob\big(R_{\rm or}(1+u)\subseteq R_{c^*}\subseteq R_{\rm or}(1-u)\big)\ge\Prob(A_n)\to1. \]

\end{proof}

\subsection*{Proof of Theorem~\ref{thm:bonf}}

\begin{proof}[]

Let $C$ be nonempty, $\gamma=\alpha/|C|$, and $w>t_B(C)$. Then $w>c^*_j(\gamma)$, so $w$ does not belong to the set from \eqref{eq:cstar} for the $j$-th column at level $\gamma$, i.e.,

\[ g_j:=q_\gamma\big\{|\bar e_j(\varepsilon)|+|\bareps|\,w\big\}<w,\qquad j\in C. \]

The infimum is achieved in \eqref{eq:q} because the function $x\mapsto\#\{\varepsilon:\ v(\varepsilon)>x\}$ is nonincreasing and right-continuous, so no more than $k_\gamma$ patterns satisfy $|\bar e_j(\varepsilon)|+|\bareps|w>g_j$. With $g = \max_{j \in C} g_j < w$ and $|C| \lfloor \gamma N \rfloor \le \lfloor \alpha N \rfloor$,

\begin{align*}
\#\Big\{\varepsilon:\ \max_{j\in C}|\bar e_j(\varepsilon)|+|\bareps|\,w>g\Big\}&\le\sum_{j\in C}\#\big\{\varepsilon:\ |\bar e_j(\varepsilon)|+|\bareps|\,w>g_j\big\}\le|C|\,k_\gamma\le k_\alpha\\
&\implies q_\alpha\Big\{\max_{j\in C}|\bar e_j(\varepsilon)|+|\bareps|\,w\Big\}\le g<w
\end{align*}

by \eqref{eq:q}, so no $w>t_B(C)$ belongs to the set from \eqref{eq:cstar} for $\phi_C$, and $c^*(C)\le t_B(C)$. Since $c^*(C)$ is non-decreasing in $C$, according to the proof of Theorem~\ref{thm:fwer}(a), from Lemma~\ref{lem:dom} with $t_1=c^*$ and $t_2=t_B$, it follows that $R_B\subseteq R_{c^*}$.

\end{proof}

\end{document}